\documentclass[10pt]{article}

\usepackage{amsmath}
\usepackage{amsthm}
\usepackage{amssymb}
\usepackage{amsfonts}
\usepackage{caption}
\usepackage{algorithm}  
\usepackage{algorithmicx} 
\usepackage{algpseudocode}
\floatstyle{ruled}
\newfloat{procedure}{htbp}{lop}
\floatname{procedure}{Procedure}  
\usepackage[numbers]{natbib}
\usepackage{booktabs}
\usepackage{multirow}
\usepackage{color,xcolor}

\newif\ifreview
\reviewfalse
\ifreview
  \usepackage[mathlines]{lineno}
  \linenumbers
\fi

\usepackage[colorlinks]{hyperref}
\usepackage{graphicx}     %

\usepackage[a4paper, top=2.5cm, bottom=2.5cm, left=2cm, right=2cm]{geometry}

\ifreview
  \newcommand*\patchAmsMathEnvironmentForLineno[1]{%
    \expandafter\let\csname old#1\expandafter\endcsname\csname #1\endcsname
    \expandafter\let\csname oldend#1\expandafter\endcsname\csname end#1\endcsname
    \renewenvironment{#1}%
      {\linenomath\csname old#1\endcsname}%
      {\csname oldend#1\endcsname\endlinenomath}}
  \newcommand*\patchBothAmsMathEnvironmentsForLineno[1]{%
    \patchAmsMathEnvironmentForLineno{#1}%
    \patchAmsMathEnvironmentForLineno{#1*}}
  \AtBeginDocument{%
    \patchBothAmsMathEnvironmentsForLineno{equation}%
    \patchBothAmsMathEnvironmentsForLineno{align}%
    \patchBothAmsMathEnvironmentsForLineno{flalign}%
    \patchBothAmsMathEnvironmentsForLineno{alignat}%
    \patchBothAmsMathEnvironmentsForLineno{gather}%
    \patchBothAmsMathEnvironmentsForLineno{multline}}
\fi

\makeatletter
\@addtoreset{equation}{section}
\makeatother

\newtheorem{theorem}{Theorem}[section]
\newtheorem{lemma}{Lemma}[section]
\newtheorem{pro}{Proposition}[section]
\newtheorem{assu}{Assumption}[section]
\newtheorem{remark}{Remark}[section]
\newtheorem{definition}{Definition}[section]
\newtheorem{corollary}{Corollary}[section]
\newcommand{\inp}[1]{\left\langle#1\right\rangle}
\newcommand{\enp}[1]{\langle#1\rangle}
\newcommand{\norm}[1]{\left\Vert#1\right\Vert}
\newcommand{\enorm}[1]{\Vert#1\Vert}

\def\bbr{{\mathbb{R}}} %
\def\bbe{{\mathbb{E}}} %

\global\long\def\norm#1{\lVert#1\rVert}%
\global\long\def\abs#1{\left |#1\right |}%
\global\long\def\bsbra#1{\big[#1\big]}%

\global\long\def\vertiii#1{\left\vert \kern-0.25ex  \left\vert \kern-0.25ex  \left\vert #1\right\vert \kern-0.25ex  \right\vert \kern-0.25ex  \right\vert }%

\global\long\def\and{\mathrm{and}}%

\global\long\def\st{\mathrm{s.t.}}%

\global\long\def\st{\operatornamewithlimits{s.t.}}%

\newif\ifcomments
\commentsfalse
\ifcomments
  \newcommand{\comm}[2][]{\textcolor{green!50!black}{\textbf{[#1:} #2\textbf{]}}}
  \newcommand{\qd}[1]{\textcolor{red}{\textbf{[QD:} #1\textbf{]}}}
  \newcommand{\bz}[1]{\textcolor{blue}{\textbf{[BZ:} #1\textbf{]}}}
\else
  \newcommand{\comm}[2][]{}
  \newcommand{\qd}[1]{}
  \newcommand{\bz}[1]{}
\fi

\title{{\bf An Accelerated Dynamic Stochastic Primal--Dual Method for Multi-stage Stochastic Optimization}}
\author{Bin Zhou\footnote{
		School of Science, Nanjing University of Posts and Telecommunications/Ministry of Education Key Laboratory of NSLSCS,
		Nanjing, 210023, China, \texttt{415281@njupt.edu.cn}}   
\qquad 
Zhiyun Guo  \footnote{Antai College of Economics and Management, Shanghai Jiao Tong University, \texttt{gzy7777@sjtu.edu.cn}}
\qquad
	Hailin Sun\footnote{
		School of Mathematical Sciences, Jiangsu Key Laboratory for NSLSCS, Nanjing Normal University, Nanjing 210023, China, \texttt{hlsun@njnu.edu.cn}}
		\qquad 
	Qi Deng\footnote{Corresponding author,  Antai College of Economics and Management,
		Shanghai Jiao Tong University
		\texttt{qdeng24@sjtu.edu.cn}}
}

\begin{document}
\maketitle

\begin{abstract}

We consider multi-stage stochastic optimization problems with smooth convex stage costs and linear linking constraints.
Dynamic stochastic approximation (DSA), a recently proposed stochastic first-order method for such problems, uses a recursive structure whose memory requirement grows only linearly with the number of stages.
Its implementation, however, requires an exact solution of the proximal subproblem associated with the stage cost, and this subproblem often admits no closed-form solution. To resolve this issue, we develop an accelerated dynamic stochastic primal-dual (A-DSPD) framework that retains the recursive value-function structure and only requires the gradient of the stage cost.
Within this framework, each stagewise saddle-point subproblem is solved by an inexact accelerated stochastic primal-dual method (iASPD), which applies Nesterov acceleration to the smooth stage cost and uses recursively generated stochastic $\epsilon$-subgradients for the nonsmooth downstream value function.
For the stochastic composite saddle-point problem, after $N$ primal--dual updates, iASPD yields optimality and feasibility bounds in which the contribution of the smooth stage cost decays at a faster rate $O(1/N^2)$, while the terms associated with stochastic $\epsilon$-subgradients of the nonsmooth value function decay as $O(1/\sqrt{N})$, up to the prescribed oracle inexactness.
We also develop bounds for the stochastic $\epsilon$-subgradient estimators and the dual iterates, which yield complexity guarantees for A-DSPD.
Specifically, A-DSPD computes an $\epsilon$-solution using $O(1/\epsilon^4)$ sampled scenarios and $O(1/\epsilon^5)$ primal--dual updates for three-stage problems. We further extend the method and its convergence analysis to problems with an arbitrary finite number of stages.
Numerical experiments on three-stage stochastic portfolio models illustrate the empirical advantage of the proposed method.

\end{abstract}


\newcommand{\iaspd}{iASPD}
\newcommand{\ADSPD}{A-DSPD}
\newcommand{\DSA}{DSA}

\section{Introduction}
We consider the following multi-stage stochastic optimization problem:
\begin{equation}\label{multi}
\begin{array}{lll}
\underset{x_1\in X_1}{\min} f_1\left(x_1, c_1\right)+&\bbe_{\xi_2\mid\xi_{[1]}}\left[\underset{x_2\in X_2}{\min} f_2\left(x_2, c_2\right)\right.+\cdots & \left.\bbe_{\xi_T\mid\xi_{[T-1]}}\left[\underset{x_T\in X_T}{\min} f_T\left(x_T, c_T\right)\right]\right] \\
\text { s.t. } A_1 x_1- b_1\leq 0  & \text { s.t. } A_2 x_2+B_2 x_1-b_2\leq 0 , & \text { s.t. } A_T x_T+B_T x_{T-1}-b_T \leq 0
\end{array}
\end{equation}
where $T$ denotes the number of stages; the feasible sets $X_t \subset \bbr^{n_t}$ are compact and convex; for each cost parameter $c_t$, the functions $f_t(\cdot,c_t):X_t\to\mathbb{R}$ are smooth and convex; and $A_t\in \bbr^{m_t\times n_t}$ and $b_t\in \bbr^{m_t}$ are stage data for $t=1,\ldots,T$, while the linking matrix $B_t\in \bbr^{m_t\times n_{t-1}}$ is defined only for $t=2,\ldots,T$. Moreover, $\xi_1=(A_1,b_1,c_1)$ is deterministic, while $\xi_t=(A_t,B_t,b_t,c_t)$ for $t=2,\cdots,T$ are random vectors supported on $\Xi_t$ at stage $t$. We use $\xi_{[t]}=(\xi_1,\xi_2,\cdots,\xi_t)$
to denote the scenario history through stage $t\geq 1$, with support $\Xi_{[t]}$. In particular, $\xi_{[1]}=\xi_1$; conditioning on this deterministic first-stage history leaves expectations over subsequent randomness unchanged. Problem~\eqref{multi} thus covers a broad class of sequential decision problems under uncertainty. The study of multi-stage stochastic programming was initiated by \citet{dantzig1955linear} and has since found applications in operational planning, finance, energy management, and natural-resource allocation~\citep{kouwenberg2001scenario, pereira1991multi, homem2016risk, tilmant2007stochastic, flach2010long}. Among these applications, multi-stage portfolio selection balances investment risk against expected return subject to budget, no-short-sale, and inter-stage turnover constraints. We study this application numerically in Section~\ref{sec:numerical}.

The classical solution framework for such multi-stage problems is dynamic programming~\citep{bellman1966dynamic}, which delivers exact solutions but can become computationally prohibitive as the state and decision spaces grow.
A range of approximation and decomposition techniques reduces this cost.
An early line of work replaces the expectation by a sample average approximation (SAA)~\citep{kleywegt2002sample}, which gives a deterministic surrogate that standard optimization methods can solve.
Decomposition algorithms such as the L-shaped method~\citep{van1969shaped} and stochastic dual dynamic programming (SDDP)~\citep{pereira1991multi} exploit problem structure by iteratively constructing cutting-plane approximations of cost-to-go functions, while progressive hedging~\citep{rockafellar1991scenarios} uses scenario-based decomposition with penalty coordination and is therefore well suited to parallel computation.

Recently, \citet{lan2021dynamic} introduced the {dynamic stochastic approximation} (\DSA{}) method, which extends stochastic approximation techniques~\citep{polyak1990new,nemirovskii1983problem} from single-stage to multi-stage settings.
In \DSA{}, each stage is formulated as a constrained optimization problem and solved approximately by a primal--dual algorithm, where biased subgradient estimates of the value (cost-to-go) function are generated by approximately solving the next-stage problem.
For three-stage problems, \DSA{} attains a sample complexity of $\mathcal{O}(1/\epsilon^4)$ for general convex objectives and $\mathcal{O}(1/\epsilon^2)$ under strong convexity.
Its sampling bounds do not scale with the size of an explicitly enumerated scenario tree. 
However, the method relies on efficient evaluation of the proximal mapping associated with $f_t(x_t, c_t)$, an assumption that is restrictive in many practical models.
For instance, when $f_t$ is quadratic, evaluating the proximal operator typically requires solving a quadratic program (or an equivalent second-order cone program), often through an external solver, which can substantially increase the per-iteration cost.

To address this limitation of \DSA{}, we develop an accelerated dynamic stochastic primal-dual (\ADSPD) algorithm for the multi-stage stochastic optimization problem~\eqref{multi} (MSOP).
Following the \DSA{} framework, we reformulate MSOP as a sequence of single-stage saddle-point problems (SPPs) parameterized by value functions from subsequent stages.
Whereas \DSA{} requires efficient evaluation of the proximal mapping associated with $f_t$, \ADSPD{} assumes only access to $\nabla f_t$, which is often available in closed form even when the proximal operator is not, for example when $f_t$ is quadratic.
At each stage, the corresponding SPP is solved by an inexact accelerated stochastic primal-dual method (\iaspd), which is a variant of the accelerated primal-dual (APD) method~\citep{chen2014optimal,liu2023single}.
For any stage $t<T$, the \iaspd{} update requires an approximate $\epsilon$-subgradient of the next-stage value function, and this quantity is obtained recursively by applying the same \iaspd{} procedure to the subsequent stage.
Computing an $\epsilon$-solution to MSOP still requires repeated iterations, but each iteration traverses a single sampled root-to-leaf scenario path instead of the full scenario tree, and the memory requirement grows linearly with the number of stages.

Our three-stage analysis proceeds backward: we first control the third-stage stochastic $\epsilon$-subgradient and dual iterate, and then propagate these bounds through the second stage to establish first-stage optimality and feasibility. Under the assumptions of Theorem~\ref{th:ADSPD-three}, A-DSPD computes an $\epsilon$-solution to a three-stage problem using $\mathcal{O}(1/\epsilon^4)$ sampled scenarios and $\mathcal{O}(1/\epsilon^5)$ primal--dual updates. Table~\ref{Table:complexity-new} summarizes the resulting bounds, while Corollaries~\ref{cor:last-stage}--\ref{cor:first-stage} give the formal statements.
\begin{itemize}
		\item \textbf{Third stage:} The terminal problem is deterministic and has no downstream stochastic or oracle error. We therefore use the choice in Corollary~\ref{cor:last-stage}, which gives sharper subgradient accuracy and tighter expected dual control than the settings required at earlier stages.
	
		\item \textbf{Second stage:} Intermediate stages must absorb stochastic and inexact downstream information while passing a valid subgradient backward. We use the setting in Corollary~\ref{cor:middle-stage}, which ensures the boundedness of the dual solution and supports recursive subgradient construction from stage $2$ through stage $T-1$.
	
		\item \textbf{First stage:} No subgradient must be passed further backward, so the choice in Corollary~\ref{cor:first-stage} directly strengthens the optimality and feasibility bounds. It improves the multiplier-dependent term from $\frac{\norm{A}(2+4\left\|y^*\right\|^2)}{\sqrt{N}}$ to $\frac{\norm{A}(2+4\left\|y^*\right\|^2)}{N}$, while the remaining terms keep the same order.
\end{itemize}

\begin{table}[htbp]
	\caption{Convergence rates of the A-DSPD algorithm in three-stage stochastic programs. Here, $N$ denotes the number of primal--dual updates in a single iASPD call for the corresponding stage subproblem. For simplicity, we omit the stage subscript $t$.}
	\label{Table:complexity-new}
	\centering
	\resizebox{\textwidth}{!}{
		\begin{tabular}{ccc}
			\toprule
			{\bf Stage} & {\bf Metric} & {\bf A-DSPD Convergence Rates} \\
			\midrule
			\multirow{2}{*}{Third}  & $\bbe[\operatorname{gap}(\tilde{z}^{N+1})]$ & $\tfrac{2\Omega^2_X(L_f+1)}{N^2}+\tfrac{\enorm{A}(2\Omega_X^2+\|y^*\|^2)}{N}+\frac{4M^2(2N+1)}{3}+\epsilon$ \\[0.3ex]
			(Cor.~\ref{cor:last-stage}) & $\bbe[\enorm{\tilde{y}^{N+1}}^2]$ & $8\Omega_X^2+\frac{8\Omega_X^2L_f}{N\enorm{A}}
			+\frac{8\Omega_X^2}{N\enorm{A}}+6\enorm{y^*}^2+\frac{16M^2(N+1)(2N+1)}{3\enorm{A}}+\frac{4(N+1)\epsilon}{\enorm{A}}$\\
            [0.5ex]
			\midrule
				\multirow{4}{*}{Second } & $\bbe[\operatorname{gap}(\tilde{z}^{N+1})]$ &$\tfrac{2\Omega^2_XL_f}{N^2}+\tfrac{2M}{3N^{3/2}}+\tfrac{2\Omega_X^2\|A\|}{N}+\tfrac{2M\Omega^2_X+4M+\enorm{A}\enorm{y^*}^2}{N^{1/2}}+\epsilon$ \\[0.5ex]
				& $\bbe[\enorm{\tilde{y}^{N+1}}^2]$ & $\frac{8L_f\Omega_X^2}{N^{3/2}\enorm{A}}+\frac{8\Omega_X^2}{\sqrt{N}}+\frac{8M\Omega_X^2+16M}{\enorm{A}}+6\enorm{y^*}^2
			+\frac{64M}{3\enorm{A}}+\frac{8\sqrt{N}\epsilon}{\enorm{A}}$\\[0.8ex]
            (Cor.~\ref{cor:middle-stage}) & $\bbe[|f(\tilde{x}^{N+1},c)+v(\tilde{x}^{N+1})-f(x^*,c)-v(x^*)|]$ &
			\multirow{2}{*}{
				$\tfrac{2\Omega^2_XL_f}{N^2}+\tfrac{2M}{3N^{3/2}}+\tfrac{\enorm{A}(2\Omega^2_X+2+4\enorm{y^*}^2)}{N^{1/2}}+\tfrac{2M\Omega^2_X+4M}{N^{1/2}}+\epsilon$
			} \\[0.8ex]
			   & $\bbe[\enorm{[A\tilde{x}^{N+1}-b+Bu]_+}]$ \\
            [0.5ex]
			\midrule
			\multirow{2}{*}{First} & $\bbe[|f(\tilde{x}^{N+1},c)+v(\tilde{x}^{N+1})-f(x^*,c)-v(x^*)|]$ &
			\multirow{2}{*}{
				$\tfrac{2\Omega^2_XL_f}{N^2}+\tfrac{2M}{3N^{3/2}}+\tfrac{\enorm{A}(2\Omega^2_X+2+4\enorm{y^*}^2)}{N}+\tfrac{2M\Omega^2_X+4M}{N^{1/2}}+\epsilon$
			} \\[0.8ex]
			(Cor.~\ref{cor:first-stage}) & $\bbe[\enorm{[A\tilde{x}^{N+1}-b+Bu]_+}]$ \\
			[0.3ex]
			\bottomrule
	\end{tabular}}
\end{table}

A central element of the A-DSPD framework is the \iaspd{} method, designed to address stochastic saddle-point problems and inspired by the accelerated primal-dual (APD) algorithm~\citep{chen2014optimal}.
While APD achieves optimal complexity for a wide range of saddle-point problems, it cannot be directly employed in our context, as it relies on a differentiable stochastic component and the availability of unbiased stochastic gradients.
In contrast, the stochastic component in our formulation arises from a value function, which is generally nonsmooth and for which an exact subgradient oracle is typically unavailable.
To address this mismatch, we develop \iaspd{} for constrained stochastic problems whose objective consists of a deterministic smooth component and a nonsmooth stochastic component. Its finite-sample optimality and feasibility bounds separate the effects of these two components: the contribution of the smooth term decays at the accelerated rate $\mathcal{O}(1/N^2)$, while the terms associated with stochastic $\epsilon$-subgradients of the nonsmooth term decay as $\mathcal{O}(1/\sqrt{N})$, up to the prescribed oracle inexactness.
Compared with the inexact primal--dual stochastic approximation (I-PDSA) method of \citet{lan2021dynamic} for the single-stage subproblems arising in MSOP, \iaspd{} retains the overall $\mathcal{O}(1/\sqrt{N})$ rate while replacing the proximal mapping of $f_t$ with a gradient evaluation.

\paragraph{Outline of the paper.}
The rest of this paper is organized as follows. In Section~\ref{sec:three-stage}, we provide the framework of the A-DSPD algorithm by integrating a multi-step acceleration scheme and establishing some convergence properties for solving a three-stage stochastic optimization problem. In Section~\ref{sec:multi-stage}, we develop the A-DSPD
algorithm for solving more general multi-stage stochastic optimization problems. 
In Section~\ref{sec:numerical}, we report numerical studies on two three-stage stochastic portfolio models.
Finally, we conclude the paper in Section~\ref{sec:conclusion}.

\subsection{Notation and Preliminaries}

Let $\mathbb{R}^{n}$ denote the $n$-dimensional Euclidean space. For vectors $x, y \in \mathbb{R}^{n}$, we write $\enp{x,y}$ for the standard inner product and $\enorm{x} := \sqrt{\enp{x,x}}$ for the Euclidean norm. Unless stated otherwise, all norms are Euclidean.
For a scalar $a\in\mathbb{R}$, define $[a]_+:=\max\{a,0\}$; for vectors, this operation and all inequalities are interpreted componentwise. We write $\mathbb{R}_+^m$ for the nonnegative orthant, $\mathbf{1}$ for the all-ones vector of the required dimension, $\|x\|_1:=\sum_i|x_i|$, and $\operatorname{ri}(X)$ for the relative interior of a convex set $X$. When a function $f$ is differentiable at $x$, we denote its gradient by $\nabla f(x)$. Expectation with respect to all scenario and oracle randomness is denoted by $\bbe[\cdot]$.
For the stagewise recursion, $\mathbb{E}_{\xi_t\mid\xi_{[t-1]}}[\cdot]$ denotes expectation with respect to the stage-$t$ scenario conditional on $\xi_{[t-1]}$, whereas $\mathbb{E}[\cdot\mid\xi_{[t]}]$ averages over all remaining randomness conditional on the realized scenario history. 
For the complexity analysis in each stage, $\mathcal{F}_{k-1}$ contains the current problem data (including the realized scenario history in a stagewise subproblem), the initialization, and all oracle observations preceding the $k$-th oracle call. We write $\bbe_k[\cdot]:=\bbe[\cdot\mid\mathcal{F}_{k-1}]$; in particular, the iterates available at the beginning of iteration $k$ are $\mathcal{F}_{k-1}$-measurable.
For a nonempty bounded set $X\subset\mathbb{R}^n$, let $\Omega_X = \sup_{x,x'\in X}\enorm{x-x'}$. For stagewise sets, we write $\Omega_t:=\Omega_{X_t}$.
For a matrix $A$, $\enorm{A}$ denotes the Euclidean operator norm.
A function $f: \mathbb{R}^n \to \mathbb{R}$ is said to be $L$-smooth if it is differentiable and satisfies the Lipschitz condition: $\enorm{\nabla f(x) - \nabla f(y)} \leq L\enorm{x-y}$ for all $x, y \in \mathbb{R}^n$.
 A function $f: \mathbb{R}^n \to \mathbb{R}$ is said to be convex if $f(y) \geq f(x) + \enp{\nabla f(x), y-x}$ for all $x, y \in \mathbb{R}^n$.
 For a convex function $f$, the $\epsilon$-subdifferential at $x$ is defined as $\partial_\epsilon f(x) = \{g \in \mathbb{R}^n : f(y) \geq f(x) + \enp{g, y-x} - \epsilon \text{ for all } y\in \mathbb{R}^n\}$.

\section{The A-DSPD Algorithm for the Three-Stage Problem}\label{sec:three-stage}

We develop the \ADSPD{} method for the three-stage stochastic optimization problem.
In Section~\ref{subsec:three-stage-reformulation}, we recast the original model as a sequence of stagewise saddle-point problems through value-function decomposition.
In Section~\ref{subsec:three-stage-iaspd}, we present the \iaspd{} method for the resulting SPPs and establish the properties needed for recursive implementation.
In Section~\ref{subsec:three-stage-convergence}, we combine these ingredients to derive convergence guarantees for the full \ADSPD{} algorithm in the three-stage setting.
Section~\ref{subsec:bounded-multipliers} gives a sufficient condition for the bounded optimal multiplier selections used in the analysis.
The same construction will then be extended to multi-stage problems with $T>3$ in the next section.

\subsection{Reformulation of the Three-Stage Problem}\label{subsec:three-stage-reformulation}
In this subsection, we reformulate the three-stage model as nested saddle-point problems and prepare the notation needed for the algorithmic development.
Throughout the stagewise reformulations in this paper, we assume that, for every realized scenario history and every predecessor decision in $X_{t-1}$ when $t\geq 2$, the corresponding stage problem is feasible, strong duality holds, and the primal and dual optima are attained. This condition justifies the nested saddle-point formulations introduced below.
We consider the following three-stage stochastic optimization problem:
\begin{equation}\label{three}
\begin{array}{lll}
\underset{x_1\in X_1}{\min} f_1(x_1, c_1)+&\bbe_{\xi_2\mid\xi_{[1]}}\left[\underset{x_2\in X_2}{\min} f_2\left(x_2, c_2\right)\right.+\quad & \left.\bbe_{\xi_3\mid\xi_{[2]}}\left[\underset{x_3\in X_3}{\min} f_3\left(x_3, c_3\right)\right]\right] \\
\text { s.t. } A_1 x_1- b_1\leq 0  & \text { s.t. } A_2 x_2+B_2 x_1-b_2\leq 0 , & \text { s.t. } A_3 x_3+B_3 x_2-b_3 \leq 0, \\
\end{array}
\end{equation}

We reformulate Problem~\eqref{three} as a sequence of nested saddle-point problems by proceeding backward from the final stage.
For a given second-stage decision $x_2$ and realization $\xi_{[3]}$, we define $V_3(x_2,\xi_{[3]})$ as the third-stage stochastic value function and $v_3(x_2,\xi_{[2]})$ as its conditional expectation given $\xi_{[2]}$:
	\begin{align}
	V_3(x_2,\xi_{[3]})&:= \min_{x_3\in X_3}\max_{y_3\in \bbr^{m_3}_+}\{\enp{A_3 x_3-b_3+B_3 x_2,y_3}+f_3(x_3,c_3)\}
	\label{third}\\
	v_3(x_2,\xi_{[2]})&:= \bbe_{\xi_3\mid\xi_{[2]}}[V_3(x_2,\xi_{[3]})]\notag.
	\end{align}
Substituting $v_3(x_2,\xi_{[2]})$ into the second-stage problem yields the stochastic value function $V_2(x_1,\xi_{[2]})$ and the corresponding conditional expectation $v_2(x_1,\xi_{[1]})$:
\begin{align}
V_2(x_1,\xi_{[2]})&:= \min_{x_2\in X_2}\max_{y_2\in \bbr^{m_2}_+}\{\enp{A_2 x_2-b_2+B_2 x_1,y_2}+f_2(x_2,c_2)+v_3(x_2,\xi_{[2]})\}
\label{second}\\
v_2(x_1,\xi_{[1]})&:= \bbe_{\xi_2\mid\xi_{[1]}}[V_2(x_1,\xi_{[2]})]\notag.
\end{align}
Finally, substituting $v_2(x_1,\xi_{[1]})$ into the first-stage problem gives an equivalent reformulation of Problem~\eqref{three}:
\begin{equation}\label{first}
	\begin{aligned}
	&\min_{x_1\in X_1} f_1(x_1,c_1)+v_2(x_1,\xi_{[1]})\\
	&\text { s.t. } A_1 x_1- b_1\leq 0.
	\end{aligned}
\end{equation}

For a target tolerance $\bar{\epsilon}>0$, we call a random first-stage decision a $\bar{\epsilon}$-solution if its expected absolute objective error relative to an optimal solution of~\eqref{first} and its expected constraint violation $\bbe[\|[A_1x_1-b_1]_+\|]$ are both at most $\bar{\epsilon}$. The symbol $\epsilon$ below denotes the inexactness of an individual value-function oracle. Whenever a scenario history is fixed, we suppress it in the arguments of the corresponding value functions.

To solve Problem~\eqref{first}, it is essential to compute first-order information for the value functions $v_2$ and $v_3$. However, since both $v_2$ and $v_3$ are defined by (conditional) expectations, evaluating their exact gradients or subgradients is often computationally intractable. To bypass this difficulty, we introduce the concept of a stochastic $\epsilon$-subgradient, as defined in \cite{lan2021dynamic}.
\begin{definition}\label{defi}
    Fix a scenario history $\xi_{[t-1]}$, a closed convex domain $X$ for the predecessor decision, a point $u\in X$, and $\epsilon\geq0$. For an integrable random vector $G(u,\xi_{[t]})$, define $g(u,\xi_{[t-1]}):=\bbe[G(u,\xi_{[t]})\mid\xi_{[t-1]}]$, where the expectation includes both the next-stage scenario and any internal oracle randomness.
    A random vector $G(u, \xi_{[t]})$ is called a stochastic $\epsilon$-subgradient of the value function $v_t(u, \xi_{[t-1]})$ if $g(u, \xi_{[t-1]})$ is an $\epsilon$-subgradient of $v_t(u, \xi_{[t-1]})$ with respect to $u$; that is,
	    \begin{equation*}
	    v_t(u^{\prime}, \xi_{[t-1]})\geq v_t(u, \xi_{[t-1]})+\enp{g(u, \xi_{[t-1]}),u^{\prime}-u}-\epsilon\quad \forall u^{\prime}\in X.
	    \end{equation*}
\end{definition}

Equations~\eqref{third}, \eqref{second}, and \eqref{first} reduce Problem~\eqref{three} to a sequence of stagewise problems, each coupled to the next stage only through its value function.
The subproblems in \eqref{third} and \eqref{second} motivate the analysis of the following generic saddle-point problem:
\begin{equation}\label{consistent}
\begin{aligned}
V(u,\xi)&:= \min_{x\in X}\max_{y\in \bbr^{m}_+}\{\enp{A x-b+Bu,y}+f(x,c)+v(x)\}\\
\end{aligned}
\end{equation}
Here, $u\in\mathbb{R}^{n_0}$ is a fixed predecessor decision, $\xi=(A,B,b,c)$ is the current data tuple, $X\subset\mathbb{R}^n$ is nonempty, compact, and convex, $A\in\mathbb{R}^{m\times n}$, $B\in\mathbb{R}^{m\times n_0}$, and $b\in\mathbb{R}^m$. The function $f(\cdot,c)$ is the smooth stage cost, and $v$ is the convex downstream value function, with $v=0$ at the terminal stage.
As an illustration, \eqref{third} is a special case of \eqref{consistent} with $n_0=n_2$, $n=n_3$, $m=m_3$, $u=x_2$, $y=y_3$, $X=X_3$, $b=b_3$, $B=B_3$, $A=A_3$, $f=f_3$, and $v=0$.
The first-stage problem can also be cast in this form, since \eqref{first} is equivalent to
\begin{equation}\label{thfirst}
	\min_{x_1\in X_1}\max_{y_1\in \bbr^{m_1}_+}\{\enp{A_1 x_1-b_1,y_1}+f_1(x_1,c_1)+v_2(x_1)\}
\end{equation}
A practical challenge in solving \eqref{consistent} is that the next-stage value function $v(x)$ generally has no closed-form expression, hence its subgradients are not directly available.
In each stagewise application, the realized scenario history is fixed and suppressed in the generic notation. Thus, \(v(x)=\bbe_{\xi'}[V'(x,\xi')]\) denotes an expectation over the next-stage scenario conditional on that history, where \(V'(\cdot,\cdot)\) has the same saddle-point structure as \eqref{consistent}.
Consequently, we must construct an approximate subgradient for the value-function term.
For $\bar z=(\bar x,\bar y)$ and $z=(x,y)$ in $X\times\mathbb{R}_+^m$, define the saddle-point residual
\begin{equation}\label{eq:saddle-residual}
Q(\bar z,z):=\enp{A\bar x-b+Bu,y}+f(\bar x,c)+v(\bar x)-\enp{Ax-b+Bu,\bar y}-f(x,c)-v(x).
\end{equation}
We also write $Q(\bar z;x,y):=Q(\bar{z},z)$. The following result of \citet{lan2021dynamic} provides a criterion for this construction:
\begin{lemma}[{\cite[Lemma 1]{lan2021dynamic}}]\label{subgradient}
	Let $(x^*,y^*)$ be a primal--dual optimal pair of \eqref{consistent}. For fixed $(u,\xi)$ and $(\bar{x},\bar{y})\in X\times\mathbb{R}_+^m$, if
	\begin{equation}\label{imp1}
	\begin{aligned}
	Q(\bar z;x,y^*)\leq\epsilon\qquad\text{for all }x\in X,
	\end{aligned}
	\end{equation}
	then $B^{\top}\bar{y}$ is an $\epsilon$-subgradient of $V(u,\xi)$ at $u$.
\end{lemma}
To facilitate the subsequent analysis, we define the following gap function:
\begin{equation}\label{gapdef}
\begin{aligned}
\operatorname{gap}(\bar{z}) := \max\{ Q(\bar{z}; x, y^*) : x \in X\},
\end{aligned}
\end{equation}
where $\bar{z} = (\bar{x}, \bar{y})$ and $y^*$ denotes an optimal dual solution.
If $\bbe[\operatorname{gap}(\bar{z})]\leq \epsilon$, then Lemma~\ref{subgradient} and Definition~\ref{defi} imply that $B^{\top}\bar{y}$ is a stochastic $\epsilon$-subgradient of $v(u)=\bbe_{\xi}[V(u,\xi)]$ at $u$.
Here the expectation includes both the next-stage scenario and internal oracle randomness, conditional on the fixed predecessor history and query point. Hence, controlling the gap provides the computable subgradient information required by the backward recursion in A-DSPD.
In Subsection~\ref{subsec:three-stage-iaspd}, we show how to compute a point $(\bar{x},\bar{y})$ satisfying this condition, and we establish bounds on dual second moments together with bounds on the optimality gap and primal infeasibility.

\subsection{The \iaspd{} Method}\label{subsec:three-stage-iaspd}

We develop \iaspd{} for solving the generic saddle-point problem~\eqref{consistent}.
We then establish finite-sample guarantees for $\epsilon$-subgradient computation, dual second moments, functional optimality, and primal infeasibility.
Finally, we provide three stage-specific parameter choices that are used later in the three-stage recursion.

\subsubsection{Algorithmic Update and iASPD Scheme}

The \iaspd{} method extends the deterministic accelerated primal--dual method introduced by~\citet{chen2014optimal}.
In this framework, each iteration consists of one call to the Accelerated Stochastic Primal--Dual Transformation (ASPDT) in Procedure~\ref{proc:ASPDT}.
At iteration $k$, $x^k$ and $y^k$ are the primal and dual iterates, $\tilde{x}^k$ and $\tilde{y}^k$ are their averaged outputs, $\bar{x}^k$ is the extrapolated primal point, and $w^k$ is the point at which the smooth gradient is evaluated. The averaging weight is $\beta_k\geq1$, the extrapolation coefficient is $\theta_{k+1}$, and the primal and dual step sizes are $\eta_k>0$ and $\tau_k>0$. The oracle vector $G(x^k)$ is a stochastic $\epsilon$-subgradient of $v$, with conditional mean $g(x^k):=\bbe_k[G(x^k)]$.
For first-stage calls, the input $\xi_1=(A_1,b_1,c_1)$ is interpreted as the generic data tuple $(A_1,0,b_1,c_1)$, with $u=0$ and $B=0$.
For compactness, the ASPDT call notation suppresses $\beta_k$ and $\theta_{k+1}$; each call uses these coefficients from the prescribed parameter schedule of its current stage.

\begin{procedure}[htbp]
	\caption{ASPDT for Generic Saddle-Point Problem~\eqref{consistent}}\label{proc:ASPDT}
	\begin{algorithmic}[0]
		\Require $(\tilde{x}^k,\tilde{y}^k,x^k,y^k,\bar{x}^k,G(x^k),u,\xi,X,\tau_k,\eta_k)$, where $\xi=(A,B,b,c)$
	\end{algorithmic}
	\begin{align}
	&w^k=(1-\beta_k^{-1}) \tilde{x}^k+\beta_k^{-1} x^k\label{PDHGx1}\\
	&y^{k+1}=\underset{y\in \bbr^m_+}{\arg\min} -\enp{A\bar{x}^{k}-b+B u,y}+\frac{1}{2\tau_{k}}\enorm{y-y^{k}}^2\label{PDHGity}\\
	&x^{k+1}= \underset{x\in X}{\arg\min}\enp{A x-b+Bu,y^{k+1}}+\enp{\nabla f(w^k,c)+G(x^{k}),x}+\frac{1}{2\eta_{k}}\enorm{x-x^{k}}^2
	\label{PDHGitx}\\
	&\tilde{x}^{k+1}=(1-\beta_k^{-1}) \tilde{x}^k+\beta_k^{-1} x^{k+1}\label{PDHGx2}\\
	&\tilde{y}^{k+1}=(1-\beta_k^{-1}) \tilde{y}^k+\beta_k^{-1} y^{k+1}\label{PDHGy1}\\
	&\bar{x}^{k+1}=\theta_{k+1}(x^{k+1}-x^k)+x^{k+1}\label{PDHGbarx}
	\end{align}
	\begin{algorithmic}[0]
		\Ensure $(\tilde{x}^{k+1},\tilde{y}^{k+1},\bar{x}^{k+1},x^{k+1},y^{k+1})$
	\end{algorithmic}
\end{procedure}

In Procedure~\ref{proc:ASPDT}, the stochastic $\epsilon$-subgradient $G(x^{k})$ (see Definition~\ref{defi}) and its conditional expectation $g(x^{k}) =\bbe_k[G(x^k)]$ provide inexact first-order information about the value function $v$ at the current search point $x^k$. The parameters $(u, \xi, X)$ describe the problem structure, while $(\beta_k, \theta_k, \eta_k,\tau_k)$ are algorithmic parameters that control the convergence behavior.

The procedure has four components:
\begin{itemize}
\item \textbf{Primal acceleration:} Equations \eqref{PDHGx1}, \eqref{PDHGitx}, and \eqref{PDHGx2} implement a multi-step Nesterov acceleration scheme \cite{tseng2010approximation} for the primal variable $\tilde{x}^k$.
\item \textbf{Dual update:} Equation \eqref{PDHGity} performs a dual space update to compute $y^{k+1}$.
\item \textbf{Dual averaging:} Equation \eqref{PDHGy1} defines a convex combination of $\tilde{y}^k$ and $y^{k+1}$ to generate the averaged dual variable $\tilde{y}^{k+1}$.
\item \textbf{Primal extrapolation:} Equation \eqref{PDHGbarx} defines a primal extrapolation (or prediction) step to estimate the primal variable $\bar{x}^{k+1}$ for the next dual update.
\end{itemize}

We assume that the subproblems in~\eqref{PDHGity} and~\eqref{PDHGitx} can be solved efficiently. When $v=0$ and $G(x^k)=0$, ASPDT reduces to the standard accelerated primal--dual method~\citep{chen2014optimal}. For additional details on primal--dual methods, refer to~\citet{zhu2024first,luo2021accelerated,xu2020primal,zhang2022primal}.
Building on ASPDT, we develop the \iaspd{} method in Algorithm~\ref{alg:A-ISPD}.

\begin{algorithm}[htbp]
	\caption{\iaspd{} Method for Generic Saddle-Point Problem \eqref{consistent}}\label{alg:A-ISPD}
		Given $\xi=(A,B,b,c)$.
	\begin{algorithmic}[0]
        \Require iteration budget $N\geq1$, initial points $x^1=\tilde{x}^1=\bar{x}^1\in X$, $y^1=\tilde{y}^1=0$, and parameter schedules
		\For{$k=1,\cdots,N$}
		\State Obtain a stochastic $\epsilon$-subgradient $G(x^k)$ of $v$ at $x^k$.
		\State $(\tilde{x}^{k+1},\tilde{y}^{k+1},\bar{x}^{k+1},x^{k+1},y^{k+1})=\operatorname{ASPDT}(\tilde{x}^k,\tilde{y}^k,x^{k},y^{k},\bar{x}^k,G(x^{k}),u,\xi,X,\tau_{k},\eta_k)$
		\EndFor		
		\Ensure $\tilde{z}^{N+1}\equiv (\tilde{x}^{N+1},\tilde{y}^{N+1})$.
	\end{algorithmic}  
\end{algorithm}  

Relative to~\citet{chen2014optimal}, \iaspd{} accommodates inexact first-order information because, in general, $g(x^k) \notin \partial v(x^k)$ and $G(x^k)$ is treated as a stochastic $\epsilon$-subgradient at the current iterate.
When $\epsilon=0$, the conditional mean $g(x^k)$ is an exact subgradient of $v$. If the oracle itself returns an exact subgradient without randomness, no stochastic oracle error remains. At the terminal stage, $v=0$ and we take $G(x^k)=0$.
Our analysis also derives complexity guarantees for $\epsilon$-subgradient construction, functional optimality, and primal infeasibility, together with bounds on the expected squared norms of the dual iterates.

\subsubsection{Assumptions and One-Step Recursion}

In order to analyze the \iaspd{} method, we impose the following assumptions.
\begin{assu}\label{ass1}
	Fix an oracle tolerance $\epsilon\geq0$, constants $L_f,M\geq0$, and a diameter bound $\Omega_X$. The following conditions hold:
	\begin{itemize}
	\item[$\operatorname{(i)}$] $f(\cdot,c)$ is a convex and smooth function, i.e.,
	\begin{equation}\label{Lif}
	0 \le f(x',c)-f(x,c)-\enp{\nabla f(x,c),x'-x} \le \frac{L_f}{2}\enorm{x'-x}^2,\quad \forall x,x'\in X.
	\end{equation}
	
	\item[$\operatorname{(ii)}$] At iteration $k$, the randomness used to construct $G(x^k)$ is conditionally independent of the preceding oracle randomness given $x^k$ and the current problem data, and $g(x^k):=\mathbb{E}_k[G(x^k)]\in\partial_\epsilon v(x^k)$ with $\mathbb{E}_k[\norm{G(x^k)}^2]\leq M^2$ almost surely.
	Moreover, $v$ is proper, closed, convex and Lipschitz continuous:
	\begin{equation}\label{Gbound}
	\abs{v(x')-v(x)}\leq M\enorm{x'-x},\quad \forall x,x'\in X.
	\end{equation} 
	
	\item[$\operatorname{(iii)}$] $X$ is bounded, i.e.,
	$
	\sup _{x_1, x_2 \in X} \enorm{x_1-x_2} \leq \Omega_X.
	$
	\end{itemize}
\end{assu}
Assumption~\ref{ass1}(ii), together with Jensen's inequality, implies that $\enorm{g(x^k)} \leq M$ almost surely.
For notational simplicity, we use the same constant $M$ in the second-moment bound of $G(x^k)$ and in the Lipschitz bound for $v$.

Write $\tilde z^k:=(\tilde x^k,\tilde y^k)$ and define the centered oracle error $\Delta^k:=g(x^k)-G(x^k)$, where $g(x^k)=\bbe_k[G(x^k)]$. Below, we establish the relation between $(\tilde{x}^{k+1},\tilde{y}^{k+1},\bar{x}^{k+1})$ and $(\tilde{x}^{k},\tilde{y}^{k},\bar{x}^{k})$ after one step of ASPDT.
\begin{pro}
   Let Assumption \ref{ass1} hold, $z=(x,y)$ and $\beta_k \geq 1$ for all $k\geq 1$. If $\tilde{z}^{k+1}=(\tilde{x}^{k+1}, \tilde{y}^{k+1})$ is generated by Algorithm~\ref{alg:A-ISPD}, then, for any $y\geq 0$ and $x\in X$, we have 
 		\begin{equation}\label{teq9}
 \begin{aligned}
 &\beta_k Q(\tilde{z}^{k+1}, z)-(\beta_k-1) Q(\tilde{z}^k, z)\\
 &\leq \frac{1}{2\eta_k}\enorm{x-x^k}^2 -\frac{1}{2\eta_k}\enorm{x-x^{k+1}}^2-\left(\frac{1}{2\eta_k}-\frac{L_f}{2\beta_k}\right)\enorm{x^{k+1}-x^k}^2\\
 &\quad + \frac{1}{2\tau_k}\norm{y-y^k}^2-\frac{1}{2\tau_k}\norm{y-y^{k+1}}^2-\frac{1}{2\tau_k}\norm{y^{k+1}-y^k}^2 \\
 &\quad - \enp{x^{k+1}-x, A^\top y^{k+1}}+\enp{A\bar{x}^k,y^{k+1}-y}+\enp{Ax^{k+1},y}-\enp{Ax,y^{k+1}}\\
 &\quad + (M+\norm{G(x^k)})\enorm{x^{k}-x^{k+1}}+\enp{\Delta^{k},x^{k}-x}+\epsilon,
 \end{aligned}
 \end{equation}
\end{pro}	
\begin{proof}
		From \eqref{PDHGx1} and \eqref{PDHGx2}, $\tilde{x}^{k+1}-w^k=\beta^{-1}_k (x^{k+1}-x^k)$. Using this observation and the Lipschitz continuity of $\nabla f(\cdot,c)$ in \eqref{Lif}, we have 
	\begin{equation}\label{teq1}
	\begin{aligned}
	\beta_k f(\tilde{x}^{k+1},c)&\leq \beta_k f(w^k,c)+\beta_k\enp{\tilde{x}^{k+1}-w^k,\nabla f(w^k,c)}+\frac{\beta_kL_f}{2}\enorm{\tilde{x}^{k+1}-w^k}^2\\
	&= \beta_k f(w^k,c)+\beta_k\enp{\tilde{x}^{k+1}-w^k,\nabla f(w^k,c)}+\frac{L_f}{2\beta_k}\norm{x^{k+1}-x^k}^2\\
	&= \beta_k f(w^k,c)+(\beta_k-1)\enp{\tilde{x}^{k}-w^k,\nabla f(w^k,c)}\\
	&\quad + \enp{x^{k+1}-w^k,\nabla f(w^k,c)}+\frac{L_f}{2\beta_k}\enorm{x^{k+1}-x^k}^2\\
	&= (\beta_k-1) [f(w^k,c)+\enp{\tilde{x}^{k}-w^k,\nabla f(w^k,c)}]\\
	&\quad + f(w^k,c)+\enp{x^{k+1}-w^k,\nabla f(w^k,c)}+\frac{L_f}{2\beta_k}\enorm{x^{k+1}-x^k}^2\\
	&= (\beta_k-1) [f(w^k,c)+\enp{\tilde{x}^{k}-w^k,\nabla f(w^k,c)}]\\
	&\quad + f(w^k,c)+\enp{x-w^k,\nabla f(w^k,c)}+\enp{x^{k+1}-x,\nabla f(w^k,c)}+\frac{L_f}{2\beta_k}\enorm{x^{k+1}-x^k}^2\\
	&\leq (\beta_k-1)f(\tilde{x}^k,c)+f(x,c)+\enp{x^{k+1}-x,\nabla f(w^k,c)}+\frac{L_f}{2\beta_k}\enorm{x^{k+1}-x^k}^2,
	\end{aligned}
	\end{equation}
	where the second equality is from \eqref{PDHGx2} and the last inequality comes from the convexity of $f(\cdot,c)$. 
    Moreover, by~\eqref{PDHGx2} and the convexity of $v(\cdot)$, we have
	\begin{equation}\label{teq2}
	\begin{aligned}
	\beta_k v(\tilde{x}^{k+1})-\beta_k v(x) &\leq (\beta_k-1) v(\tilde{x}^k)+v(x^{k+1})-\beta_k v(x)\\
	&= (\beta_k-1)(v(\tilde{x}^k)-v(x))+v(x^{k+1})-v(x).
	\end{aligned}
	\end{equation}	
	With straightforward algebraic manipulation using \eqref{PDHGy1}, we obtain
	\begin{equation}\label{teq3}
		\begin{aligned}
		& \beta_k\enp{-b+B u, y-\tilde{y}^{k+1}} - (\beta_k-1)\enp{-b+B u, y-\tilde{y}^k} \\
		&= \enp{-b+B u, y - (\beta_k \tilde{y}^{k+1} - (\beta_k-1)\tilde{y}^k)} \\
		&= \enp{-b+B u, y - y^{k+1}}.
		\end{aligned}
	\end{equation}
	From \eqref{teq1}, \eqref{teq2}, \eqref{teq3}, and the definition of $Q$ in~\eqref{eq:saddle-residual}, we obtain
	\begin{equation}\label{teq4}
		\begin{aligned}
		& \beta_k Q(\tilde{z}^{k+1}, z)-(\beta_k-1) Q(\tilde{z}^k, z)\\
		&= \beta_k[\enp{A\tilde{x}^{k+1}-b+B u,y}+f(\tilde{x}^{k+1},c)+v(\tilde{x}^{k+1})]-\beta_k[\enp{Ax-b+B u,\tilde{y}^{k+1}}+f(x,c)+v(x)]\\
		&\quad - (\beta_k-1)[\enp{A\tilde{x}^k-b+B u,y}+f(\tilde{x}^k,c)+v(\tilde{x}^k)]+(\beta_k-1)[\enp{Ax-b+B u,\tilde{y}^k}+f(x,c)+v(x)]\\
		&= \beta_k f(\tilde{x}^{k+1},c)-(\beta_k-1)f(\tilde{x}^k,c)-f(x,c)+\beta_k(v(\tilde{x}^{k+1})-v(x))-(\beta_k-1)(v(\tilde{x}^k)-v(x))\\
		&\quad + \enp{A(\beta_k \tilde{x}^{k+1}-(\beta_k-1)\tilde{x}^k),y}-\enp{Ax,\beta_k\tilde{y}^{k+1}-(\beta_k-1)\tilde{y}^k}\\
		&\quad + \beta_k\enp{-b+B u,y-\tilde{y}^{k+1}}-(\beta_k-1)\enp{-b+B u,y-\tilde{y}^k}\\
		&\leq \enp{\nabla f(w^k,c),x^{k+1}-x}+\frac{L_f}{2\beta_k}\enorm{x^{k+1}-x^k}^2+v(x^{k+1})-v(x)\\
		&\quad + \enp{Ax^{k+1},y}-\enp{Ax,y^{k+1}}+\enp{-b+B u,y-y^{k+1}}.
		\end{aligned}
	\end{equation}
	Moreover, utilizing the Lipschitz continuity of $v(x)$ and the definition of the $\epsilon$-subgradient, we obtain
	\begin{equation}\label{heq1}
	\begin{aligned}
	v(x^{k+1}) &\leq v(x^k) + M\enorm{x^{k+1} - x^k} \\
	&\leq v(x) + \enp{g(x^k), x^k - x} + M\enorm{x^{k+1} - x^k} + \epsilon.
	\end{aligned}
	\end{equation}
	By the definition of $\Delta^k$, we have
	\begin{equation}\label{heq2}
	\begin{aligned}
	\enp{g(x^{k}), x^{k} - x}
	&= \enp{G(x^{k}), x^{k} - x} + \enp{\Delta^{k}, x^{k} - x} \\
	&= \enp{G(x^{k}), x^{k+1} - x} + \enp{G(x^{k}), x^{k} - x^{k+1}} + \enp{\Delta^{k}, x^{k} - x} \\
	&\leq \enp{G(x^{k}), x^{k+1} - x} + \enorm{G(x^{k})}\enorm{x^{k} - x^{k+1}} + \enp{\Delta^{k}, x^{k} - x},
	\end{aligned}
	\end{equation}
	where the final inequality follows from the Cauchy--Schwarz inequality.
	Summing up \eqref{heq2} and~\eqref{heq1}, we have
	\begin{equation}\label{heq5}
	\begin{aligned}
	v(x^{k+1}) - v(x)
	\leq \enp{G(x^{k}), x^{k+1} - x}+(M + \enorm{G(x^k)}) \enorm{x^{k} - x^{k+1}} + \enp{\Delta^{k}, x^{k} - x} + \epsilon.
	\end{aligned}
	\end{equation}
	Combining \eqref{heq5} with~\eqref{teq4}, we obtain
	\begin{equation}\label{teq6}
		\begin{aligned}
		&\beta_k Q(\tilde{z}^{k+1}, z)-(\beta_k-1) Q(\tilde{z}^k, z)\\
		&\leq \enp{\nabla f(w^k,c)+G(x^{k}), x^{k+1}-x}+\frac{L_f}{2\beta_k}\enorm{x^{k+1}-x^k}^2\\	
		&\quad + \enp{Ax^{k+1},y}-\enp{Ax,y^{k+1}}+\enp{-b+B u,y-y^{k+1}}\\
		&\quad + (M+\enorm{G(x^k)})\enorm{x^{k}-x^{k+1}}+\enp{\Delta^{k},x^{k}-x}+\epsilon.
		\end{aligned}
	\end{equation}	
	By invoking the optimality conditions of \eqref{PDHGity}, we obtain
	\begin{equation}\label{teq7}
	\begin{aligned}
	\enp{A\bar{x}^{k}-b+B u, y-y^{k+1}} \leq \frac{1}{2\tau_k}\enorm{y-y^k}^2 - \frac{1}{2\tau_k}\enorm{y-y^{k+1}}^2 - \frac{1}{2\tau_k}\enorm{y^{k+1}-y^k}^2,
	\end{aligned}
	\end{equation}
	which holds for any $y\in\mathbb{R}_+^m$. Similarly, the optimality condition of \eqref{PDHGitx} yields
	\begin{equation}\label{teq8}
	\begin{aligned}
	&\enp{\nabla f(w^k,c)+G(x^{k}),x^{k+1}-x} + \enp{x^{k+1}-x,A^\top y^{k+1}} \\
	&\leq \frac{1}{2\eta_k}\enorm{x-x^k}^2 - \frac{1}{2\eta_k}\enorm{x-x^{k+1}}^2 - \frac{1}{2\eta_k}\enorm{x^{k+1}-x^k}^2, \quad \forall x\in X.
	\end{aligned}
	\end{equation}
	Summing \eqref{teq6}, \eqref{teq7}, and \eqref{teq8}, we arrive at
	\begin{equation*}
	\begin{aligned}
	&\beta_k Q(\tilde{z}^{k+1}, z) - (\beta_k-1) Q(\tilde{z}^k, z) \\
	&\leq \frac{1}{2\eta_k}\enorm{x-x^k}^2 - \frac{1}{2\eta_k}\enorm{x-x^{k+1}}^2 - \left(\frac{1}{2\eta_k} - \frac{L_f}{2\beta_k}\right)\enorm{x^{k+1}-x^k}^2 \\
	&\quad + \frac{1}{2\tau_k}\enorm{y-y^k}^2 - \frac{1}{2\tau_k}\enorm{y-y^{k+1}}^2 - \frac{1}{2\tau_k}\enorm{y^{k+1}-y^k}^2 \\
	&\quad - \enp{x^{k+1}-x, A^\top y^{k+1}} + \enp{A\bar{x}^k, y^{k+1}-y} + \enp{Ax^{k+1}, y} - \enp{Ax, y^{k+1}} \\
	&\quad + (M + \enorm{G(x^k)})\enorm{x^{k} - x^{k+1}} + \enp{\Delta^{k}, x^{k} - x} + \epsilon,
	\end{aligned}
	\end{equation*}
	which completes the proof.
\end{proof}

\subsubsection{Finite-Sample Guarantees for iASPD}

We are now ready to present the main finite-sample results for the \iaspd{} method.
\begin{theorem}\label{th1}
	Suppose Assumption~\ref{ass1} holds. Let $\eta_k,\tau_k>0$ and $\beta_k\geq1$ for $1\leq k\leq N$, and set $\theta_1=0$. Suppose the parameters in Algorithm~\ref{alg:A-ISPD} and the auxiliary weights $\gamma_k$ satisfy the following relations over the indicated indices:
	\begin{align}
	&\theta_k = \frac{\tau_{k-1}}{\tau_k} \leq \frac{\eta_{k-1}}{\eta_k},\quad2\leq k\leq N, \label{PDHGtheta}\\
	&\gamma_k= \begin{cases}1, & k=1 \\ \theta_k^{-1} \gamma_{k-1}, & k \geq 2\end{cases}, \label{PDHGgamma}\\
	&\beta_1=1, \quad \beta_{k+1}-1 = \beta_k \theta_{k+1},\quad 1\leq k\leq N, \label{PDHGbeta}\\
	&\frac{1}{\eta_k} - \frac{L_f}{\beta_k} - 2\norm{A}^2\tau_k > 0,\quad 1\leq k\leq N. \label{PDHGka}
	\end{align}
	Then, for every $z=(x,y)\in X\times\mathbb{R}_+^m$ and $N\geq1$, the following inequality holds:
	\begin{equation}\label{thimq1}
	\begin{aligned}
		Q(\tilde{z}^{N+1}, z)
		\leq \frac{\Omega_X^2}{2\beta_N \eta_N}
		+ \frac{1}{2\beta_N \tau_N}\norm{y}^2 - \frac{1}{4\beta_N \tau_N}\norm{y - y^{N+1}}^2
		+ \frac{1}{\beta_N \gamma_N} \sum_{k=1}^{N} \Lambda_k(x) + \epsilon,
	\end{aligned}
	\end{equation}
	where 
	\begin{equation}\label{Lambda}
		\Lambda_k(x) = \frac{ \gamma_k (M + \norm{G(x^k)})^2 }{ 2( \eta_k^{-1} - L_f\beta_k^{-1} - 2\enorm{A}^2\tau_k ) }
		+ \gamma_k \enp{\Delta^k, x^k - x}.
	\end{equation}
\end{theorem}
\begin{proof}
For the initial index, use the conventions $x^0=x^1$, $\gamma_0=0$, and $\tau_0=0$. Thus all terms involving $x^1-x^0$ vanish. Relations between consecutive positive step sizes below are used only for $k\geq2$. In view of \eqref{PDHGbarx}, we have
\begin{equation*}
\begin{aligned}
	&\quad - \enp{x^{k+1}-x, A^\top y^{k+1}}+\enp{A\bar{x}^k,y^{k+1}-y}+\enp{Ax^{k+1},y}-\enp{Ax,y^{k+1}}\\
	&= \enp{A(x^{k+1}-x^k),y-y^{k+1}}-\theta_k\enp{A(x^{k}-x^{k-1}),y-y^{k+1}}\\
	&= \enp{A(x^{k+1}-x^k),y-y^{k+1}}-\theta_k\enp{A(x^{k}-x^{k-1}),y-y^{k}}-\theta_k\enp{A(x^{k}-x^{k-1}),y^k-y^{k+1}}.
\end{aligned}
\end{equation*}
	Multiplying both sides of \eqref{teq9} by $\gamma_k$, using the above identity and the fact that $\gamma_k \theta_k=\gamma_{k-1}$ (i.e., \eqref{PDHGgamma}), we obtain
	\begin{equation}\label{teq10}
		\begin{aligned}
		&\beta_k\gamma_k(Q(\tilde{z}^{k+1}, z)-\epsilon) -(\beta_k-1)\gamma_k (Q(\tilde{z}^k, z)-\epsilon)\\
		 &\leq \frac{\gamma_k}{2\eta_k}\enorm{x-x^k}^2-\frac{\gamma_k}{2\eta_k}\enorm{x-x^{k+1}}^2
		 +\frac{\gamma_k}{2\tau_k}\norm{y-y^k}^2-\frac{\gamma_k}{2\tau_k}\enorm{y-y^{k+1}}^2\\
		&\quad - \gamma_k\left(\frac{1}{2\eta_k}-\frac{L_f}{2\beta_k}\right)\norm{x^{k+1}-x^k}^2-\frac{\gamma_k}{2\tau_k}\norm{y^{k+1}-y^k}^2 
		-\gamma_{k-1}\enp{A(x^{k}-x^{k-1}),y^k-y^{k+1}}\\
		&\quad + \gamma_k\enp{A(x^{k+1}-x^k),y-y^{k+1}}-\gamma_{k-1}\enp{A(x^{k}-x^{k-1}),y-y^{k}}\\
		&\quad + \gamma_k(M+\enorm{G(x^k)})\enorm{x^{k}-x^{k+1}}+\gamma_k\enp{\Delta^{k},x^{k}-x}.
		\end{aligned}
	\end{equation}
	Applying the Cauchy--Schwarz inequality to the term $-\gamma_{k-1}\enp{A(x^{k}-x^{k-1}),y^k-y^{k+1}}$ in the above inequality and noticing that $\gamma_{k-1} / \gamma_k=\theta_k = \tau_{k-1} / \tau_k $ from \eqref{PDHGtheta}, we have
	\begin{equation*}
		\begin{aligned}
		-\gamma_{k-1}\enp{A(x^{k}-x^{k-1}),y^k-y^{k+1}}
		&\leq \gamma_{k-1}\enorm{A}\norm{x^{k}-x^{k-1}}\enorm{y^k-y^{k+1}}\\
		&\leq \frac{\enorm{A}^2\gamma^2_{k-1}\tau_{k}}{\gamma_{k}}\enorm{x^k-x^{k-1}}^2+\frac{\gamma_{k}}{4\tau_{k}}\enorm{y^k-y^{k+1}}^2\\
		&\leq (\enorm{A}^2\gamma_{k-1}\tau_{k-1})\enorm{x^k-x^{k-1}}^2+\frac{\gamma_{k}}{4\tau_{k}}\enorm{y^k-y^{k+1}}^2.
		\end{aligned}
	\end{equation*}
		Noting that $\theta_{k+1}=\gamma_k / \gamma_{k+1}$ and using \eqref{PDHGbeta}, we have $\left(\beta_{k+1}-1\right) \gamma_{k+1}=\beta_k \gamma_k$. Combining the above two relations with \eqref{teq10}, we get the following recursion
			\begin{equation*}
		\begin{aligned}
		&(\beta_{k+1}-1)\gamma_{k+1}(Q(\tilde{z}^{k+1}, z)-\epsilon) -(\beta_k-1)\gamma_k (Q(\tilde{z}^k, z)-\epsilon)\\
		 &\leq \frac{\gamma_k}{2\eta_k}\enorm{x-x^k}^2-\frac{\gamma_k}{2\eta_k}\enorm{x-x^{k+1}}^2
		+\frac{\gamma_k}{2\tau_k}\enorm{y-y^k}^2-\frac{\gamma_k}{2\tau_k}\enorm{y-y^{k+1}}^2\\
		&\quad - \gamma_k\left(\frac{1}{2\eta_k}-\frac{L_f}{2\beta_k}\right)\enorm{x^{k+1}-x^k}^2
		+\enorm{A}^2\gamma_{k-1}\tau_{k-1}\enorm{x^k-x^{k-1}}^2\\
		&\quad + \gamma_k\enp{A(x^{k+1}-x^k),y-y^{k+1}}-\gamma_{k-1}\enp{A(x^{k}-x^{k-1}),y-y^{k}}\\
		&\quad + \gamma_k(M+\norm{G(x^k)})\enorm{x^{k}-x^{k+1}}+\gamma_k\enp{\Delta^{k},x^{k}-x}.
		\end{aligned}
		\end{equation*}
		Summing up the above inequality from 1 to $N$ and assuming that $x^0=x^1$, we conclude that
		\begin{equation}\label{zq3}
			\begin{aligned}
			& \beta_{N} \gamma_{N} (Q(\tilde{z}^{N+1}, z)-\epsilon)\\
			&= (\beta_{N+1}-1) \gamma_{N+1}(Q(\tilde{z}^{N+1}, z)-\epsilon)
			-(\beta_1-1) \gamma_1( Q(\tilde{z}^1, z)-\epsilon)\\
			&\leq \sum_{k=1}^{N}\frac{\gamma_k}{2\eta_k}(\enorm{x-x^k}^2-\enorm{x-x^{k+1}}^2)
			+\sum_{k=1}^{N}\frac{\gamma_k}{2\tau_k}(\enorm{y-y^k}^2-\enorm{y- y^{k+1}}^2)\\
			&\quad - \gamma_N\left(\frac{1}{2\eta_N}-\frac{L_f}{2\beta_N}\right)\norm{x^{N+1}-x^N}^2
			-\sum_{k=1}^{N-1}\gamma_k \left(\frac{1}{2\eta_k}-\frac{L_f}{2\beta_k}-\enorm{A}^2\tau_{k}\right)\enorm{x^{k+1}-x^k}^2\\
			&\quad + \sum_{k=1}^{N}\gamma_k(M+\norm{G(x^k)})\enorm{x^{k}-x^{k+1}}+\sum_{k=1}^{N}\gamma_k\enp{\Delta^{k},x^{k}-x}\\
			&\quad + \gamma_N\enp{A(x^{N+1}-x^N),y-y^{N+1}},
			\end{aligned}
		\end{equation}
		where the first inequality follows from $(\beta_{N+1}-1) \gamma_{N+1}=\beta_{N} \gamma_{N}$ and $\beta_1=1$. By the definition of $\gamma_{k}$ and 
		the relation \eqref{PDHGtheta}, we have $\frac{\gamma_k}{\eta_k}\leq\frac{\gamma_{k+1}}{\eta_{k+1}}$ 
			and $\frac{\gamma_k}{\tau_k}=\frac{\gamma_{k+1}}{\tau_{k+1}}$. Using this observation and $y^1=0$, we conclude that
		\begin{equation}\label{zq1}
			\begin{aligned}
				&\sum_{k=1}^{N}\frac{\gamma_k}{2\eta_k}(\enorm{x-x^k}^2-\enorm{x-x^{k+1}}^2)
			+\sum_{k=1}^{N}\frac{\gamma_k}{2\tau_k}(\enorm{y-y^k}^2-\enorm{y- y^{k+1}}^2)\\
			&= \frac{\gamma_1}{2\eta_1}\enorm{x-x^1}^2-\sum_{k=1}^{N-1}\left(\frac{\gamma_k}{2\eta_k}-\frac{\gamma_{k+1}}{2\eta_{k+1}}\right)\enorm{x-x^{k+1}}^2
			-\frac{\gamma_N}{2\eta_N}\enorm{x-x^{N+1}}^2\\
			&\quad + \frac{\gamma_N}{2\tau_N}\enorm{y}^2-\frac{\gamma_N}{2\tau_N}\enorm{y-y^{N+1}}^2\\
			&\leq \frac{\gamma_1\Omega^2_X}{2\eta_1}-\sum_{k=1}^{N-1}\left(\frac{\gamma_k}{2\eta_k}-\frac{\gamma_{k+1}}{2\eta_{k+1}}\right)\Omega^2_X
			+\frac{\gamma_N}{2\tau_N}\enorm{y}^2-\frac{\gamma_N}{2\tau_N}\enorm{y-y^{N+1}}^2\\
			&\leq \frac{\gamma_{N}\Omega^2_X}{2\eta_{N}}
			+\frac{\gamma_N}{2\tau_N}\enorm{y}^2-\frac{\gamma_N}{2\tau_N}\enorm{y-y^{N+1}}^2.
			\end{aligned}
		\end{equation}
		Now applying the Cauchy--Schwarz inequality to $\gamma_N\enp{A(x^{N+1}-x^N),y-y^{N+1}}$, we get
		\begin{equation}\label{zq2}
			\begin{aligned}
			\gamma_N\enp{A(x^{N+1}-x^N),y-y^{N+1}}
			&\leq \gamma_{N}\enorm{A}\norm{x^{N+1}-x^N}\enorm{y-y^{N+1}}\\
			&\leq \enorm{A}^2\gamma_{N}\tau_N\enorm{x^{N+1}-x^N}^2+\frac{\gamma_N}{4\tau_N}\enorm{y-y^{N+1}}^2.
			\end{aligned}
		\end{equation}
		Combining \eqref{zq1} and \eqref{zq2} with \eqref{zq3}, we have 
			\begin{equation*}
		\begin{aligned}
			\beta_{N} \gamma_{N} (Q(\tilde{z}^{N+1}, z)-\epsilon)
		&\leq \frac{\gamma_N\Omega^2_X}{2\eta_N}
		+\frac{\gamma_N}{2\tau_N}\enorm{y}^2-\frac{\gamma_N}{4\tau_N}\enorm{y-y^{N+1}}^2\\
		&\quad - \sum_{k=1}^{N}\gamma_k \left(\frac{1}{2\eta_k}-\frac{L_f}{2\beta_k}-\norm{A}^2\tau_{k}\right)\enorm{x^{k+1}-x^k}^2\\
		&\quad + \sum_{k=1}^{N}\gamma_k(M+\enorm{G(x^k)})\enorm{x^{k}-x^{k+1}}+\sum_{k=1}^{N}\gamma_k\enp{\Delta^{k},x^{k}-x}\\
		&\leq \frac{\gamma_N\Omega_X^2}{2\eta_N}
		+\frac{\gamma_N}{2\tau_N}\norm{y}^2-\frac{\gamma_N}{4\tau_N}\norm{y-y^{N+1}}^2\\
		&\quad + \sum_{k=1}^{N} \frac{\gamma_k(M+\enorm{G(x^k)})^2}{2(\eta_k^{-1}-L_f\beta_k^{-1}-2\enorm{A}^2\tau_{k})}
		+\sum_{k=1}^{N}\gamma_k\enp{\Delta^{k},x^{k}-x},
		\end{aligned}
		\end{equation*}
		where the second inequality follows from $-at^2+bt\leq \frac{b^2}{4a}$, the parameter setting in \eqref{PDHGka}, and the boundedness of $X$. Dividing both sides of the above inequality by $\beta_{N}\gamma_{N}$ yields \eqref{thimq1}.
\end{proof}

\citet{lan2021dynamic} establish the following technical result, which bounds the cumulative error due to stochastic subgradient estimation.
\begin{lemma}\label{shmmmle}
	Let $\alpha_k>0$ be auxiliary parameters, define $\hat{x}^1=x^1$ and
	\begin{equation}\label{sgmlemma1}
	\hat{x}^{k+1} = \underset{x \in X}{\arg\min}\left\{\gamma_k\enp{\Delta^k, x} + \frac{\alpha_k}{2}\norm{x - \hat{x}^k}^2\right\}
	\end{equation}
	for any $k \geq 1$. Then, for any $x \in X$, the following inequality holds:
	\begin{equation}\label{Delta}
		\sum_{k=1}^N \gamma_k\enp{\Delta^k, \hat{x}^{k} - x} \leq \sum_{k=1}^N \frac{\alpha_k}{2}(\enorm{x - \hat{x}^{k}}^2 - \enorm{x - \hat{x}^{k+1}}^2) + \sum_{k=1}^N \frac{\gamma_k^2 \enorm{\Delta^k}^2}{2 \alpha_k}.
	\end{equation}
\end{lemma}

The following theorem provides explicit finite-sample upper bounds for the gap function defined in \eqref{gapdef}:

\begin{theorem}\label{thzb}
Suppose that all the assumptions in Theorem \ref{th1} hold, that $\alpha_k>0$ for $1\leq k\leq N$, and that $\alpha_k \geq \alpha_{k+1}$, $\beta_{k}\tau_{k}\leq \beta_{k+1}\tau_{k+1}$,
$\tau_k/\gamma_{k}\leq \tau_{k+1}/\gamma_{k+1}$ for all $1\leq k<N$.
Then, the following results hold:
	\begin{itemize}
		\item[$\operatorname{(a)}$] For any $N \geq 1$, it holds that
		\begin{equation*}
		\begin{aligned}
		\bbe[\operatorname{gap}(\tilde{z}^{N+1})]&
		\leq \frac{\Omega^2_X}{2\beta_{N}\eta_N}+\frac{\alpha_1 \Omega^2_X}{2\beta_{N}\gamma_{N}}+\frac{1}{2\beta_{N}\tau_N}\enorm{y^*}^2\\
		&\quad + \frac{1}{\beta_{N}\gamma_{N}}\sum_{k=1}^{N}\frac{2\gamma_k M^2}{\eta_k^{-1}-L_f\beta_k^{-1}-2\norm{A}^2\tau_{k}}
		+\frac{1}{\beta_{N}\gamma_{N}}\sum_{k=1}^{N}\frac{2\gamma^2_k M^2}{\alpha_k} +\epsilon,
		\end{aligned}
		\end{equation*}
		
		\item[$\operatorname{(b)}$] For any $N \geq 1$, it holds that
		\begin{equation*}
		\begin{aligned}
		\bbe[\enorm{\tilde{y}^{N+1}}^2]
		&\leq \frac{4\Omega^2_X\tau_{N}}{\eta_N}+\frac{4\alpha_1 \Omega^2_X\tau_{N}}{\gamma_{N}}+6\enorm{y^*}^2\\
		&\quad + \frac{8\tau_{N}}{\gamma_{N}}\sum_{k=1}^{N}\frac{2\gamma_k M^2}{\eta_k^{-1}-L_f\beta_k^{-1}-2\enorm{A}^2\tau_{k}}
		+\frac{8\tau_{N}}{\gamma_{N}}\sum_{k=1}^{N}\frac{2\gamma^2_k M^2}{\alpha_k}+8\beta_{N}\tau_{N}\epsilon.
		\end{aligned}
		\end{equation*}

	\end{itemize}
\end{theorem}
\begin{proof}
	Taking $y=y^*$ in \eqref{thimq1}, we have 
	\begin{equation}\label{bbeq1}
		\begin{aligned}
		Q(\tilde{z}^{N+1}; x, y^*)
		\leq \frac{\Omega^2_X}{2\beta_{N}\eta_N}+\frac{1}{2\beta_{N}\tau_N}\enorm{y^*}^2-\frac{1}{4\beta_{N}\tau_N}\enorm{y^*-y^{N+1}}^2
		+\frac{1}{\beta_{N}\gamma_{N}}\sum_{k=1}^{N}\Lambda_k +\epsilon.
		\end{aligned}
	\end{equation}
	Taking the supremum of \eqref{bbeq1} over $x\in X$ and a full expectation, we have 
	\begin{equation}\label{bbzq2}
		\begin{aligned}
		\bbe[\operatorname{gap}(\tilde{z}^{N+1})]\leq \frac{\Omega^2_X}{2\beta_{N}\eta_N}+\frac{1}{2\beta_{N}\tau_N}\enorm{y^*}^2
		+\frac{1}{\beta_{N}\gamma_{N}}\bbe\left[\sup_{x\in X}\sum_{k=1}^{N}\Lambda_k(x)\right] +\epsilon.
		\end{aligned}
	\end{equation}
	Now it follows from \eqref{Lambda} that 
	\begin{equation}\label{sseq1}
		\begin{aligned}
		\sum_{k=1}^{N}\Lambda_k(x)&= \sum_{k=1}^{N}
		\frac{\gamma_k(M+\enorm{G(x^k)})^2}{2(\eta_k^{-1}-L_f\beta_k^{-1}-2\enorm{A}^2\tau_{k})}
		+\sum_{k=1}^{N}\gamma_k\enp{\Delta^{k},x^{k}-x}\\
		&= \sum_{k=1}^{N}
		\frac{\gamma_k(M+\enorm{G(x^k)})^2}{2(\eta_k^{-1}-L_f\beta_k^{-1}-2\enorm{A}^2\tau_{k})}
		+\sum_{k=1}^{N}\gamma_k\enp{\Delta^{k},x^{k}-\hat{x}^k}+\sum_{k=1}^{N}\gamma_k\enp{\Delta^{k},\hat{x}^k-x}.
		\end{aligned}
	\end{equation}
	The oracle randomness at iteration $k$ is conditionally independent of the preceding oracle randomness given $x^k$ and the current problem data, while $x^k$ and $\hat{x}^k$ are $\mathcal{F}_{k-1}$-measurable. Hence $\mathbb{E}_k[\Delta^k]=0$ and
	$\mathbb{E}_k\bsbra{\enp{\Delta^{k},x^{k}-\hat{x}^k}}=0$. Moreover, conditional Jensen's inequality and Assumption~\ref{ass1}(ii) give $\norm{g(x^{k})}\leq M$ and
	\begin{equation*}
		\begin{aligned}
		\mathbb{E}_k[\enorm{\Delta^k}^2]
		&= \mathbb{E}_k[\enorm{g(x^k)-G(x^k)}^2]\\
		&\leq 2(\enorm{g(x^k)}^2+\mathbb{E}_k[\enorm{G(x^k)}^2]) \leq 4 M^2.
		\end{aligned}
	\end{equation*}
	The same conditional moment bound also gives $\mathbb{E}_k[(M+\norm{G(x^k)})^2]\leq 4M^2$.
	Taking the expectation on both sides of \eqref{sseq1} over all the randomness, we therefore obtain
	\begin{equation}\label{lammmdd}
		\begin{aligned}
		\bbe\left[\sup_{x\in X}\sum_{k=1}^{N}\Lambda_k(x)\right]
		&\leq \sum_{k=1}^{N}\bbe\left[\frac{\gamma_k(M+\enorm{G(x^k)})^2}{2\left(\eta_k^{-1}-L_f\beta_k^{-1}-2\enorm{A}^2\tau_{k}\right)}\right]
		+\bbe\left[\sup_{x\in X}\sum_{k=1}^{N}\gamma_k\enp{\Delta^{k},\hat{x}^k-x}\right]\\
		&\leq \sum_{k=1}^{N}\frac{2\gamma_k M^2}{\eta_k^{-1}-L_f\beta_k^{-1}-2\enorm{A}^2\tau_{k}}
		+\frac{\alpha_1 \Omega^2_X}{2}+\sum_{k=1}^{N}\frac{2\gamma^2_k M^2}{\alpha_k},
		\end{aligned}
	\end{equation}
	where the last inequality comes from Lemma \ref{shmmmle}. Combining the above inequality with \eqref{bbzq2}, 
	we obtain
	\begin{equation*}
	\begin{aligned}
	\bbe[\operatorname{gap}(\tilde{z}^{N+1})]&\leq \frac{\Omega^2_X}{2\beta_{N}\eta_N}+\frac{\alpha_1 \Omega^2_X}{2\beta_{N}\gamma_{N}}+\frac{1}{2\beta_{N}\tau_N}\enorm{y^*}^2\\
	&\quad + \frac{1}{\beta_{N}\gamma_{N}}\sum_{k=1}^{N}\frac{2\gamma_k M^2}{\eta_k^{-1}-L_f\beta_k^{-1}-2\enorm{A}^2\tau_{k}}
	+\frac{1}{\beta_{N}\gamma_{N}}\sum_{k=1}^{N}\frac{2\gamma^2_k M^2}{\alpha_k} +\epsilon.
	\end{aligned}
	\end{equation*}
	Fixing $x=x^*$ in~\eqref{bbeq1} and using the fact that $Q\left(\tilde{z}^{N+1}; x^*, y^*\right) \geq 0$, we have
		\begin{equation*}
	\begin{aligned}
	\frac{1}{4\beta_{N}\tau_N}\norm{y^*-y^{N+1}}^2
	\leq \frac{\Omega^2_X}{2\beta_{N}\eta_N}+\frac{1}{2\beta_{N}\tau_N}\norm{y^*}^2
	+\frac{1}{\beta_{N}\gamma_{N}}\sum_{k=1}^{N}\Lambda_k(x) +\epsilon.
	\end{aligned}
	\end{equation*}
	Taking expectation over all the randomness, multiplying both sides by $4\beta_{N}\tau_N$ and then using \eqref{lammmdd},  we conclude
	\begin{equation}\label{abeq1}
\begin{aligned}
\bbe[\enorm{y^{N+1}-y^*}^2]
&\leq \frac{2\Omega^2_X\tau_{N}}{\eta_N}+\frac{2\alpha_1 \Omega^2_X\tau_{N}}{\gamma_{N}}+2\enorm{y^*}^2\\
&\quad + \frac{4\tau_{N}}{\gamma_{N}}\sum_{k=1}^{N}\frac{2\gamma_k M^2}{\eta_k^{-1}-L_f\beta_k^{-1}-2\enorm{A}^2\tau_{k}}
+\frac{4\tau_{N}}{\gamma_{N}}\sum_{k=1}^{N}\frac{2\gamma^2_k M^2}{\alpha_k}+4\beta_{N}\tau_{N}\epsilon.
\end{aligned}
\end{equation}
Using the fact that $b^2-2a^2 \leq 2(b-a)^2$, we have for any $N \geq 1$ that
\begin{equation}\label{zbeq1}
\begin{aligned}
\bbe[\enorm{y^{N+1}}^2]
&\leq \frac{4\Omega^2_X\tau_{N}}{\eta_N}+\frac{4\alpha_1 \Omega^2_X\tau_{N}}{\gamma_{N}}+6\enorm{y^*}^2\\
&\quad + \frac{8\tau_{N}}{\gamma_{N}}\sum_{k=1}^{N}\frac{2\gamma_k M^2}{\eta_k^{-1}-L_f\beta_k^{-1}-2\enorm{A}^2\tau_{k}}
+\frac{8\tau_{N}}{\gamma_{N}}\sum_{k=1}^{N}\frac{2\gamma^2_k M^2}{\alpha_k}+8\beta_{N}\tau_{N}\epsilon.
\end{aligned}
\end{equation}
Observe that \eqref{zbeq1} holds for any $y^{k+1}, k=1, \cdots, N$, and hence that
\begin{equation}\label{zbeq2}
\begin{aligned}
\bbe[\enorm{y^{k+1}}^2]
&\leq \frac{4\Omega^2_X\tau_{k}}{\eta_k}+\frac{4\alpha_1 \Omega^2_X\tau_{k}}{\gamma_{k}}+6\enorm{y^*}^2\\
&\quad + \frac{8\tau_{k}}{\gamma_{k}}\sum_{i=1}^{k}\frac{2\gamma_i M^2}{\eta_i^{-1}-L_f\beta_i^{-1}-2\enorm{A}^2\tau_{i}}
+\frac{8\tau_{k}}{\gamma_{k}}\sum_{i=1}^{k}\frac{2\gamma^2_i M^2}{\alpha_i}+8\beta_{k}\tau_{k}\epsilon.
\end{aligned}
\end{equation}
Denote 
$ 
R_{k+1}:=\frac{4\Omega^2_X\tau_{k}}{\eta_k}+\frac{4\alpha_1 \Omega^2_X\tau_{k}}{\gamma_{k}}+6\enorm{y^*}^2
 +\frac{8\tau_{k}}{\gamma_{k}}\sum_{i=1}^{k}\frac{2\gamma_i M^2}{\eta_i^{-1}-L_f\beta_i^{-1}-2\enorm{A}^2\tau_{i}}
+\frac{8\tau_{k}}{\gamma_{k}}\sum_{i=1}^{k}\frac{2\gamma^2_i M^2}{\alpha_i}+8\beta_{k}\tau_{k}\epsilon.
$
It follows from \eqref{PDHGtheta} that $\tau_k/\eta_k\leq \tau_{k+1}/\eta_{k+1}$. Combining this with
\eqref{PDHGka}, $\tau_k/\gamma_{k}\leq \tau_{k+1}/\gamma_{k+1}$ and $\beta_{k}\tau_{k}\leq\beta_{k+1}\tau_{k+1}$, we obtain $R_{k+1}\leq R_{N+1}$.

Next, we bound $\bbe[\norm{\tilde{y}^{N+1}}^2]$ by induction. For the case of $N=1$,
	it readily follows from  $\tilde{y}^{1}=0$, $\beta_k\geq 1$, the update $\tilde{y}^{k+1}=(1-\beta_k^{-1})\tilde{y}^{k}+\beta_k^{-1}y^{k+1}$ in \eqref{PDHGy1} and \eqref{zbeq2} that
$\bbe[\enorm{\tilde{y}^{2}}^2]\leq \bbe[\enorm{y^{2}}^2]\leq R_2\leq R_{N+1}$. Suppose that $\bbe[\enorm{\tilde{y}^{k+1}}^2] \leq R_{N+1}$ holds for $k<N$. It then follows from \eqref{PDHGy1} and \eqref{zbeq2} (placing $k=N$) that
$\bbe[\enorm{\tilde{y}^{N+1}}^2] \leq R_{N+1}$, which proves result (b).
\end{proof}

Obtaining separate bounds for the optimality gap and the primal infeasibility is often essential in constrained optimization analysis.
The following lemma extends the results for equality constraints in~\citep{xu2017accelerated} to the case of inequality constraints, and we use it to derive these bounds.
\begin{lemma}\label{rate2}
	Consider the auxiliary convex optimization problem
	\begin{equation}\label{aop}
	\min_{x \in X} f(x) \quad \text{subject to} \quad Ax - b \leq 0,
	\end{equation}
	where $f: X \rightarrow \mathbb{R}$ is a convex function, $A \in \mathbb{R}^{m \times n}$, and $b \in \mathbb{R}^m$.
	Let $\left(x^*, y^*\right)$ be a Karush--Kuhn--Tucker (KKT) point of \eqref{aop}.
	For $\bar{x}\in X$, $\rho>\|y^*\|$, and $\epsilon\geq0$, suppose that
	\begin{equation*}
		f(\bar{x})-f(x^*)+\rho\enorm{[A \bar{x}-b]_{+}} \leq \epsilon.
	\end{equation*}
	Then it holds that
		\begin{equation*}
		\enorm{[A \bar{x}-b]_{+}} \leq \frac{\epsilon}{\rho-\enorm{y^*}} \quad \text{and} \quad -\frac{\enorm{y^*} \epsilon}{\rho-\enorm{y^*}} \leq f(\bar{x})-f(x^*) \leq \epsilon.
		\end{equation*}
\end{lemma}
\begin{proof}
    It suffices to replace $\norm{A\bar{x}-b}$ with $\norm{[A\bar{x}-b]_+}$ and the unrestricted dual space with $\bbr^m_+$ in the proofs of Lemmas 2.2-2.3 in~\citep{xu2017accelerated}. The desired result then follows directly.
\end{proof}

Using the preceding lemmas, we derive explicit convergence rates for both the functional optimality gap and the primal infeasibility.

\begin{theorem}\label{optgapth}
	Suppose Assumption~\ref{ass1} holds, and let $(x^*,y^*)$ be a primal--dual optimal pair of~\eqref{consistent}. Let $\tilde z^{k+1}=(\tilde x^{k+1},\tilde y^{k+1})$ be generated by Algorithm~\ref{alg:A-ISPD}. Under the parameter conditions of Theorem~\ref{thzb}, the following bounds hold:
\begin{equation*}
\begin{aligned}
&\bbe[|f(\tilde{x}^{N+1},c)+v(\tilde{x}^{N+1})-f(x^*,c)-v(x^*)|]\\
&\leq \frac{\Omega^2_X}{2\beta_{N}\eta_N}+\frac{\alpha_1 \Omega^2_X}{2\beta_N\gamma_{N}}
+\frac{1+2\left\|y^*\right\|^2}{\beta_{N}\tau_N}\\
&+ \frac{1}{\beta_{N}\gamma_{N}}\sum_{k=1}^{N}\frac{2\gamma_k M^2}{\eta_k^{-1}-L_f\beta_k^{-1}-2\norm{A}^2\tau_{k}}
+\frac{1}{\beta_{N}\gamma_{N}}\sum_{k=1}^{N}\frac{2\gamma^2_k M^2}{\alpha_k}+\epsilon,
\end{aligned}
\end{equation*}
and
\begin{equation*}
\begin{aligned}
&\bbe[\enorm{[A \tilde{x}^{N+1}-b+B u]_+}]\\
&\leq \frac{\Omega^2_X}{2\beta_{N}\eta_N}+\frac{\alpha_1 \Omega^2_X}{2\beta_N\gamma_{N}}
+\frac{1+2\left\|y^*\right\|^2}{\beta_{N}\tau_N}\\
&+ \frac{1}{\beta_{N}\gamma_{N}}\sum_{k=1}^{N}\frac{2\gamma_k M^2}{\eta_k^{-1}-L_f\beta_k^{-1}-2\norm{A}^2\tau_{k}}
+\frac{1}{\beta_{N}\gamma_{N}}\sum_{k=1}^{N}\frac{2\gamma^2_k M^2}{\alpha_k}+\epsilon.
\end{aligned}
\end{equation*}
\end{theorem}
\begin{proof}
	Since $-\enp{Ax^*+Bu-b,\tilde{y}^{N+1}}\geq 0$, we have 
	\begin{equation*}
	Q(\tilde{z}^{N+1}; x^*,y)\geq \enp{A \tilde{x}^{N+1}-b+B u,y}+f(\tilde{x}^{N+1},c)+v(\tilde{x}^{N+1})-f(x^*,c)-v(x^*).
	\end{equation*}
	Letting $x=x^*$ in~\eqref{thimq1}, with each $\Lambda_k$ evaluated at $x^*$, we obtain
	\begin{equation*}
		\begin{aligned}
		&\quad\;f(\tilde{x}^{N+1},c)+v(\tilde{x}^{N+1})-f(x^*,c)-v(x^*)+\enp{A \tilde{x}^{N+1}-b+B u,y}\\
		&\leq \frac{\Omega_X^2}{2\beta_N\eta_N}
		+\frac{1}{2\beta_{N}\tau_N}\norm{y}^2
		+\frac{1}{\beta_{N}\gamma_{N}}\sum_{k=1}^{N}\Lambda_k(x)+\epsilon.
		\end{aligned}
	\end{equation*}
	Let $r^{N+1}:=[A \tilde{x}^{N+1}-b+B u]_+$. If $r^{N+1}=0$, set $y=0$; otherwise, take $y=\rho r^{N+1}/\norm{r^{N+1}}$ with $\rho>0$. Substituting this choice into the above inequality, we have
		\begin{equation*}
	\begin{aligned}
	&f(\tilde{x}^{N+1},c)+v(\tilde{x}^{N+1})-f(x^*,c)-v(x^*)+\rho\norm{[A \tilde{x}^{N+1}-b+B u]_+}\\
	&\leq \frac{\Omega_X^2}{2\beta_N\eta_N}
	+\frac{\rho^2}{2\beta_{N}\tau_N}
	+\frac{1}{\beta_{N}\gamma_{N}}\sum_{k=1}^{N}\Lambda_k(x)+\epsilon.
	\end{aligned}
	\end{equation*}	
	Applying Lemma~\ref{rate2} with $\rho=\max\{1+\enorm{y^*},2\enorm{y^*}\}$, and using $\max\{(1+a)^2,4a^2\}\leq2+4a^2$ for $a\geq0$, we have
	\begin{equation*}
		\begin{aligned}
			&|f(\tilde{x}^{N+1},c)+v(\tilde{x}^{N+1})-f(x^*,c)-v(x^*)|\\
            &\leq \frac{\Omega_X^2}{2\beta_N\eta_N}
		+\frac{\max\{(1+\enorm{y^*})^2, 4\enorm{y^*}^2\}}{2\beta_{N}\tau_N}
		+\frac{1}{\beta_{N}\gamma_{N}}\sum_{k=1}^{N}\Lambda_k(x)+\epsilon\\
        &\leq \frac{\Omega_X^2}{2\beta_N\eta_N}
		+\frac{1+2\enorm{y^*}^2}{\beta_{N}\tau_N}
		+\frac{1}{\beta_{N}\gamma_{N}}\sum_{k=1}^{N}\Lambda_k(x)+\epsilon
		\end{aligned}
	\end{equation*}
	and 
	\begin{equation*}
	\begin{aligned}
	    \enorm{[A \tilde{x}^{N+1}-b+B u]_+} \leq \frac{\Omega_X^2}{2\beta_N\eta_N}
	+\frac{1+2\enorm{y^*}^2}{\beta_{N}\tau_N}
	+\frac{1}{\beta_{N}\gamma_{N}}\sum_{k=1}^{N}\Lambda_k(x)+\epsilon.
	\end{aligned}
	\end{equation*}
	Taking a full expectation of the above two inequalities and applying~\eqref{lammmdd} with $x=x^*$, we have the desired result.
\end{proof}

In contrast to previous work~\citep[Lemma 8]{lan2021dynamic}, which employs the perturbed gap function to simultaneously measure the functional optimality gap and primal infeasibility, we use Lemma~\ref{rate2} to establish complexity bounds for both quantities without introducing additional perturbation variables.

\subsubsection{Stage-Specific Parameter Settings}
Under Assumption~\ref{ass1}, we specify three parameter choices for different stage roles. Throughout these corollaries, assume $\|A\|>0$; the intermediate- and first-stage choices also require $M>0$, which may be any positive upper bound satisfying Assumption~\ref{ass1}(ii). The weights $\gamma_k$ and $\alpha_k$ are used only in the analysis. The displayed formulas for $\beta_k$, $\theta_k$, and $\gamma_k$ also define these quantities at $k=N+1$, as needed for the recurrence and final extrapolation. Proofs are deferred to Appendix~\ref{app:proofs-corollaries}.

\begin{corollary}\label{cor:last-stage}
	Consider Algorithm~\ref{alg:A-ISPD} with a given $N \geq 1$, and suppose the parameters are chosen as
	\begin{equation}\label{paraco1}
	\begin{aligned}
	\gamma_{k}=k,\beta_{k}=\frac{k+1}{2},\theta_k=\frac{k-1}{k},
	\eta_k=\frac{k}{2(L_f+N\norm{A})+1},\tau_k=\frac{k}{N\norm{A}},\alpha_k=1
	\end{aligned} 
	\end{equation}
	 for any $1\leq k\leq N$. Then the following results hold:
	\begin{itemize}
		\item[$\operatorname{(a)}$] 	
		\begin{equation*}
		\begin{aligned}
		\bbe[\operatorname{gap}(\tilde{z}^{N+1})]
	\leq \frac{2\Omega^2_XL_f}{N^2}+\frac{2\Omega_X^2\enorm{A}}{N}+\frac{2\Omega^2_X}{N^2}+\frac{\enorm{A}\enorm{y^*}^2}{N}
     +\frac{4M^2(2N+1)}{3}+\epsilon
		\end{aligned}
		\end{equation*}

		\item[$\operatorname{(b)}$]
		\begin{equation*}
		\begin{aligned}
		\bbe[\enorm{\tilde{y}^{N+1}}^2]
		\leq 8\Omega_X^2+\frac{8\Omega_X^2L_f}{N\enorm{A}}
	+\frac{8\Omega_X^2}{N\enorm{A}}+6\enorm{y^*}^2
	 +\frac{16M^2(N+1)(2N+1)}{3\enorm{A}}+\frac{4(N+1)\epsilon}{\enorm{A}}
		\end{aligned}
		\end{equation*}
		
		\item[$\operatorname{(c)}$] 	
		\begin{equation*}
\begin{aligned}
		&\bbe[|f(\tilde{x}^{N+1},c)+v(\tilde{x}^{N+1})-f(x^*,c)-v(x^*)|]\\
&\leq \frac{2\Omega^2_XL_f}{N^2}+\frac{2\Omega_X^2\enorm{A}}{N}+\frac{2\Omega^2_X}{N^2}
+\frac{\enorm{A}(2+4\enorm{y^*}^2)}{N} +\frac{4M^2(2N+1)}{3}+\epsilon.
\end{aligned}
\end{equation*}
		and
\begin{equation*}
\begin{aligned}
	& \bbe[\enorm{[A \tilde{x}^{N+1}-b+B u]_+}]\\
    &\leq \frac{2\Omega^2_XL_f}{N^2}+\frac{2\Omega_X^2\enorm{A}}{N}+\frac{2\Omega^2_X}{N^2}
+\frac{\enorm{A}(2+4\enorm{y^*}^2)}{N} +\frac{4M^2(2N+1)}{3}+\epsilon.
\end{aligned}
\end{equation*}
	\end{itemize}
\end{corollary}

When $M=0$ and $\epsilon=0$, Corollary~\ref{cor:last-stage} gives an $O(1/N)$ gap while keeping the expected squared dual norm $O(1)$, making it suitable for the deterministic terminal stage. With stochastic variability or oracle inexactness, however, its gap and dual bounds contain terms that grow with $N$, so a different setting is needed at intermediate stages.

\begin{corollary}\label{cor:middle-stage}
	Consider Algorithm~\ref{alg:A-ISPD} with a given $N \geq 1$, and suppose the parameters are chosen as
	\begin{equation}\label{paraco2}
	\begin{aligned}
	\gamma_{k}=k,\beta_{k}=\frac{k+1}{2},\theta_k=\frac{k-1}{k},
	\eta_k=\frac{k}{2(L_f+N\enorm{A})+k^{\frac{3}{2}}M},\tau_k=\frac{k}{N^{\frac{3}{2}} \enorm{A}}, \alpha_k=N^{\frac{3}{2}}M
	\end{aligned}
	\end{equation}
    for any $1\leq k\leq N$. 
		Then the following results hold:
	\begin{itemize}
		\item[$\operatorname{(a)}$] 
		\begin{equation*}
		\begin{aligned}
		\bbe[\operatorname{gap}(\tilde{z}^{N+1})]
	\leq \frac{2\Omega^2_XL_f}{N^2}+\frac{2\Omega_X^2\enorm{A}}{N}+\frac{2M\Omega^2_X}{\sqrt{N}}+\frac{\enorm{A}\enorm{y^*}^2}{\sqrt{N}}
	 +\frac{12M}{3\sqrt{N}}+\frac{2M}{3N^{\frac{3}{2}}}+\epsilon,
		\end{aligned}
		\end{equation*}
		
		\item[$\operatorname{(b)}$] 
		\begin{equation*}
		\begin{aligned}
		\bbe[\enorm{\tilde{y}^{N+1}}^2]
		\leq \frac{8L_f\Omega_X^2}{N^{\frac{3}{2}}\norm{A}}+\frac{8\Omega_X^2}{\sqrt{N}}+\frac{8M\Omega_X^2}{\enorm{A}}+6\enorm{y^*}^2
		+\frac{64M}{3\enorm{A}}+\frac{16M}{\enorm{A}}+\frac{8\sqrt{N}\epsilon}{\enorm{A}},
		\end{aligned}
		\end{equation*}
		
		\item[$\operatorname{(c)}$] 	
		\begin{equation*}
		\begin{aligned}
		& \bbe[|f(\tilde{x}^{N+1},c)+v(\tilde{x}^{N+1})-f(x^*,c)-v(x^*)|]\\
	&\leq \frac{2\Omega^2_XL_f}{N^2}+\frac{2\Omega_X^2\enorm{A}}{N}+\frac{2M\Omega^2_X}{\sqrt{N}}
    +\frac{\enorm{A}(2+4\enorm{y^*}^2)}{\sqrt{N}} +\frac{12M}{3\sqrt{N}}+\frac{2M}{3N^{\frac{3}{2}}}+\epsilon
		\end{aligned}
		\end{equation*}
		and 
		\begin{equation*}
		\begin{aligned}
		& \bbe[\enorm{[A \tilde{x}^{N+1}-b+B u]_+}]\\
	&\leq \frac{2\Omega^2_XL_f}{N^2}+\frac{2\Omega_X^2\enorm{A}}{N}+\frac{2M\Omega^2_X}{\sqrt{N}}
    +\frac{\enorm{A}(2+4\enorm{y^*}^2)}{\sqrt{N}}  +\frac{12M}{3\sqrt{N}}+\frac{2M}{3N^{\frac{3}{2}}}+\epsilon.
		\end{aligned}
		\end{equation*}
	\end{itemize}
\end{corollary}

With $N=\mathcal{O}(1/\epsilon^2)$, the setting in \eqref{paraco2} controls the expected dual moment and is therefore used for recursive subgradient construction at intermediate stages.

\begin{corollary}\label{cor:first-stage}
	Consider Algorithm~\ref{alg:A-ISPD} with a given $N\geq 1$ and parameters
	\begin{equation}\label{paraco3}
	\begin{aligned}
	\gamma_{k}=k, \beta_{k}=\frac{k+1}{2},\theta_k=\frac{k-1}{k},
	\eta_k=\frac{k}{2(L_f+N\enorm{A})+k^{\frac{3}{2}}M},\tau_k=\frac{k}{N\enorm{A}},\alpha_k=N^{\frac{3}{2}}M
	\end{aligned}
	\end{equation}
    for $1\leq k\leq N$. 
		Then the following results hold:
	\begin{itemize}
		\item[$\operatorname{(a)}$] 
		\begin{equation*}
		\begin{aligned}
		\bbe[\operatorname{gap}(\tilde{z}^{N+1})]
		\leq \frac{2\Omega^2_XL_f}{N^2}+\frac{2\Omega^2_X\enorm{A}}{N}+\frac{2M\Omega^2_X}{\sqrt{N}}+\frac{\enorm{A}\enorm{y^*}^2}{N}
		+\frac{12M}{3\sqrt{N}}+\frac{2M}{3N^{\frac{3}{2}}}+\epsilon,
		\end{aligned}
		\end{equation*}
		
		\item[$\operatorname{(b)}$]
		\begin{equation*}
		\begin{aligned}
		\bbe[\enorm{\tilde{y}^{N+1}}^2]
		&\leq 8\Omega^2_X+\frac{8\Omega^2_X L_f}{N\enorm{A}}+\frac{8\Omega_X^2M\sqrt{N}}{\enorm{A}}+6\enorm{y^*}^2\\
		&\quad + \frac{32M(N+1)}{3\sqrt{N}\enorm{A}}
		+\frac{8M(N+1)(2N+1)}{3N^{\frac{3}{2}}\enorm{A}}+\frac{4(N+1)\epsilon}{\enorm{A}}
		\end{aligned}
		\end{equation*}
		
		\item[$\operatorname{(c)}$]
		\begin{equation*}
		\begin{aligned}
		&\bbe[|f(\tilde{x}^{N+1},c)+v(\tilde{x}^{N+1})-f(x^*,c)-v(x^*)|]\\
	&\leq \frac{2\Omega^2_XL_f}{N^2}+\frac{2\Omega^2_X\enorm{A}}{N}+\frac{2M\Omega^2_X}{\sqrt{N}}
    +\frac{\enorm{A}(2+4\enorm{y^*}^2)}{N} +\frac{12M}{3\sqrt{N}}+\frac{2M}{3N^{\frac{3}{2}}}+\epsilon,
		\end{aligned}
		\end{equation*}
		and
		\begin{equation*}
		\begin{aligned}
		&\bbe[\enorm{[A \tilde{x}^{N+1}-b+B u]_+}]\\
		&\leq \frac{2\Omega^2_XL_f}{N^2}+\frac{2\Omega^2_X\enorm{A}}{N}+\frac{2M\Omega^2_X}{\sqrt{N}}
        +\frac{\enorm{A}(2+4\enorm{y^*}^2)}{N} +\frac{12M}{3\sqrt{N}}+\frac{2M}{3N^{\frac{3}{2}}}+\epsilon.
		\end{aligned}
		\end{equation*}
	\end{itemize}
\end{corollary}

\begin{remark}
For the optimality and feasibility bounds, Corollary~\ref{cor:first-stage} improves the multiplier-dependent term in Corollary~\ref{cor:middle-stage} by replacing its denominator $\sqrt{N}$ with $N$, while the remaining terms keep the same order. Its gap, optimality, and feasibility bounds also avoid the growing stochastic term in Corollary~\ref{cor:last-stage}. Its dual second-moment bound can still grow with $N$, so the intermediate-stage setting is used when a subgradient must be passed to a preceding stage.
\end{remark}

\subsection{Convergence Analysis of A-DSPD}\label{subsec:three-stage-convergence}

We now apply the stage-specific \iaspd{} bounds to the three nested loops of Algorithm~\ref{A-DSPD}. The innermost and middle loops construct stochastic $\epsilon$-subgradients of $v_3$ and $v_2$, respectively, and the outer loop uses them to bound first-stage optimality and infeasibility. Along this recursion, we verify \eqref{Gbound} and use the terminal-, intermediate-, and first-stage settings according to their roles described above.

Fix a target accuracy $0<\bar{\epsilon}\leq1$. We round each displayed iteration budget up to an integer, with a minimum of one. Let $N_t$ be the number of updates per stage-$t$ call, with indices $i$, $j$, and $k$ for stages $1$, $2$, and $3$, respectively. A superscript $t$ on an algorithmic parameter identifies its stage; it is suppressed in the stage-specific proofs. Each inner call starts from its prescribed primal point in $X_t$ and zero dual variables. When an outer-loop index is suppressed, the corresponding history, predecessor decision, and sampled matrices remain those of the current call. We now present the framework of the A-DSPD algorithm for the three-stage problem~\eqref{three} in Algorithm~\ref{A-DSPD}.
At each nonterminal value-function oracle call, a fresh scenario is sampled from the distribution conditional on the current scenario history, independently of preceding oracle draws given that history. The full history is retained when constructing the next-stage value-function oracle.
\begin{algorithm}[htbp]  
	\caption{A-DSPD Algorithm for Three-Stage Problem \eqref{three}}\label{A-DSPD}
	\begin{algorithmic}[0]
		\Require initial points $(x_1^1,x_2^1,x_3^1)=(\tilde{x}_1^1,\tilde{x}_2^1,\tilde{x}_3^1)=(\bar{x}^1_1,\bar{x}^1_2,\bar{x}^1_3)$, $(y^1_1,y^1_2,y^1_3)=(\tilde{y}^1_1,\tilde{y}^1_2,\tilde{y}^1_3)=0$
		\For{$i=1,\cdots,N_1$}
		\State Sample $\xi_2^i$ from the conditional distribution of $\xi_2$ given $\xi_{[1]}$.
        \State Set $\xi_{[2]}^i=(\xi_{[1]},\xi_2^i)$.
		\For{$j=1,\cdots,N_2$}
		\State Sample $\xi_3^j$ from the conditional distribution of $\xi_3$ given $\xi_{[2]}^i$.
        \State Set $\xi_{[3]}^j=(\xi_{[2]}^i,\xi_3^j)$.
		\For{$k=1,\cdots,N_3$}
		\State $(\tilde{x}^{k+1}_3,\tilde{y}^{k+1}_3,\bar{x}^{k+1}_3,x_3^{k+1},y_3^{k+1})=\operatorname{ASPDT}(\tilde{x}^k_3,\tilde{y}^k_3,x^{k}_3,y^{k}_3,\bar{x}^k_3,0,x^j_2,\xi^j_3,X_3,\tau^3_k,\eta^3_k)$
		\EndFor
		\State $G_2(x^j_2)=B_3^\top\tilde{y}^{N_3+1}_3$
		\State $(\tilde{x}^{j+1}_2,\tilde{y}^{j+1}_2,\bar{x}^{j+1}_2,x^{j+1}_2,y_2^{j+1})=\operatorname{ASPDT}(\tilde{x}^j_2,\tilde{y}^j_2,x^{j}_2,y^{j}_2,\bar{x}^j_2,G_2(x^j_2),x^i_1,\xi^i_2,X_2,\tau^2_j,\eta^2_j)$
		\EndFor
		\State $G_1(x^i_1)=B_2^\top\tilde{y}^{N_2+1}_2$
		\State $(\tilde{x}^{i+1}_1,\tilde{y}^{i+1}_1,\bar{x}^{i+1}_1,x^{i+1}_1,y_1^{i+1})=\operatorname{ASPDT}(\tilde{x}^i_1,\tilde{y}^i_1,x^{i}_1,y^{i}_1,\bar{x}^i_1,G_1(x^i_1),0,\xi_1,X_1,\tau^1_i,\eta^1_i)$
		\EndFor
		\Ensure $(\tilde{x}_1^{N_1+1}, \tilde{y}_1^{N_1+1})$
	\end{algorithmic}  
\end{algorithm}  
The A-DSPD algorithm is organized as three nested loops. In the innermost (third) loop, $N_3$ iterations of the ASPDT are performed to compute a stochastic $\frac{\bar{\epsilon}}{3}$-subgradient of the value function $v_3$. This stochastic subgradient is then used in the second-stage (middle) loop, where $N_2$ iterations of ASPDT, together with the third-stage calls, yield a stochastic $\frac{2\bar{\epsilon}}{3}$-subgradient of the value function $v_2$. Finally, the outermost (first) loop utilizes the subgradient obtained from the second stage and performs $N_1$ iterations to approximate the optimal solution and control the primal infeasibility of the first-stage problem within $\bar{\epsilon}$.

To support the subsequent analysis, we introduce the following assumption.
\begin{assu}\label{ass2}
	The following conditions are assumed to hold:
	\begin{itemize}
		\item[$\operatorname{(i)}$] For each $t=1,2,3$, the function $f_t(\cdot,c_t)$ is smooth and convex, and there exist deterministic constants $L_{f_t},G_{f_t}\geq 0$ such that, almost surely,
		\begin{equation}
		f_t(x',c_t)\leq f_t(x,c_t)+\enp{\nabla f_t(x,c_t),x'-x}+\frac{L_{f_t}}{2}\enorm{x'-x}^2,\quad \enorm{\nabla f_t(\tilde{x},c_t)}\leq G_{f_t},
		\end{equation}
        for all $x,x',\tilde{x}\in X_t$, uniformly over the scenario data. Here $L_{f_t}$ is a smoothness constant, whereas $G_{f_t}$ bounds the gradient norm.
		\item[$\operatorname{(ii)}$] The stochastic linear mappings $B_t$ are uniformly bounded, i.e., there exist constants $L_{B_t} > 0$ such that
		\begin{equation}\label{bound}
             \enorm{B_t} \leq L_{B_t}, \quad t=2,3.
		\end{equation}
		\item[$\operatorname{(iii)}$] For each $t=1,2,3$, there exist deterministic constants $\underline{a}_t$ and $\overline{a}_t$ such that, almost surely,
        \begin{equation}\label{Abound-three}
            0<\underline{a}_t\leq\enorm{A_t}\leq\overline{a}_t<\infty.
        \end{equation}
        These bounds control the operator norms in the step sizes and the estimates for the dual iterates.
        \item[$\operatorname{(iv)}$] For each $t=1,2,3$, there exists a deterministic constant $Y_t<\infty$ and a jointly measurable selection of optimal dual solutions $y_t^*(u,\xi_{[t]})\in\bbr^{m_t}_+$ such that, outside a common null set of scenario histories,
        \begin{equation*}
            \enorm{y_t^*(u,\xi_{[t]})}\leq Y_t
            \quad\text{for every }u\in X_{t-1}\quad(t\geq2).
        \end{equation*}
        At $t=1$, this means that a fixed optimal multiplier $y_1^*$ satisfies $\enorm{y_1^*}\leq Y_1$. The constants $Y_t$ are uniform over predecessor decisions and scenario histories and do not depend on the requested accuracy or the iteration counts. Only one such selection is required; the entire dual optimal set need not be bounded. Sufficient conditions are discussed in Subsection~\ref{subsec:bounded-multipliers}.
    \end{itemize}
\end{assu}
Under Assumption~\ref{ass2}, the estimates below provide a common constant $M_t$ for the Lipschitz continuity of $v_t$ and the second moment of its stochastic subgradient oracle. The multiplier bound $Y_t$ is an input to these estimates; it is not inferred from Lipschitz continuity on $X_t$.
The constants $M_t$ are chosen deterministically and uniformly over predecessor decisions and scenario histories, independently of the target accuracy. The same convention applies in the general $T$-stage analysis.

Before deriving stagewise bounds, we recall the following elementary lemma.
\begin{lemma}\label{subg}
	Let $f(x)$ be a convex function defined on a convex set $X$. Suppose there exists a constant $\bar{M}>0$ such that, for any $x \in X$, there is a subgradient $f'(x) \in \partial f(x)$ satisfying $\norm{f'(x)} \leq \bar{M}$. Then, for all $x, y \in X$, it holds that
	\[
	|f(x) - f(y)| \leq \bar{M}\norm{x - y}.
	\]
\end{lemma}
\begin{proof}
	By the definition of subgradient and the Cauchy--Schwarz inequality, we obtain
	\begin{equation*}
		\begin{aligned}
		f(x)-f(y)&\leq \enp{f'(x),x-y}\leq \norm{f'(x)}\norm{x-y}\leq \bar{M}\norm{x-y},\\
		f(x)-f(y)&\geq \enp{f'(y),x-y}\geq -\norm{f'(y)}\norm{y-x}\geq -\bar{M}\norm{x-y}.\\
		\end{aligned}
	\end{equation*}
	Combining the two inequalities yields $|f(x) - f(y)| \leq \bar{M} \norm{x-y}$.
\end{proof}

For fixed scenario data, an optimal multiplier at $u$ gives a supporting subgradient of the stage value function: by strong duality and dual feasibility,
\begin{equation*}
V_t(u',\xi_{[t]})\geq V_t(u,\xi_{[t]})
+\enp{B_t^\top y_t^*(u,\xi_{[t]}),u'-u}.
\end{equation*}
Applying this inequality at both $u$ and $u'$ gives a Lipschitz bound $L_{B_t}Y_t$ for $V_t$, and hence for its conditional expectation $v_t$. The second-moment estimates below may require a larger constant; taking the maximum preserves the use of the same $M_t$ for both purposes.

In the stagewise estimates below, $x_t^*$ denotes an optimal primal solution for the current predecessor decision and scenario history, and $y_t^*$ is the multiplier selected in Assumption~\ref{ass2}(iv). We write $\tilde z_t^{k+1}:=(\tilde x_t^{k+1},\tilde y_t^{k+1})$ and use $\operatorname{gap}$ for the corresponding instance of~\eqref{gapdef}. All conditional estimates are uniform in the fixed predecessor decision; for an adaptive query, they apply conditional on the information available before that oracle call. Any oracle bound $M_t$ below may be enlarged to a positive constant when used in~\eqref{paraco2} or~\eqref{paraco3}. We begin with the third (innermost) loop.
\begin{lemma}\label{A-DSAle1}
	Suppose Assumption~\ref{ass2} holds and let the sequence $\{\tilde{x}_3^{k}, \tilde{y}_3^{k}\}_{k=1}^{N_3+1}$ be generated by Algorithm~\ref{A-DSPD}. Suppose further that the
        parameters $\alpha_k$, $\beta_k$, $\theta_k$, $\eta_k$, and $\tau_k$ are chosen according to \eqref{paraco1} (with $M=0$ and $\epsilon=0$), and 
        the number of iterations $N_3\geq 1$ is set as
	\begin{equation}\label{N_3}
		N_3=\max\left\{\frac{2\sqrt{6L_{f_3}}\Omega_{3}}{\sqrt{\bar{\epsilon}}},
		\frac{24\Omega_3^2\overline{a}_3}{\bar{\epsilon}},\frac{2\sqrt{6}\Omega_{3}}{\sqrt{\bar{\epsilon}}},
		\frac{12\overline{a}_3Y_3^2}{\bar{\epsilon}},
		\frac{12\overline{a}_3}{\bar{\epsilon}} \cdot \left(2
			+4Y_3^2\right)\right\},
	\end{equation}
	then $B_3^\top \tilde{y}_3^{N_3+1}$ is a stochastic $\frac{\bar{\epsilon}}{3}$-subgradient of the value function $v_3$ at the point $x_2^j$. Furthermore, for any given $\xi_{[2]}$, there exists a constant $M_3$ such that $|v_3(x_2,\xi_{[2]})-v_3(x_2',\xi_{[2]})|\leq M_3\enorm{x_2-x_2'}$ for all $x_2, x_2' \in X_2$, with
		\begin{equation*}
		\bbe[\enorm{B_3^\top\tilde{y}^{N_3+1}_3}^2\vert\xi_{[2]}]\leq M^2_3.
		\end{equation*}
	In addition, the following bounds hold:
	\begin{equation}\label{result1}
	\begin{aligned}
	&|f_3(\tilde{x}_3^{N_3+1},c_3)-f_3(x^*_3,c_3)|\leq \frac{\bar{\epsilon}}{3}\\
	&\enorm{[A_3 \tilde{x}_3^{N_3+1}-b_3+B_3 x^j_2]_+}\leq \frac{\bar{\epsilon}}{3}.
	\end{aligned}
	\end{equation}
\end{lemma}
\begin{proof}
		We apply \iaspd{} to the third-stage problem \eqref{third} in the innermost loop of Algorithm~\ref{A-DSPD}. Since there are no subsequent stages, it is a deterministic saddle-point problem for a given $x^j_2$ and $\xi^j_{[3]}$. As a result, $v=0$, $M=0$, and $\epsilon=0$. Let $u=x^j_2$, $\xi=\xi^j_3$, $A=A_3$,
	$B=B_3$, $L_f=L_{f_3}$, $\Omega_X=\Omega_{3}$, $\tilde{x}^{N+1}=\tilde{x}^{N_3+1}$, $x^*=x^*_3$, {$\tilde{y}^{N+1}=\tilde{y}^{N_3+1}$}, $y^*=y^*_3$ and $y^1=y^1_3=0$. 
	    By Assumption~\ref{ass2}(iv), choose the optimal multiplier $y_3^*$ from the stated measurable selection, so that
\begin{equation}\label{yy3b}
    \enorm{y^*_3}\leq Y_3.
\end{equation} 
	Immediately, in view of Corollary \ref{cor:last-stage}(a) with $M=0$ and $\epsilon=0$, the definition of 
	$N_3$ and \eqref{yy3b}, we have
	\begin{equation*}
		\begin{aligned}
		\operatorname{gap}(\tilde{z}^{N_3+1}_3)
		&\leq \frac{2\Omega^2_3L_{f_3}}{N_3^2}+\frac{2\Omega_3^2\enorm{A_3}}{N_3}+\frac{2\Omega^2_3}{N_3^2}+\frac{\enorm{A_3}\enorm{y_3^*}^2}{N_3}\\
		&\leq \frac{2\Omega^2_3L_{f_3}}{N_3^2}+\frac{2\Omega_3^2\overline{a}_3}{N_3}+\frac{2\Omega^2_3}{N_3^2}+\frac{\overline{a}_3Y_3^2}{N_3}\\
		&\leq \frac{\bar{\epsilon}}{12}+\frac{\bar{\epsilon}}{12}+\frac{\bar{\epsilon}}{12}+\frac{\bar{\epsilon}}{12}
		\leq \frac{\bar{\epsilon}}{3},
		\end{aligned}
	\end{equation*}
	which, by Lemma~\ref{subgradient}, implies that $B^\top_3 \tilde{y}^{N_3+1}_3$ is a $\bar{\epsilon}/3$-subgradient of the value function $V_3$ at $x_2^{j}$ (i.e., a stochastic $\bar{\epsilon}/3$-subgradient of $v_3$ at $x_2^{j}$).
	Moreover, it follows from Corollary \ref{cor:last-stage}(b) (with $M=0$ and $\epsilon=0$) and \eqref{yy3b} that 
			\begin{equation*}
	\begin{aligned}
	\enorm{\tilde{y}_3^{N_3+1}}^2
	\leq 8\Omega_3^2+\frac{8\Omega_3^2L_{f_3}}{N_3\norm{A_3}}+\frac{8\Omega_3^2}{N_3\enorm{A_3}}+6\norm{y_3^*}^2
	\leq 8\Omega_3^2+\frac{8\Omega_3^2L_{f_3}}{\underline{a}_3}+\frac{8\Omega_3^2}{\underline{a}_3}+6Y_3^2.
	\end{aligned}
	\end{equation*}
	By noting that $\enorm{B_3}\leq L_{B_3}$ from \eqref{bound}, we obtain
	\begin{equation}\label{thirdbound1}
	\bbe[\enorm{B^\top_3\tilde{y}_3^{N_3+1}}^2\vert\xi_{[2]}]
	\leq 2L^2_{B_3}\bbe\left[4\Omega_3^2+\frac{4\Omega_3^2L_{f_3}}{\underline{a}_3}+\frac{4\Omega_3^2}{\underline{a}_3}
	+3Y_3^2\vert\xi_{[2]}\right].
	\end{equation} 
	Lemma~\ref{subg} and \eqref{yy3b} show that, to bound the Lipschitz constant of $v_3$, it suffices to choose $M_3$ so that
	\begin{equation}\label{thirdbound2}
		\bbe[\enorm{B^\top_3y^*_3}\vert\xi_{[2]}]\leq L_{B_3}Y_3 \leq M_3.
	\end{equation}
	Combining \eqref{thirdbound1} and \eqref{thirdbound2} together, we set
	\begin{equation*}
		M_3=L_{B_3}\max\left\{Y_3,\sqrt{8\Omega_3^2+\frac{8\Omega_3^2(L_{f_3}+1)}{\underline{a}_3}+6Y_3^2}\right\}.
	\end{equation*}
	Finally, in view of Corollary \ref{cor:last-stage}(c) (with $M=0$ and $\epsilon=0)$, the definition of $N_3$ and \eqref{yy3b}, we conclude that 
		\begin{equation*}
	\begin{aligned}
	|f_3(\tilde{x}_3^{N_3+1},c_3)-f_3(x^*_3,c_3)|
	&\leq \frac{2\Omega^2_3L_{f_3}}{N_3^2}+\frac{2\Omega_3^2\norm{A_3}}{N_3}+\frac{2\Omega^2_3}{N_3^2}
    +\frac{\enorm{A_3}(2+4\left\|y_3^*\right\|^2)}{N_3}\\
	&\leq \frac{2\Omega^2_3L_{f_3}}{N_3^2}+\frac{2\Omega_3^2\overline{a}_3}{N_3}+\frac{2\Omega^2_3}{N_3^2}
    +\frac{\overline{a}_3\left(2+4Y_3^2\right)}{N_3}\\
	&\leq \frac{\bar{\epsilon}}{12}+\frac{\bar{\epsilon}}{12}+\frac{\bar{\epsilon}}{12}+\frac{\bar{\epsilon}}{12}\leq \frac{\bar{\epsilon}}{3}
	\end{aligned}
	\end{equation*}
	and
	\begin{equation*}
	\begin{aligned}
	\norm{[A_3 \tilde{x}_3^{N_3+1}-b_3+B_3 x^j_2]_+}\leq \frac{\bar{\epsilon}}{3}.
	\end{aligned}
	\end{equation*}
	This completes the proof.
\end{proof}

Lemma~\ref{A-DSAle2} establishes the convergence properties of the middle loop within the A-DSPD algorithm.
\begin{lemma}\label{A-DSAle2}
	Suppose the parameters for the third-stage problem are specified as in Lemma~\ref{A-DSAle1}, and let Assumption~\ref{ass2} hold. Consider the sequence $\{\tilde{x}^{j}_2, \tilde{y}^{j}_2\}_{j=1}^{N_2+1}$ generated by Algorithm~\ref{A-DSPD}. If the parameters $\alpha_j$, $\beta_j$, $\theta_j$, $\eta_j$, and $\tau_j$ are chosen according to \eqref{paraco2}, and the number of iterations $N_2\geq 1$ is set as
	\begin{equation*}
	N_2=\max\left\{
	\begin{aligned}
		&\frac{\sqrt{30L_{f_2}}\Omega_{2}}{\sqrt{\bar{\epsilon}}},
	     \frac{30\Omega_2^2\overline{a}_2}{\bar{\epsilon}},\left(\frac{30M_3\Omega_{2}^2}{\bar{\epsilon}}\right)^2,
	    \left(\frac{15\overline{a}_2Y_2^2}{\bar{\epsilon}}\right)^2,\left(\frac{210M_3}{3\bar{\epsilon}}\right)^2,\\
	     &~~~~~~~~~~~~~~~~~~~~~~~~
	     \left(\frac{15\overline{a}_2}{\bar{\epsilon}} \cdot \left(2+4Y_2^2\right)\right)^2
	\end{aligned}
	\right\},
	\end{equation*}
	then $B^\top_2\tilde{y}^{N_2+1}_2$ is a stochastic $\frac{2\bar{\epsilon}}{3}$-subgradient of the value function $v_2$ at $x^{i}_1$. Furthermore, for a given $\xi_{[1]}$, there exists a constant $M_2$ such that for any $x_1, x_1' \in X_1$, 
	$|v_2(x_1,\xi_{[1]})-v_2(x'_1,\xi_{[1]})| \leq M_2\enorm{x_1-x'_1}$,
	and 
	\[
	\bbe[\enorm{B_2^\top\tilde{y}^{N_2+1}_2}^2\vert\xi_{[1]}] \leq M_2^2.
	\]
	In addition,
   \begin{equation*}
\begin{aligned}
	&\bbe[|f_2(\tilde{x}_2^{N_2+1},c_2)+v_3(\tilde{x}_2^{N_2+1})-f_2(x_2^*,c_2)-v_3(x_2^*)|\vert\xi_{[2]}]
\leq \frac{2\bar{\epsilon}}{3},\\
&\bbe[\enorm{[A_2 \tilde{x}_2^{N_2+1}-b_2+B_2 x^i_1]_+}\vert\xi_{[2]}]\leq \frac{2\bar{\epsilon}}{3}.
\end{aligned}
\end{equation*}
\end{lemma}
\begin{proof}
	We apply \iaspd{} to the second-stage problem \eqref{second} within the second loop of Algorithm~\ref{A-DSPD}.
	In this context, for any given $(x^i_1, \xi^i_{[2]})$, let $v = v_3$, $M = M_3$, $\epsilon = \bar{\epsilon}/3$, $u = x^i_1$, $\xi = \xi^i_2$, $A = A_2$, $B = B_2$, $L_f = L_{f_2}$, $\Omega_X = \Omega_{2}$, $\tilde{x}^{N+1} = \tilde{x}^{N_2+1}$, $x^* = x^*_2$, $\tilde{y}^{N+1} = \tilde{y}^{N_2+1}$, $y^* = y^*_2$, and $y^1 = y^1_2 = 0$.
	By Assumption~\ref{ass2}(iv), choose the selected optimal multiplier $y_2^*$, which satisfies
	\begin{equation}\label{yy2b}
    \enorm{y^*_2}\leq Y_2.
\end{equation}
	In view of Corollary \ref{cor:middle-stage}(a), the definition of $N_2$ and \eqref{yy2b}, we have
		\begin{equation*}
	\begin{aligned}
	\bbe[\operatorname{gap}(\tilde{z}^{N_2+1}_2)\vert\xi_{[2]}]
	&\leq \frac{2\Omega^2_2L_{f_2}}{N_2^2}+\frac{2\Omega_2^2\norm{A_2}}{N_2}+\frac{2M_3\Omega^2_2}{\sqrt{N_2}}
	+\frac{\norm{A_2}\norm{y_2^*}^2}{\sqrt{N_2}}
	 +\frac{12M_3}{3\sqrt{N_2}}+\frac{2M_3}{3N_2^{\frac{3}{2}}}+\frac{\bar{\epsilon}}{3}\\
		&\leq \frac{2\Omega^2_2L_{f_2}}{N_2^2}+\frac{2\Omega_2^2\overline{a}_2}{N_2}+\frac{2M_3\Omega^2_2}{\sqrt{N_2}}
	+\frac{\overline{a}_2Y_2^2}{\sqrt{N_2}} +\frac{14M_3}{3\sqrt{N_2}}+\frac{\bar{\epsilon}}{3}\\
	&\leq \frac{\bar{\epsilon}}{15}+\frac{\bar{\epsilon}}{15}+\frac{\bar{\epsilon}}{15}+\frac{\bar{\epsilon}}{15}
	+\frac{\bar{\epsilon}}{15}+\frac{\bar{\epsilon}}{3}\leq \frac{2\bar{\epsilon}}{3},
	\end{aligned}
	\end{equation*}
    which implies from Lemma \ref{subgradient} that $B^\top_2\tilde{y}_2^{N_2+1}$ is a $2\bar{\epsilon}/3$-subgradient of the value function $V_2$ at $x^{i}_1$ (i.e. a stochastic $2\bar{\epsilon}/3$-subgradient of the value function $v_2$ at $x^{i}_1$).  Moreover, it follows from Corollary~\ref{cor:middle-stage}(b) (with $M=M_3$ and $\epsilon=\bar{\epsilon}/3$) and \eqref{yy2b} that
	\begin{equation}\label{bineq1}
	\begin{aligned}
	\bbe[\norm{\tilde{y}^{N_2+1}_2}^2\vert\xi_{[2]}]
	&\leq \frac{8L_{f_2}\Omega_2^2}{N_2^{\frac{3}{2}}\norm{A_2}}+\frac{8\Omega_2^2}{\sqrt{N_2}}+\frac{8M_3\Omega_2^2}{\norm{A_2}}+6\norm{y_2^*}^2
   +\frac{64M_3}{3\norm{A_2}}+\frac{16M_3}{\norm{A_2}}+\frac{8\sqrt{N_2}\bar{\epsilon}}{\norm{A_2}}\\
	&\leq \frac{8L_{f_2}\Omega_2^2}{\underline{a}_2}+8\Omega_2^2+\frac{8M_3\Omega_2^2}{\underline{a}_2}+6Y_2^2
    +\frac{64M_3}{3\underline{a}_2}+\frac{16M_3}{\underline{a}_2}+\frac{8\sqrt{N_2}\bar{\epsilon}}{\underline{a}_2}.
	\end{aligned}
	\end{equation}
	Let us define
    \begin{equation*}
    \begin{aligned}
    H:=1+\max\bigl\{&1,\sqrt{30L_{f_2}}\Omega_2,30\Omega_2^2\overline{a}_2,
    (30M_3\Omega_2^2)^2,\;(15\overline{a}_2Y_2^2)^2,\\
    &(70M_3)^2,\;[15\overline{a}_2(2+4Y_2^2)]^2\bigr\}.
    \end{aligned}
    \end{equation*}
    For $0<\bar{\epsilon}\leq1$, the rounded budget satisfies $\bar{\epsilon}^2N_2\leq H$.
   Using this bound in \eqref{bineq1} together with $\norm{B_2}\leq L_{B_2}$ in \eqref{bound}, we obtain
	\begin{equation}\label{secondbound1}
	\begin{aligned}
	\bbe[\enorm{B^\top_2\tilde{y}_2^{N_2+1}}^2\vert\xi_{[1]}]
	\leq 2L^2_{B_2}\bbe\left[4\Omega_2^2+\frac{12\Omega_2^2(L_{f_2}+M_3)+56M_3+12\sqrt{H}}{3\underline{a}_2}+3Y_2^2\vert\xi_{[1]}\right].
	\end{aligned}
	\end{equation} 
	Lemma~\ref{subg} and \eqref{yy2b} show that, to bound the Lipschitz constant of $v_2$, it suffices to choose $M_2$ so that
	\begin{equation}\label{secondbound2}
	\bbe[\enorm{B^\top_2y^*_2}\vert\xi_{[1]}]\leq  L_{B_2}Y_2\leq M_2.
	\end{equation}
	Combining \eqref{secondbound1} and \eqref{secondbound2}, we may choose $M_2$ as follows: 
	\begin{small}
			\begin{equation*}
		M_2=L_{B_2}\max\left\{Y_2,\sqrt{8\Omega_2^2+\frac{24\Omega_2^2(L_{f_2}+M_3)+112M_3+24\sqrt{H}}{3\underline{a}_2}+6Y_2^2}\right\}.
		\end{equation*}
	\end{small}
		Finally, Corollary~\ref{cor:middle-stage}(c), the definition of $N_2$, and \eqref{yy2b} imply that
   \begin{equation*}
	\begin{aligned}
	&\bbe[|f_2(\tilde{x}_2^{N_2+1},c_2)+v_3(\tilde{x}_2^{N_2+1})-f_2(x_2^*,c_2)-v_3(x_2^*)|\vert\xi_{[2]}]\\
	&\leq \frac{2\Omega^2_2L_{f_2}}{N_2^2}+\frac{2\Omega_2^2\norm{A_2}}{N_2}+\frac{2M_3\Omega^2_2}{\sqrt{N_2}}
    +\frac{\norm{A_2}(2+4\enorm{y_2^*}^2)}{\sqrt{N_2}}
    +\frac{12M_3}{3\sqrt{N_2}}+\frac{2M_3}{3N_2^{\frac{3}{2}}}+\frac{\bar{\epsilon}}{3}\\
    &\leq \frac{2\Omega^2_2L_{f_2}}{N_2^2}+\frac{2\Omega_2^2\overline{a}_2}{N_2}+\frac{2M_3\Omega^2_2}{\sqrt{N_2}}
    +\frac{\overline{a}_2\left(2+4Y_2^2\right)}{\sqrt{N_2}} +\frac{14M_3}{3\sqrt{N_2}}+\frac{\bar{\epsilon}}{3}\\
    &\leq \frac{\bar{\epsilon}}{15}+\frac{\bar{\epsilon}}{15}+\frac{\bar{\epsilon}}{15}+\frac{\bar{\epsilon}}{15}
    +\frac{\bar{\epsilon}}{15}+\frac{\bar{\epsilon}}{3}\leq \frac{2\bar{\epsilon}}{3},
	\end{aligned}
	\end{equation*}
	and 
	\begin{equation*}
		\bbe[\enorm{[A_2\tilde{x}_2^{N_2+1}-b_2+B_2 x^i_1]_+}\vert\xi_{[2]}]\leq \frac{2\bar{\epsilon}}{3}.
	\end{equation*}
	This completes the proof.
\end{proof}

The following theorem bounds the functional optimality gap and the primal infeasibility of the A-DSPD algorithm for a three-stage stochastic optimization problem.
\begin{theorem}\label{th:ADSPD-three}
	Suppose the parameters for the third and second stage problems in the A-DSPD algorithm are selected according to Lemmas~\ref{A-DSAle1} and~\ref{A-DSAle2}, and that Assumption~\ref{ass2} holds. Let $\{\tilde{x}^{i}_1,\tilde{y}^{i}_1\}_{i=1}^{N_1+1}$ denote the sequence generated by Algorithm~\ref{A-DSPD}. If the parameters $\alpha_i$, $\beta_i$, $\theta_i$, $\eta_i$, and $\tau_i$ are set as specified in~\eqref{paraco3}, and
    the number of iterations $N_1\geq 1$ is set as
	\begin{equation*}
	N_1 = \max\left\{
	\frac{\sqrt{30L_{f_1}}\Omega_{1}}{\sqrt{\bar{\epsilon}}},
	\frac{30\Omega_1^2\overline{a}_1}{\bar{\epsilon}},\left(\frac{30M_2\Omega_1^2}{\bar{\epsilon}}\right)^2
	,\left(\frac{210M_2}{3\bar{\epsilon}}\right)^2,
	\frac{15\overline{a}_1}{\bar{\epsilon}} \cdot \left(2+4Y_1^2\right)
	\right\},
	\end{equation*}
	then the following performance guarantees hold:
	\begin{equation*}
	\begin{aligned}
	&\bbe[|f_1(\tilde{x}_1^{N_1+1},c_1)+v_2(\tilde{x}_1^{N_1+1})-f_1(x_1^*,c_1)-v_2(x_1^*)|\vert\xi_{[1]}] \leq \bar{\epsilon},\\
	&\bbe[\enorm{[A_1 \tilde{x}_1^{N_1+1}-b_1]_+}\vert\xi_{[1]}] \leq \bar{\epsilon}.
	\end{aligned}
	\end{equation*}
\end{theorem}
\begin{proof}
	We apply \iaspd{} to the first-stage problem \eqref{first} in the outer loop of Algorithm~\ref{A-DSPD}. We set $v=v_2$, $M=M_2$ and $\epsilon=2\bar{\epsilon}/3$, $u=0$, $\xi=\xi_1$, $A=A_1$, $B=0$, $L_f=L_{f_1}$, $\Omega_X=\Omega_{1}$, $\tilde{x}^{N+1}=\tilde{x}_1^{N_1+1}$, $x^*=x^*_1$, $\tilde{y}^{N+1}=\tilde{y}_1^{N_1+1}$, $y^*=y^*_1$ and $y^1=y^1_1=0$. By Assumption~\ref{ass2}(iv), choose the selected optimal multiplier $y_1^*$, so that
	\begin{equation}\label{yy1b}
    \enorm{y^*_1}\leq Y_1.
\end{equation}
	In view of Corollary~\ref{cor:first-stage}(c), the definition of $N_1$ and \eqref{yy1b}, we have
	\begin{equation*}
	\begin{aligned}
	&\bbe[|f_1(\tilde{x}_1^{N_1+1},c_1)+v_2(\tilde{x}_1^{N_1+1})-f_1(x_1^*,c_1)-v_2(x_1^*)|\vert\xi_{[1]}]\\
	&\leq \frac{2\Omega^2_1 L_{f_1}}{N_1^2}+\frac{2\Omega^2_1\norm{A_1}}{N_1}+\frac{2M_2\Omega^2_1}{\sqrt{N_1}}
    +\frac{\enorm{A_1}(2+4\enorm{y_1^*}^2)}{N_1} +\frac{12M_2}{3\sqrt{N_1}}+\frac{2M_2}{3N_1^{\frac{3}{2}}}+\frac{2\bar{\epsilon}}{3}\\
	&\leq \frac{2\Omega^2_1 L_{f_1}}{N_1^2}+\frac{2\Omega^2_1\overline{a}_1}{N_1}+\frac{2M_2\Omega^2_1}{\sqrt{N_1}}
	+\frac{\overline{a}_1\left(2+4Y_1^2\right)}{N_1} +\frac{14M_2}{3\sqrt{N_1}}+\frac{2\bar{\epsilon}}{3}\\
	&\leq \frac{\bar{\epsilon}}{15}+\frac{\bar{\epsilon}}{15}+\frac{\bar{\epsilon}}{15}+\frac{\bar{\epsilon}}{15}
	+\frac{\bar{\epsilon}}{15}+\frac{2\bar{\epsilon}}{3}\leq \bar{\epsilon}
	\end{aligned}
	\end{equation*}
	and
	\begin{equation*}
	\begin{aligned}
	\bbe[\enorm{[A_1 \tilde{x}_1^{N_1+1}-b_1]_+}\vert\xi_{[1]}]\leq \bar{\epsilon}.
	\end{aligned}
	\end{equation*}
	This completes the proof.
\end{proof}
Theorem~\ref{th:ADSPD-three} controls first-stage optimality and feasibility, while Lemmas~\ref{A-DSAle1}-\ref{A-DSAle2} primarily construct stochastic $\epsilon$-subgradients. Each descent in the recursion follows one sampled scenario path, and linearizing the stage costs replaces the proximal subproblems required by DSA with gradient-based steps.

\subsection{Sufficient Conditions for Bounded Optimal Multipliers}
\label{subsec:bounded-multipliers}

We give a sufficient condition for the uniform multiplier norm bounds in Assumption~\ref{ass2}(iv) for the three-stage problem; measurability of the optimal selection remains part of that assumption. We set $v_4=0$ and $M_4=0$ and suppress the history argument of $v_{t+1}$ when the scenario history is fixed. The preceding gradient and value-function bounds imply that $f_t(\cdot,c_t)+v_{t+1}$ is $(G_{f_t}+M_{t+1})$-Lipschitz on $X_t$.

\paragraph{Uniform Slater condition.}
Slater's condition is a standard constraint qualification for convex optimization; see \citet[Section~5.2.3]{boyd2004convex}. A uniform positive strict feasibility margin yields a uniform multiplier bound.
\begin{pro}[Slater bound]\label{prop:slater-multipliers}
For each $t=1,2,3$, suppose there exists a deterministic $\kappa_t>0$ such that, outside a common null set of scenario histories, for every $x_{t-1}\in X_{t-1}$ there is some $\bar{x}_t\in X_t$ satisfying
\begin{equation}\label{eq:uniform-slater}
A_t\bar{x}_t+B_tx_{t-1}-b_t\leq-\kappa_t\mathbf{1}.
\end{equation}
At $t=1$, omit the predecessor decision and the term $B_tx_{t-1}$. Under the standing strong-duality and attainment assumptions, every optimal multiplier satisfies
\begin{equation}\label{eq:slater-Y}
\enorm{y_t^*}\leq\enorm{y_t^*}_1
\leq\frac{(G_{f_t}+M_{t+1})\Omega_t}{\kappa_t}.
\end{equation}
\end{pro}
\begin{proof}
For an optimal primal solution $x_t^*$ and an optimal multiplier $y_t^*\geq0$, strong duality and evaluation of the Lagrangian at $\bar{x}_t$ give
\begin{equation*}
\begin{aligned}
f_t(x_t^*,c_t)+v_{t+1}(x_t^*)
&\leq f_t(\bar{x}_t,c_t)+v_{t+1}(\bar{x}_t)
+\enp{A_t\bar{x}_t+B_tx_{t-1}-b_t,y_t^*}\\
&\leq f_t(\bar{x}_t,c_t)+v_{t+1}(\bar{x}_t)
-\kappa_t\enorm{y_t^*}_1.
\end{aligned}
\end{equation*}
Rearranging and applying Lipschitz continuity yields
\begin{equation*}
\begin{aligned}
\kappa_t\enorm{y_t^*}_1
&\leq f_t(\bar{x}_t,c_t)-f_t(x_t^*,c_t)
+v_{t+1}(\bar{x}_t)-v_{t+1}(x_t^*)\\
&\leq(G_{f_t}+M_{t+1})\enorm{\bar{x}_t-x_t^*}
\leq(G_{f_t}+M_{t+1})\Omega_t.
\end{aligned}
\end{equation*}
The result follows from $\enorm{y_t^*}\leq\enorm{y_t^*}_1$.
\end{proof}

Thus one may take $Y_t=(G_{f_t}+M_{t+1})\Omega_t/\kappa_t$. The condition is readily checked when a uniform strict feasibility margin can be constructed. Pointwise strict feasibility over infinitely many scenarios need not give a uniform margin. If strong duality and dual attainment are not assumed in advance, the usual Slater condition with a strictly feasible point in $\operatorname{ri}(X_t)$ ensures them here.

The bounds can be verified backward, starting with $Y_3=G_{f_3}\Omega_3/\kappa_3$ and then using the established estimates for $M_3$ and $M_2$ to obtain $Y_2$ and $Y_1$. The same backward argument extends directly to $T$ stages.

\section{The A-DSPD Algorithm for the Multi-stage Problem}\label{sec:multi-stage}
In this section, we generalize the A-DSPD algorithm to multi-stage stochastic optimization problems with $T > 3$ stages. Our approach utilizes a backward decomposition strategy, beginning at the final stage $T$. For a given $x_{T-1} \in X_{T-1}$ and scenario $\xi_{[T]} \in \Xi_{[T]}$, we define the stochastic value function $V_T(x_{T-1},\xi_{[T]})$ as well as the corresponding conditional expectation $v_{T}\left(x_{T-1}, \xi_{[T-1]}\right)$:
\begin{align}
V_T(x_{T-1},\xi_{[T]})&:= \min_{x_T\in X_T}\max_{y_T\in\bbr^{m_T}_+}\{\enp{A_T x_T - b_T + B_T x_{T-1}, y_T} + f_T(x_T,c_T)\},\label{mu-last}\\
v_T(x_{T-1},\xi_{[T-1]})&:= \bbe_{\xi_T \vert \xi_{[T-1]}}[V_T(x_{T-1}, \xi_{[T]})].\notag
\end{align}
Proceeding recursively for $t = T-1, \ldots, 2$, the value function $v_{t+1}(x_t, \xi_{[t]})$ is substituted into the $t$-th stage problem. At each stage, we define the stochastic value function $V_t(x_{t-1},\xi_{[t]})$ and its expected value $v_t(x_{t-1}, \xi_{[t-1]})$ as follows:
\begin{align}
V_t(x_{t-1},\xi_{[t]})&:= \min_{x_t\in X_t}\max_{y_t\in\bbr^{m_t}_+}\{\enp{A_t x_t-b_t+B_t x_{t-1}, y_t} + f_t(x_t, c_t) + v_{t+1}(x_t, \xi_{[t]})\},\label{mu-middle}\\
v_t(x_{t-1},\xi_{[t-1]})&:= \bbe_{\xi_t\vert \xi_{[t-1]}}[V_t(x_{t-1}, \xi_{[t]})].\notag
\end{align}
Finally, the expected value function $v_2(x_1, \xi_{[1]})$ obtained from the second stage is substituted into the first-stage problem in~\eqref{multi}, yielding the following reformulation for the multi-stage problem:
\begin{equation}\label{mufi}
\begin{aligned}
&\min_{x_1 \in X_1} \ f_1(x_1, c_1) + v_2(x_1, \xi_{[1]})\\
&\st \quad A_1 x_1 - b_1 \leq 0.
\end{aligned}
\end{equation}

In what follows, we extend the A-DSPD algorithm introduced in Section \ref{sec:three-stage} to the multi-stage stochastic optimization problem \eqref{mufi}. 
We develop a backward recursion that begins at the terminal stage by computing a stochastic $\epsilon$-subgradient of $v_{T}$ at $x_{T-1}$ and then propagates these subgradients backward to obtain $\epsilon$-subgradients of $v_{t}$ at $x_{t-1}$ for $t = T-1, \ldots, 2$. 
The resulting $\epsilon$-subgradient of $v_{2}$ at $x_{1}$ is the input to the first-stage problem~\eqref{mufi}, which closes the recursion.

We now present the A-DSPD algorithm for solving the multi-stage optimization problem~\eqref{mufi}. For notational consistency, the algorithm follows the same local initialization and stagewise notation as the three-stage algorithm. Specifically, we use the following notation:
   \begin{itemize}
	\item[--] $N_{t}$ is the number of iterations for stage $t$ and $k_{t}$ is the corresponding index, i.e., $k_{t}=1, \ldots, N_{t}$.
	\item[--] $\xi^{k_t}_{[t+1]}:=(\xi_{[t]},\xi^{k_t}_{t+1})$ is the full scenario history obtained by appending the $k_t$-th stage-$(t+1)$ scenario to $\xi_{[t]}$.
\end{itemize}

\begin{algorithm}[htbp]  
	\caption{A-DSPD Algorithm for Multi-stage Problem \eqref{multi}}  \label{A-DSPD-mu}
	\begin{algorithmic}[0]
		\Require initial points $\{x_t^1\}=\{\tilde{x}_t^1\}=\{\bar{x}_t^1\}$, $\{y_t^1\}=\{\tilde{y}_t^1\}=0$
		\State Start with procedure A-DSPD$(1,0,\xi_{[1]})$
		\State \textbf{Procedure:} A-DSPD$(t,u,\xi_{[t]})$
		\State Initialize the local primal variables at the prescribed point in $X_t$ and set $y_t^1=\tilde y_t^1=0$.
			\For{$k_t=1,\ldots,N_t$}
		\If {$t<T$}
		\State Sample $\xi_{t+1}^{k_t}$ from the conditional distribution of $\xi_{t+1}$ given $\xi_{[t]}$.
        \State Set $\xi_{[t+1]}^{k_t}=(\xi_{[t]},\xi_{t+1}^{k_t})$.
		\State $(\tilde{x}_{t+1}^{N_{t+1}+1},\tilde{y}_{t+1}^{N_{t+1}+1},\bar{x}_{t+1}^{N_{t+1}+1},x_{t+1}^{N_{t+1}+1},y_{t+1}^{N_{t+1}+1})=\text{A-DSPD}(t+1,x_t^{k_t},\xi_{[t+1]}^{k_t})$
		\State $G_t(x^{k_t}_t)=B^\top_{t+1} \tilde{y}_{t+1}^{N_{t+1}+1}$
		\Else
		\State $G_T(x^{k_T}_T)=0$
		\EndIf
		\State $(\tilde{x}_t^{k_t+1},\tilde{y}_t^{k_t+1},\bar{x}_t^{k_t+1},x_t^{k_t+1},y_t^{k_t+1})
		=\operatorname{ASPDT}(\tilde{x}_t^{k_t},\tilde{y}_t^{k_t},x_t^{k_t},y_t^{k_t},\bar{x}_t^{k_t},G_t(x_t^{k_t}),u,\xi_t,X_t,\tau^t_{k_t},\eta^t_{k_t})$
		\EndFor
		\State \Return $(\tilde{x}_t^{N_t+1},\tilde{y}_t^{N_t+1},\bar{x}_t^{N_t+1},x_t^{N_t+1},y_t^{N_t+1})$
		\Ensure $(\tilde{x}^{N_1+1}_1,\tilde{y}^{N_1+1}_1)$
	\end{algorithmic}  
\end{algorithm}  
To establish the convergence of the A-DSPD algorithm for~\eqref{mufi}, we impose the following assumptions.

\begin{assu}\label{ass3}
	The following conditions are assumed to hold:
	\begin{itemize}
		\item[$\operatorname{(i)}$] For each $t=1,\ldots,T$, the function $f_t(\cdot,c_t)$ is smooth and convex, and there exist deterministic constants $L_{f_t},G_{f_t}\geq 0$ such that, almost surely,
		\begin{equation}\label{Lif-mu}
		f_t(x'_t, c_t) \leq f_t(x_t, c_t) + \inp{\nabla f_t(x_t, c_t), x'_t - x_t} + \frac{L_{f_t}}{2} \norm{x'_t - x_t}^2,
        \quad \norm{\nabla f_t(\tilde{x}_t,c_t)}\leq G_{f_t},
		\end{equation}
		for all $x_t,x'_t,\tilde{x}_t\in X_t$, uniformly over the scenario data. As in Assumption~\ref{ass2}, $L_{f_t}$ controls smoothness and $G_{f_t}$ bounds the gradient norm.
		\item[$\operatorname{(ii)}$] The linear operators $B_t$ are bounded; specifically, there exists a constant $L_{B_t} > 0$ such that
		\begin{equation}\label{mubound}
		\norm{B_t} \leq L_{B_t}, \quad \text{for all } t=2,\ldots,T.
		\end{equation}
		\item[$\operatorname{(iii)}$] For each $t=1,\ldots,T$, there exist deterministic constants $\underline{a}_t$ and $\overline{a}_t$ such that, almost surely,
        \begin{equation}\label{Abound-mu}
            0<\underline{a}_t\leq\enorm{A_t}\leq\overline{a}_t<\infty.
        \end{equation}
        These bounds control the operator norms in the step sizes and the estimates for the dual iterates.
        \item[$\operatorname{(iv)}$] For each $t=1,\ldots,T$, there exists a deterministic constant $Y_t<\infty$ and a jointly measurable selection of optimal dual solutions $y_t^*(u,\xi_{[t]})\in\bbr^{m_t}_+$ such that, outside a common null set of scenario histories,
        \begin{equation*}
            \enorm{y_t^*(u,\xi_{[t]})}\leq Y_t
            \quad\text{for every }u\in X_{t-1}\quad(t\geq2).
        \end{equation*}
        At $t=1$, this means that a fixed optimal multiplier $y_1^*$ satisfies $\enorm{y_1^*}\leq Y_1$. The constants $Y_t$ are uniform over predecessor decisions and scenario histories and do not depend on the requested accuracy or the iteration counts. Only one such selection is required; the entire dual optimal set need not be bounded. Sufficient conditions are discussed in Subsection~\ref{subsec:bounded-multipliers}.
    \end{itemize}
\end{assu}

Using the same diameter bounds and stagewise conventions for $t=1,\ldots,T$, with $0<\bar{\epsilon}\leq1$ and the same rounding convention for iteration budgets, we first analyze the last-stage problem.
\begin{lemma}\label{mu-A-DSAle1}
	Suppose Assumption~\ref{ass3} holds and let $\{\tilde{x}^{k_T}_T,\tilde{y}^{k_T}_T\}_{k_T=1}^{N_T+1}$ be the sequence generated by Algorithm~\ref{A-DSPD-mu}. Assume further that the parameters $\alpha^T_{k_T}$, $\beta^T_{k_T}$, $\theta^T_{k_T}$, $\eta^T_{k_T}$, and $\tau^T_{k_T}$ are chosen according to~\eqref{paraco1} (with $M=0$ and $\epsilon=0$), and the number of iterations $N_T\geq 1$ is set as
	\begin{equation}\label{muNT}
	N_T = \max \left\{
		\frac{\sqrt{8 T L_{f_T}}\,\Omega_T}{\sqrt{\bar{\epsilon}}},\;
		\frac{8T\Omega_T^2 \overline{a}_T}{\bar{\epsilon}},\;
		\frac{\sqrt{8T}\,\Omega_T}{\sqrt{\bar{\epsilon}}},\;
		\frac{8T\overline{a}_TY_T^2}{\bar{\epsilon}},\;
		\frac{4T \overline{a}_T}{\bar{\epsilon}} \cdot \left(2+4Y_T^2\right)
	\right\},
	\end{equation}
	then $B_T^\top \tilde{y}_T^{N_T+1}$ is a stochastic $\frac{\bar{\epsilon}}{T}$-subgradient of the value function $v_T$ at $x_{T-1}^{k_{T-1}}$. Moreover, for any fixed $\xi_{[T-1]}$, there exists a constant $M_T$ such that 
	\[
		|v_T(x_{T-1}, \xi_{[T-1]}) - v_T(x'_{T-1}, \xi_{[T-1]})|\leq M_T \enorm{x_{T-1} - x'_{T-1}},
	\]
	for all $x_{T-1}, x'_{T-1} \in X_{T-1}$, and
	\[
		\bbe[\enorm{B_T^\top \tilde{y}_T^{N_T+1}}^2\vert \xi_{[T-1]} ] \leq M_T^2.
	\]
	In addition, the following bounds hold:
	\begin{equation}\label{mu-result1}
	\begin{aligned}
	& |f_T(\tilde{x}_T^{N_T+1}, c_T) - f_T(x^*_T, c_T)| \leq \frac{\bar{\epsilon}}{T} \\
	& \enorm{[A_T \tilde{x}_T^{N_T+1} - b_T + B_T x_{T-1}^{k_{T-1}}]_+} \leq \frac{\bar{\epsilon}}{T}.
	\end{aligned}
	\end{equation}
\end{lemma}
\begin{proof}
	We apply the \iaspd{} method to the final-stage saddle-point problem~\eqref{mu-last}. As there are no subsequent stages, this setting yields a deterministic saddle-point problem for a given $x^{k_{T-1}}_{T-1}$ and $\xi^{k_{T-1}}_{[T]}$. Consequently, we have $v=0$, $M=0$, and $\epsilon=0$. Denote $u = x^{k_{T-1}}_{T-1}$, $\xi = \xi^{k_{T-1}}_T$, $A = A_T$, $B = B_T$, $L_f = L_{f_T}$, $\Omega_X = \Omega_T$, $\tilde{x}^{N+1} = \tilde{x}^{N_T+1}$, $x^* = x_T^*$, $\tilde{y}^{N+1} = \tilde{y}_T^{N_T+1}$, $y^* = y^*_T$, and $y^1 = y^1_T = 0$.
	By Assumption~\ref{ass3}(iv), choose the selected optimal multiplier $y_T^*$, so that
	\begin{equation}\label{mu-yy3b}
    \enorm{y^*_T}\leq Y_T.
\end{equation}
	The remainder of the argument is identical to the proof of Lemma~\ref{A-DSAle1}, with the subscript $3$ replaced by $T$, the point $x_2^{j}$ replaced by $x_{T-1}^{k_{T-1}}$, and the error budget $\bar{\epsilon}/3$ replaced by $\bar{\epsilon}/T$, so that each of the four terms in the gap bound is now controlled by $\bar{\epsilon}/(4T)$ instead of $\bar{\epsilon}/12$, as guaranteed by the definition of $N_T$ in \eqref{muNT}. We record only the resulting bounds. First, Corollary~\ref{cor:last-stage}(a) (with $M=0$ and $\epsilon=0$), the definition of $N_T$, and \eqref{mu-yy3b} give
	\begin{equation*}
	\operatorname{gap}(\tilde{z}_T^{N_T+1})
	\leq \frac{2\Omega_T^2 L_{f_T}}{N_T^2} + \frac{2\Omega_T^2 \overline{a}_T}{N_T} + \frac{2\Omega_T^2}{N_T^2} + \frac{\overline{a}_TY_T^2}{N_T}
	\leq \frac{\bar{\epsilon}}{T},
	\end{equation*}
	which, by Lemma~\ref{subgradient}, establishes that $B_T^\top \tilde{y}_T^{N_T+1}$ is a $\bar{\epsilon}/T$-subgradient of the value function $V_T$ at $x_{T-1}^{k_{T-1}}$, that is, a stochastic $\bar{\epsilon}/T$-subgradient of $v_T$ at $x_{T-1}^{k_{T-1}}$. Second, repeating the derivation of \eqref{thirdbound1}--\eqref{thirdbound2} (which combines Corollary~\ref{cor:last-stage}(b), the bound $\norm{B_T}\leq L_{B_T}$ from \eqref{mubound}, Lemma~\ref{subg}, and \eqref{mu-yy3b}) shows that the Lipschitz and second-moment bounds in the statement hold with
	\begin{equation*}
	M_T=L_{B_T}\max\left\{Y_T,\sqrt{8\Omega_T^2+\frac{8\Omega_T^2(L_{f_T}+1)}{\underline{a}_T}+6Y_T^2}\right\}.
	\end{equation*}
	Finally, Corollary~\ref{cor:last-stage}(c) (with $M = 0$ and $\epsilon = 0$), the definition of $N_T$, and \eqref{mu-yy3b} yield \eqref{mu-result1} by exactly the same estimates as in the proof of Lemma~\ref{A-DSAle1}.
\end{proof}

The following lemma establishes the convergence properties of  A-DSPD for the intermediate stages.
\begin{lemma}\label{mu-A-DSAle2}
	Suppose that the parameters for the last-stage problem are chosen according to Lemma~\ref{mu-A-DSAle1}, and that Assumption~\ref{ass3} holds. Let $\{\tilde{x}_t^{k_t}, \tilde{y}_t^{k_t}\}_{k_t=1}^{N_t+1}$, for $t = T-1, \ldots, 2$, be a sequence generated by Algorithm~\ref{A-DSPD-mu}. Assume further that the parameters $\alpha^t_{k_t}$, $\beta^t_{k_t}$, $\theta^t_{k_t}$, $\eta^t_{k_t}$, and $\tau^t_{k_t}$ are set according to~\eqref{paraco2}, and the number of iterations $N_t\geq 1$ is set as
	\begin{equation}\label{muNt}
	N_t = \max\left\{
	\begin{aligned}
	&\frac{\sqrt{10T L_{f_t}}\,\Omega_t}{\sqrt{\bar{\epsilon}}},
	\frac{10T\Omega_t^2\overline{a}_t}{\bar{\epsilon}},\left(\frac{10T M_{t+1} \Omega_t^2}{\bar{\epsilon}}\right)^2,\left(\frac{5T\overline{a}_t}{\bar{\epsilon}} \cdot Y_t^2\right)^2,\left(\frac{70T M_{t+1}}{3\bar{\epsilon}}\right)^2, \\
	&~~~~~~~~~~~~~~~~~~~~~~~~~~~~~~\left(\frac{5T\overline{a}_t}{\bar{\epsilon}} \cdot \left(2+4Y_t^2\right)\right)^2
	\end{aligned}
	\right\},
	\end{equation}
	then $B_t^\top \tilde{y}_t^{N_t+1}$ is a stochastic $\frac{(T-t+1)\bar{\epsilon}}{T}$-subgradient of the value function $v_t$ at $x_{t-1}^{k_{t-1}}$. Furthermore, for any given $\xi_{[t-1]}$, there exists a constant $M_t$ such that the value function $v_t(\cdot, \xi_{[t-1]})$ is Lipschitz continuous on $X_{t-1}$ with Lipschitz constant $M_t$; that is,
	\[
	|v_t(x_{t-1}, \xi_{[t-1]}) - v_t(x_{t-1}', \xi_{[t-1]})| \leq M_t \enorm{x_{t-1} - x'_{t-1}}\quad \forall\, x_{t-1}, x_{t-1}' \in X_{t-1},
	\]
	and it holds that
	\[
	\bbe[\|B_t^\top \tilde{y}_t^{N_t+1}\|^2\vert\xi_{[t-1]}] \leq M_t^2.
	\]
	In addition,
	\begin{equation*}
	\begin{aligned}
	&\bbe[|f_t(\tilde{x}_t^{N_t+1}, c_t) + v_{t+1}(\tilde{x}_t^{N_t+1}) - f_t(x_t^*, c_t) - v_{t+1}(x_t^*)|\vert\xi_{[t]}]
	\leq \frac{(T-t+1)\bar{\epsilon}}{T},\\
	&\bbe[\enorm{[A_t \tilde{x}_t^{N_t+1} - b_t + B_t x_{t-1}^{k_{t-1}}]_+}\vert\xi_{[t]}] \leq \frac{(T-t+1)\bar{\epsilon}}{T}.
	\end{aligned}
	\end{equation*}
\end{lemma}
\begin{proof}
	We prove the result by backward induction. For $t=T-1$, Lemma~\ref{mu-A-DSAle1} provides the required oracle conditions with $M=M_T$ and $\epsilon=\bar{\epsilon}/T$. For $t\leq T-2$, the induction hypothesis at stage $t+1$ provides the same conditions with $M=M_{t+1}$ and $\epsilon=(T-t)\bar{\epsilon}/T$. Thus, the following argument establishes both the base case and the induction step.
	We apply the \iaspd{} algorithm to the intermediate-stage problem in Algorithm~\ref{A-DSPD-mu}. For fixed $(x^{k_{t-1}}_{t-1}, \xi^{k_{t-1}}_{[t]})$, denote $v = v_{t+1}$, $M = M_{t+1}$, and $\epsilon = \frac{(T-t)\bar{\epsilon}}{T}$. Let $u = x^{k_{t-1}}_{t-1}$, $\xi = \xi^{k_{t-1}}_t$, $A = A_t$, $B = B_t$, $L_f = L_{f_t}$, $\Omega_X = \Omega_t$, $\tilde{x}^{N+1} = \tilde{x}_t^{N_t+1}$, $x^* = x^*_t$, $\tilde{y}^{N+1} = \tilde{y}_t^{N_t+1}$, $y^* = y^*_t$, and $y^1 = y^1_t = 0$. By Assumption~\ref{ass3}(iv), choose the selected optimal multiplier $y_t^*$, so that
	\begin{equation}\label{muyy2b}
    \enorm{y^*_t}\leq Y_t.
\end{equation}
	The remainder of the argument is identical to the proof of Lemma~\ref{A-DSAle2}, with the subscripts $(2,3)$ replaced by $(t,t+1)$, the point $x^{i}_1$ replaced by $x^{k_{t-1}}_{t-1}$, the inherited error $\bar{\epsilon}/3$ replaced by $\frac{(T-t)\bar{\epsilon}}{T}$, and the per-term threshold $\bar{\epsilon}/15$ replaced by $\bar{\epsilon}/(5T)$, as guaranteed by the definition of $N_t$ in \eqref{muNt}. We record only the resulting bounds. First, Corollary~\ref{cor:middle-stage}(a), the definition of $N_t$, and \eqref{muyy2b} give
	\begin{equation*}
	\begin{aligned}
	\bbe[\operatorname{gap}(\tilde{z}^{N_t+1}_t)\vert\xi_{[t]}]
	&\leq \frac{2\Omega_t^2 L_{f_t}}{N_t^2} + \frac{2\Omega_t^2\overline{a}_t}{N_t} + \frac{2M_{t+1}\Omega_t^2}{\sqrt{N_t}}
	+ \frac{\overline{a}_tY_t^2}{\sqrt{N_t}} + \frac{12M_{t+1}}{3\sqrt{N_t}} + \frac{2M_{t+1}}{3 N_t^{3/2}} + \frac{(T-t)\bar{\epsilon}}{T}\\
	&\leq \frac{(T-t+1)\bar{\epsilon}}{T},
	\end{aligned}
	\end{equation*}
	which, by Lemma~\ref{subgradient}, implies that $B_t^\top\tilde{y}^{N_t+1}_t$ is a $\frac{(T-t+1)\bar{\epsilon}}{T}$-subgradient of the value function $V_t$ at $x^{k_{t-1}}_{t-1}$ (i.e., a stochastic $\frac{(T-t+1)\bar{\epsilon}}{T}$-subgradient of $v_t$ at $x^{k_{t-1}}_{t-1}$).  Second, define the deterministic constant
    \begin{equation*}
    \begin{aligned}
    H_t:=1+\max\bigl\{&1,\sqrt{10T L_{f_t}}\Omega_t,10T\Omega_t^2\overline{a}_t,
    (10T M_{t+1}\Omega_t^2)^2,\\
    &(5T\overline{a}_tY_t^2)^2,\;(70T M_{t+1}/3)^2,
    [5T\overline{a}_t(2+4Y_t^2)]^2\bigr\}.
    \end{aligned}
    \end{equation*}
    Then $\bar{\epsilon}^2N_t\leq H_t$ for $0<\bar{\epsilon}\leq1$, including the rounding of $N_t$. Repeating the derivation of \eqref{bineq1}--\eqref{secondbound2} (which combines Corollary~\ref{cor:middle-stage}(b) with $M=M_{t+1}$ and $\epsilon=\frac{(T-t)\bar{\epsilon}}{T}$, the bound $\norm{B_t}\leq L_{B_t}$ from \eqref{mubound}, Lemma~\ref{subg}, and \eqref{muyy2b}) shows that the Lipschitz and second-moment bounds in the statement hold with
	\begin{footnotesize}
	\begin{equation*}
	M_t=L_{B_t}\max\left\{Y_t,\sqrt{8\Omega_t^2+\frac{24\Omega_t^2(L_{f_t}+M_{t+1})+112M_{t+1}+24\sqrt{H_t}}{3\underline{a}_t}+6Y_t^2}\right\}.
	\end{equation*}
	\end{footnotesize}
	Finally, Corollary~\ref{cor:middle-stage}(c), the definition of $N_t$, and \eqref{muyy2b} yield the last two bounds in the statement by exactly the same estimates as in the proof of Lemma~\ref{A-DSAle2}.
	The result now follows by induction.
\end{proof}

We now present the main convergence result for the A-DSPD algorithm applied to the multi-stage stochastic optimization problem~\eqref{mufi}.
\begin{theorem}\label{th:ADSPD-multi}
	Suppose the parameters for the terminal and intermediate stages in the \text{A-DSPD} algorithm are chosen according to Lemmas~\ref{mu-A-DSAle1} and~\ref{mu-A-DSAle2}. Assume further that Assumption~\ref{ass3} holds and that the sequence $\{\tilde{x}_1^{i}, \tilde{y}_1^{i}\}_{i=1}^{N_1+1}$ is generated by Algorithm~\ref{A-DSPD-mu}. Let the parameters $\alpha^1_{k_1}$, $\beta^1_{k_1}$, $\theta^1_{k_1}$, $\eta^1_{k_1}$, and $\tau^1_{k_1}$ be set according to~\eqref{paraco3}, and the number of iterations $N_1\geq 1$ is set as
\begin{equation}\label{muN1}
	N_1 = \max\left\{
	\frac{\sqrt{10T L_{f_1}}\,\Omega_{1}}{\sqrt{\bar{\epsilon}}},
	\frac{10T\Omega_1^2\overline{a}_1}{\bar{\epsilon}},
	\left(\frac{10TM_2\Omega_1^2}{\bar{\epsilon}}\right)^2,
	\left(\frac{70TM_2}{3\bar{\epsilon}}\right)^2,
	\frac{5T\overline{a}_1}{\bar{\epsilon}} \cdot \left(2+4Y_1^2\right)
	\right\}.
\end{equation}
Then the following bounds hold:
\begin{equation*}
\begin{aligned}
	&\bbe[|f_1(\tilde{x}_1^{N_1+1},c_1)+v_2(\tilde{x}_1^{N_1+1})-f_1(x_1^*,c_1)-v_2(x_1^*)|\vert\xi_{[1]}] \leq \bar{\epsilon}, \\
	&\bbe[\enorm{[A_1 \tilde{x}_1^{N_1+1}-b_1]_+}\vert\xi_{[1]}] \leq \bar{\epsilon}. \\
\end{aligned}
\end{equation*}
\end{theorem}
\begin{proof}
	We apply the \iaspd{} algorithm to the first-stage problem \eqref{mufi} in Algorithm~\ref{A-DSPD-mu}. Let $v = v_2$, $M = M_2$, $\epsilon = \frac{(T-1)\bar{\epsilon}}{T}$, $u = 0$, $\xi = \xi_1$, $A = A_1$, $B = 0$, $L_f = L_{f_1}$, $\Omega_X = \Omega_{1}$, $\tilde{x}^{N+1} = \tilde{x}_1^{N_1+1}$, $x^* = x^*_1$, $\tilde{y}^{N+1} = \tilde{y}_1^{N_1+1}$, $y^* = y^*_1$, and $y^1 = y^1_1 = 0$. By Assumption~\ref{ass3}(iv), choose the selected optimal multiplier $y_1^*$, so that
	\begin{equation}\label{muyy1b}
    \enorm{y^*_1}\leq Y_1.
\end{equation}
	Invoking Corollary~\ref{cor:first-stage}(c), the definition of $N_1$ in \eqref{muN1}, and \eqref{muyy1b}, and repeating the estimates in the proof of Theorem~\ref{th:ADSPD-three} with the inherited error $2\bar{\epsilon}/3$ replaced by $\frac{(T-1)\bar{\epsilon}}{T}$ and the per-term threshold $\bar{\epsilon}/15$ replaced by $\bar{\epsilon}/(5T)$, we obtain
	\begin{equation*}
	\bbe[|f_1(\tilde{x}_1^{N_1+1},c_1)+v_2(\tilde{x}_1^{N_1+1})-f_1(x_1^*,c_1)-v_2(x_1^*)|\vert\xi_{[1]}]
	\leq \frac{5\bar{\epsilon}}{5T} + \frac{(T-1)\bar{\epsilon}}{T} = \bar{\epsilon}
	\end{equation*}
	and
	\begin{equation*}
	\bbe[\enorm{[A_1 \tilde{x}_1^{N_1+1} - b_1]_+}\vert\xi_{[1]}]\leq \bar{\epsilon}.
	\end{equation*}
	This completes the proof.
\end{proof}

For fixed $T$ and problem constants, including the uniform multiplier bounds $Y_t$, the constants $M_t$ above can be chosen independently of $\bar{\epsilon}$. Thus, \eqref{muNT} gives $N_T=O(1/\bar{\epsilon})$, whereas \eqref{muNt} and \eqref{muN1} give $N_t=O(1/\bar{\epsilon}^2)$ for $t<T$. Consequently, A-DSPD uses $O(1/\bar{\epsilon}^{2(T-1)})$ sampled scenarios and $O(1/\bar{\epsilon}^{2T-1})$ primal--dual updates.

\newcommand{\zg}[1]{\comm[zg]{#1}}
\section{Numerical Study}
\label{sec:numerical}

We evaluate \ADSPD{} on two three-stage portfolio selection models that share the same nested decision structure as Problem~\eqref{three}.
Both models use the same interstage $\ell_1$ turnover-linking constraints.
At each node of stage $t$, $x_t$ denotes the full decision vector and $w_t$ its portfolio-weight block; $c_t$ denotes the stage-cost data contained in the scenario tuple $\xi_t$. The two models below specify the remaining decision blocks and the corresponding cost functions.
Section~\ref{subsec:exp-protocol} summarizes the compared methods, the evaluation metrics, and the randomization scheme; Sections~\ref{subsec:exp-mv} and~\ref{subsec:exp-cvar} describe the two models and report the corresponding results. All experiments are run on a single CPU node with a 32-core AMD EPYC 7T83 processor and 120\,GB of RAM (Ubuntu 22.04, Python 3.12). 
\subsection{Experimental settings}
\label{subsec:exp-protocol}

\paragraph{Compared methods.}
We compare \ADSPD{} with \DSA{}.
The two methods are run under identical conditions: for each instance and random seed, they share the iteration budget $(N_1,N_2,N_3)$ and the same sequence of sampled scenarios.
Both methods report the averaged primal iterate \(\tilde{x}_t\) and use the most recent dual iterates to warm-start stages 2 and 3 as a practical acceleration strategy. 
For the step-size schedules,
\ADSPD{} uses~\eqref{paraco3}, \eqref{paraco2}, and~\eqref{paraco1} at stages~$1$, $2$, and~$3$, and \DSA{} uses the choices of \citet[Cor.~6 and~7]{lan2021dynamic}. 
The remaining constants are estimated from the scenario tree.
For the inner primal--dual step, \DSA{} requires a proximal update of the stage cost, which we approximate via Nesterov's accelerated gradient (NAG) method~\citep{beck2017first,nesterov1988approach} with a relative tolerance of $10^{-5}$ and a maximum of $200$ iterations.

\paragraph{Random seeds and reported statistics.}
Each run is indexed by a random seed that determines both the scenario tree and the sampling path followed within the run.
We report means and two-sided 95\% confidence intervals based on Student's $t$ distribution over the five independent seeds $\{40,41,42,43,44\}$.
In all convergence plots, solid curves represent seed averages and shaded bands depict the corresponding 95\% confidence intervals based on Student's $t$ distribution.

Let $\mathcal I_2$ index the stage-$2$ nodes and $\mathcal C(i)$ the stage-$3$ children of node $i$. Write $p_i$ for the stage-$2$ probability and $p_{j\mid i}$ for the conditional probability of child $j$. A feasible policy is $\pi=(x_1,\{x_{2,i}\}_i,\{x_{3,j\mid i}\}_{i,j})$. Its probability-weighted objective on the finite scenario tree is
\begin{equation}
\label{eq:tree-objective}
  F(\pi)
  = f_1(x_1,c_1)
  + \sum_{i\in\mathcal I_2} p_i f_2(x_{2,i},c_{2,i})
  + \sum_{i\in\mathcal I_2}\sum_{j\in\mathcal C(i)} p_i p_{j\mid i} f_3(x_{3,j\mid i},c_{3,j\mid i}),
\end{equation}
where $c_{2,i}$ and $c_{3,j\mid i}$ are the cost data at the corresponding nodes. For seed $s$, we write $F_s$ for this objective on the generated tree.
For each seed~$s$, we solve the deterministic equivalent formulation (extensive form) once using the COPT solver~\citep{copt} to obtain a certified optimum $F_s^\star$.
Given the learned first-stage portfolio $\widehat w_{1,s}$, we fix $w_1=\widehat w_{1,s}$ and reoptimize all scenario-dependent recourse decisions in the extensive form (equivalently, one quadratic program (QP) per stage-$2$ parent).
This yields
\begin{equation}
  Q_s(\widehat w_{1,s})
  :=\min_{\pi:\,w_1=\widehat w_{1,s}} F_s(\pi),
  \qquad
  \operatorname{Gap}_s
  :=\frac{Q_s(\widehat w_{1,s})-F_s^\star}{|F_s^\star|}.
  \label{eq:fixed-w1-gap}
\end{equation}

For each experiment, we report means and 95\% confidence intervals based on Student's $t$ distribution at the final checkpoint, and plot $\operatorname{Gap}_s$ against outer iteration, wall-clock time, and cumulative matrix--vector product (matvec) count.

\subsection{Mean--variance portfolio model}
\label{subsec:exp-mv}

For a portfolio of $n$ assets, $w_t$ denotes the weight vector at stage~$t$ and satisfies $w_t\in\Delta_n:=\{w\in\bbr^n\mid \mathbf{1}^\top w=1,\;w\ge0\}$, which imposes full investment and excludes short selling. We use this simplex notation for both portfolio models.
The stagewise cost is
\begin{equation}
\label{eq:mv-stage-cost}
  f_t(w_t,c_t)
  =\frac{1}{2}w_t^\top\Sigma_t w_t-\mu\, r_t^\top w_t,
  \qquad t=1,2,3,
\end{equation}
Here $c_t=(\Sigma_t,r_t)$ consists of the covariance matrix and expected-return vector at the current node, and $\mu$ controls the trade-off between risk and expected return. These cost data form part of the full scenario tuple $\xi_t$ in~\eqref{multi}.

For $t=2,3$, we measure turnover between consecutive stages by the total absolute change in portfolio weights, $\|w_t-w_{t-1}\|_1$. Given a turnover limit $\psi_t>0$, we impose $\|w_t-w_{t-1}\|_1\le\psi_t$.
Introducing an auxiliary vector $s_t\ge 0$, we write this constraint in linear form as
\[
  w_t-w_{t-1}\le s_t,\quad
  w_{t-1}-w_t\le s_t,\quad
  \mathbf{1}^\top s_t \le \psi_t.
\]
The full stage decisions are $x_1=w_1$ and $x_t=(w_t,s_t)$ for $t=2,3$. Since the stage cost depends only on the portfolio block, we use $f_t(x_t,c_t)$ and $f_t(w_t,c_t)$ interchangeably.
The linear constraints take the form $A_1x_1\leq b_1$ and
$A_tx_t+B_tx_{t-1}\leq b_t$ for $t=2,3$, where
\[
A_1=
\begin{bmatrix}
-I_n\\
\mathbf{1}^{\top}\\
-\mathbf{1}^{\top}
\end{bmatrix},
\qquad
b_1=
\begin{bmatrix}
\mathbf{0}\\1\\-1
\end{bmatrix}.
\]
In the reported instances, we set $\mu=1$ and $\psi_2=\psi_3=0.5$.

The uniform Slater condition in Proposition~\ref{prop:slater-multipliers} can be checked for these turnover constraints. For $t=2,3$, take $w_t=w_{t-1}$ and $s_t=\psi_t\mathbf{1}/(2n)$. For every $w_{t-1}\in\Delta_n$, the residuals of $w_t-w_{t-1}\leq s_t$ and $w_{t-1}-w_t\leq s_t$ equal $-\psi_t\mathbf{1}/(2n)$, while that of $\mathbf{1}^\top s_t\leq\psi_t$ equals $-\psi_t/2$. Hence $\kappa_t=\psi_t/(2n)$ works uniformly when $\psi_t>0$; the simplex constraints remain in $X_t$. At the first stage, the explicit inequalities are redundant on
$X_1=\Delta_n$. Therefore, $y_1^\star=0$ is an optimal multiplier,
and the first-stage multiplier bound in Assumption~\ref{ass2}(iv)
holds with $Y_1=0$.

At the first stage, the model parameters are estimated from historical observations, whereas the parameters at stages~$2$ and~$3$ are scenario dependent and therefore vary with the realized stochastic process.

We construct the instances from daily log returns computed from adjusted closing prices downloaded from Yahoo Finance~\footnote{Data retrieved from Yahoo Finance (\url{https://finance.yahoo.com}) via the \texttt{yfinance} application programming interface (API).}. We retain assets with at least $95\%$ return coverage from the S\&P~500 universe.
The sample covers January~4, 2023 through December~31, 2025, and the reference date is $t_0=\mathtt{2024\textrm{-}12\textrm{-}31}$.
The scenario tree has three stages, with $J_2$ stage-$2$ nodes and $J_3$ stage-$3$ children per stage-$2$ node.
At stage~$1$, $(\Sigma_1,r_1)$ is estimated from a lookback window of 252 trading days ending at $t_0$ and is attached to the single root node. At stage~$2$, we sample $J_2$ windows of 63 consecutive trading days starting after $t_0$; each window defines one stage-$2$ node with probability $p_i=1/J_2$.
At stage~$3$, each stage-$2$ node has $J_3$ children, each defined by a 63-day window starting after its parent window, with conditional probability $p_{j\mid i}=1/J_3$.
Using the run seed, starting dates are sampled uniformly without replacement when their number is at least the requested branch count, and with replacement otherwise. Distinct starting dates can still give overlapping windows.
For each window $W$, let $|W|$ be its number of trading days and let $r_d$ denote the daily log-return vector. Write $\operatorname{Cov}_{d\in W}$ for the sample covariance and $\operatorname{Diag}(C)$ for the matrix obtained by retaining only the diagonal of $C$. We estimate the annualized return and covariance by
\[
\bar r_W=\frac{252}{|W|}\sum_{d\in W}r_d,\qquad
C_W=252\,\operatorname{Cov}_{d\in W}(r_d),\qquad
\widetilde\Sigma_W=0.9C_W+0.1\operatorname{Diag}(C_W),
\]
Let $m_j$ and $\sigma_j$ be the mean and standard deviation of asset $j$'s window returns $(\bar r_W)_j$ across the tree. We normalize the moments by
\[(r_W)_j=\frac{(\bar r_W)_j-m_j}{\sigma_j},\qquad
\Sigma_W=\frac{\widetilde\Sigma_W}{\max_V\|\widetilde\Sigma_V\|_2},
\]
where the maximum in the covariance normalization is taken over all windows in the tree.
All objective values and relative gaps reported below are computed
on the normalized scenario trees, which serve as in-sample computational benchmarks for comparing the two methods.

\paragraph{Problem instances.}
We consider three instances.
The small instance has $n=20$ assets, $J_2=J_3=10$ (hence $1+10+10^{2}=111$ scenario-tree nodes and $100$ root-to-leaf paths), and budget $(N_1,N_2,N_3)=(200,200,150)$.
The medium instance has $n=100$ assets, $J_2=J_3=30$ ($931$ nodes and $900$ paths), and budget $(800,800,150)$.
The large instance has $n=200$ assets, $J_2=J_3=50$ ($2551$ nodes and $2500$ paths), and budget $(1500,1500,150)$.

\begin{figure}[H]
  \centering
  \includegraphics[width=0.32\textwidth]{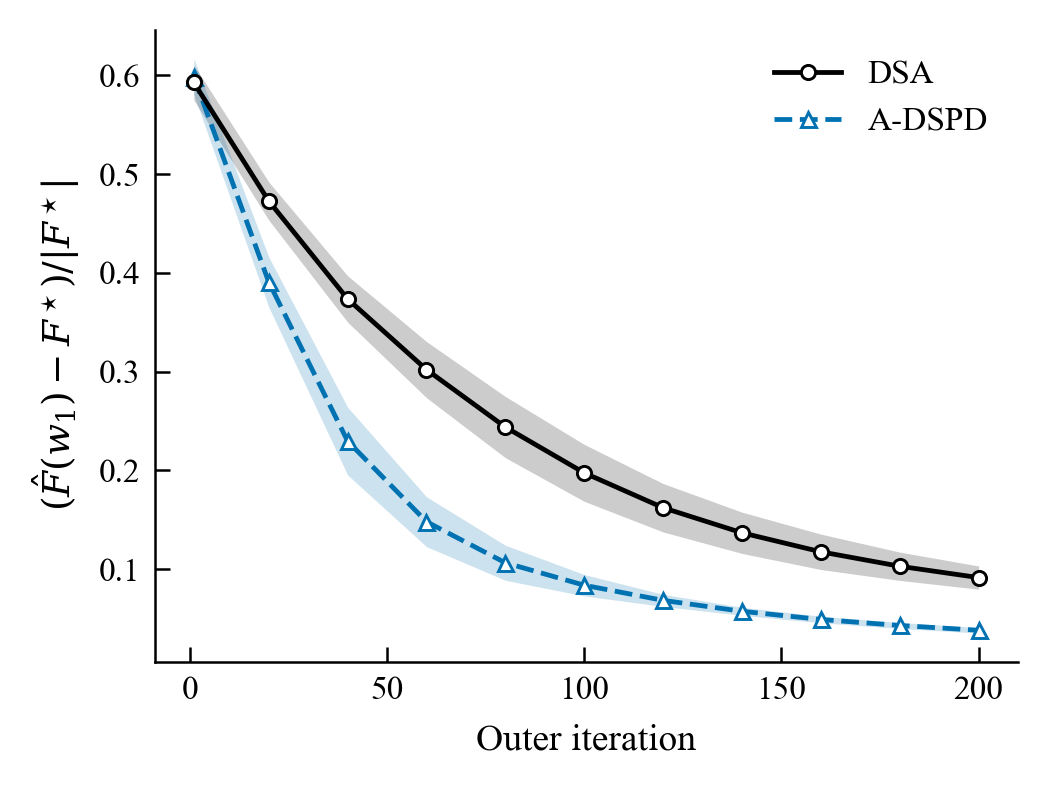}\hfill
    \includegraphics[width=0.32\textwidth]{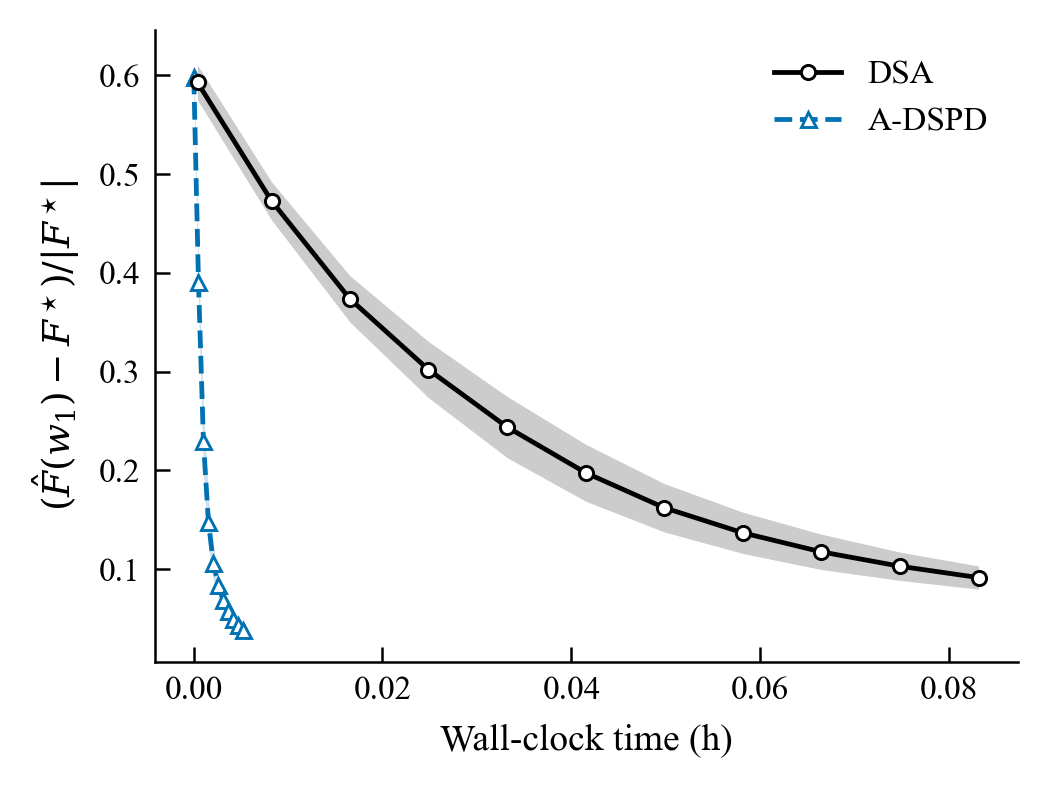}\hfill
  \includegraphics[width=0.32\textwidth]{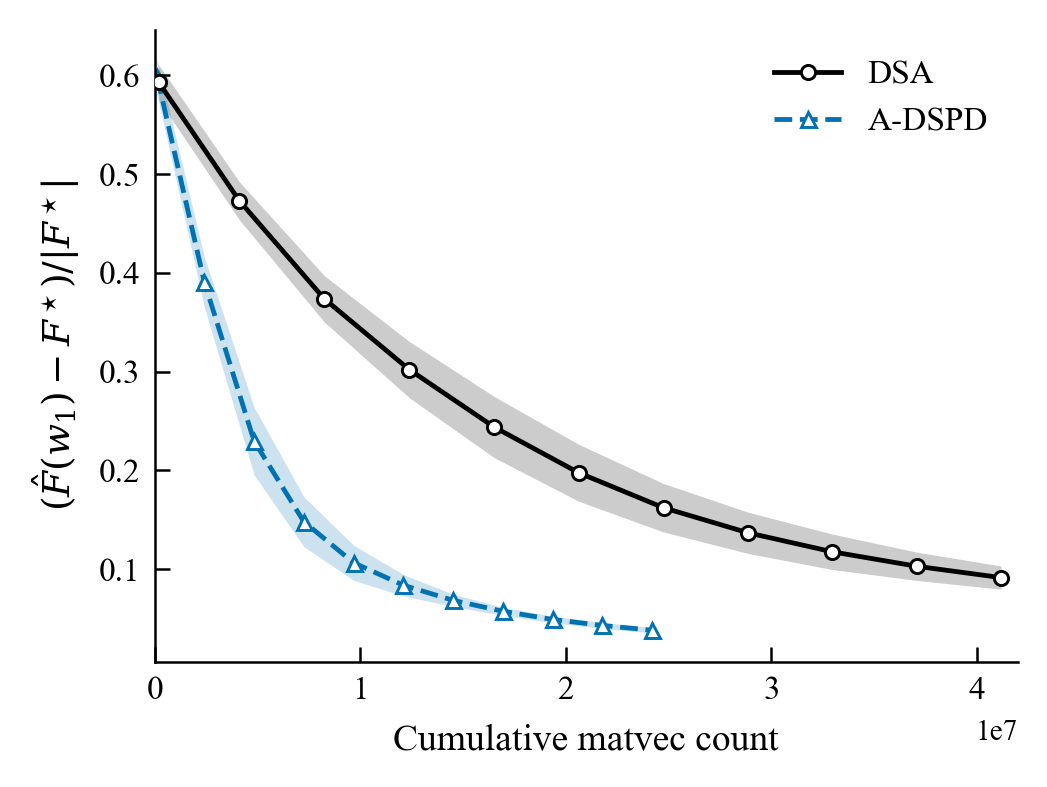}\\[-0.4em]
  {\footnotesize (a) Small instance ($n=20$ assets)}\\[0.3em]
  \includegraphics[width=0.32\textwidth]{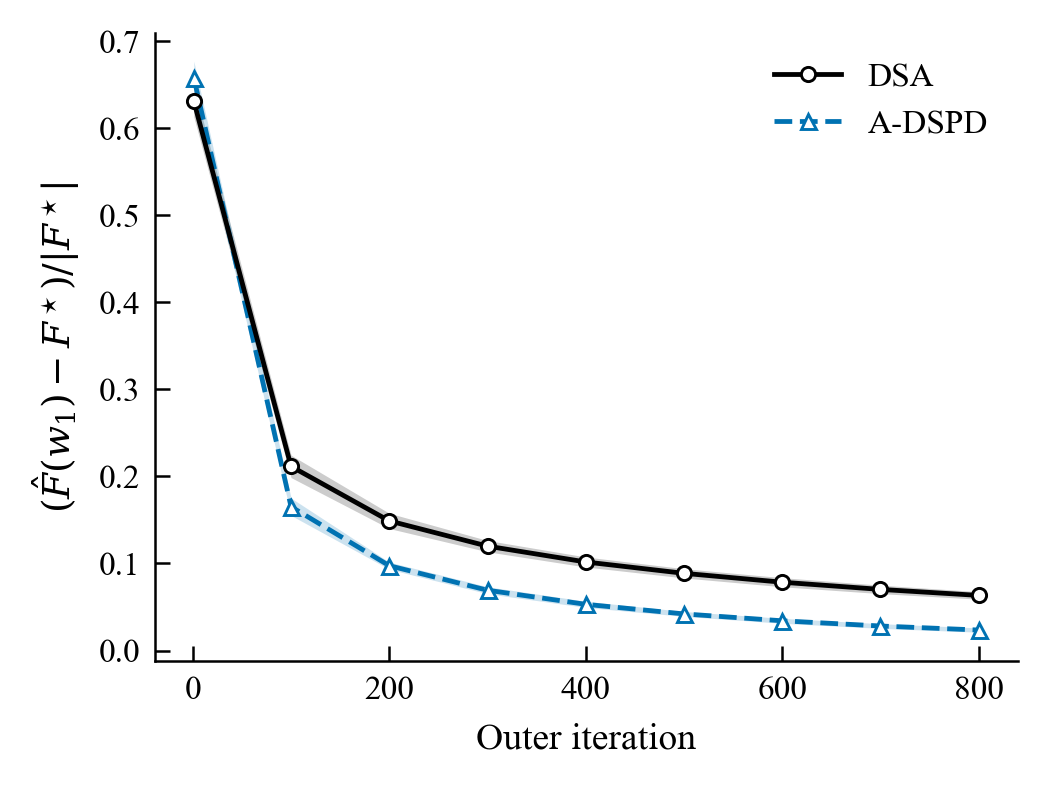}\hfill
  \includegraphics[width=0.32\textwidth]{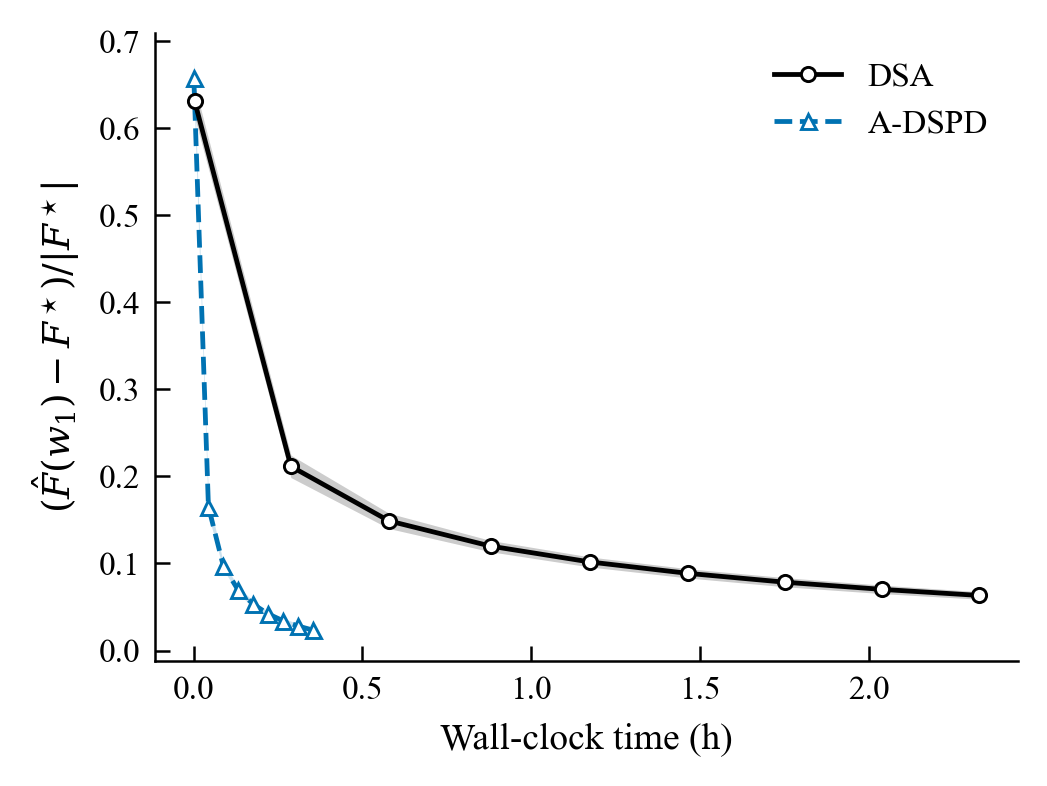}\hfill
  \includegraphics[width=0.32\textwidth]{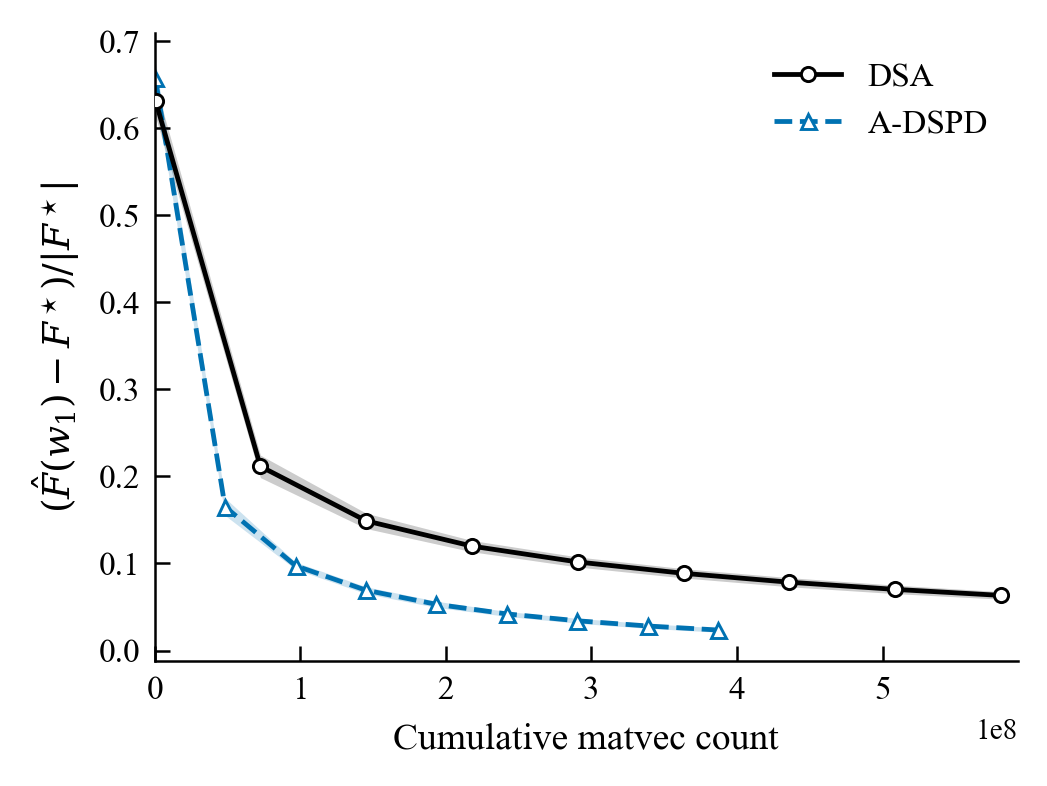}\\[-0.4em]
  {\footnotesize (b) Medium instance ($n=100$ assets)}\\[0.3em]
  \includegraphics[width=0.32\textwidth]{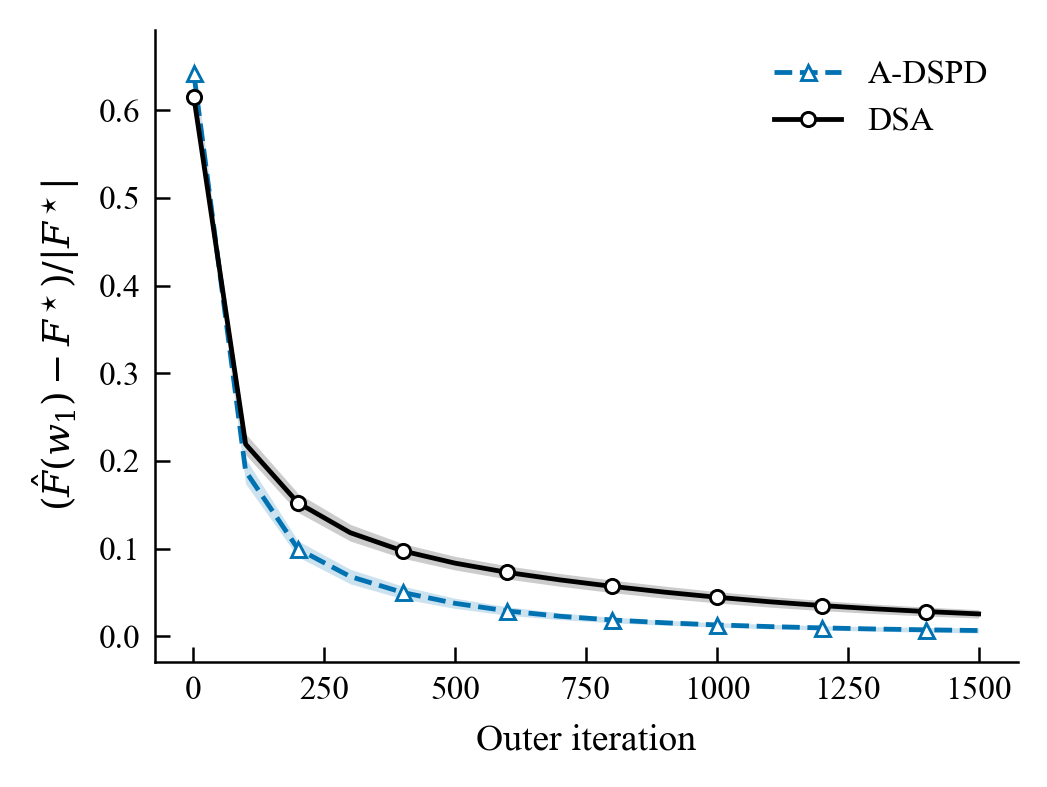}\hfill
  \includegraphics[width=0.32\textwidth]{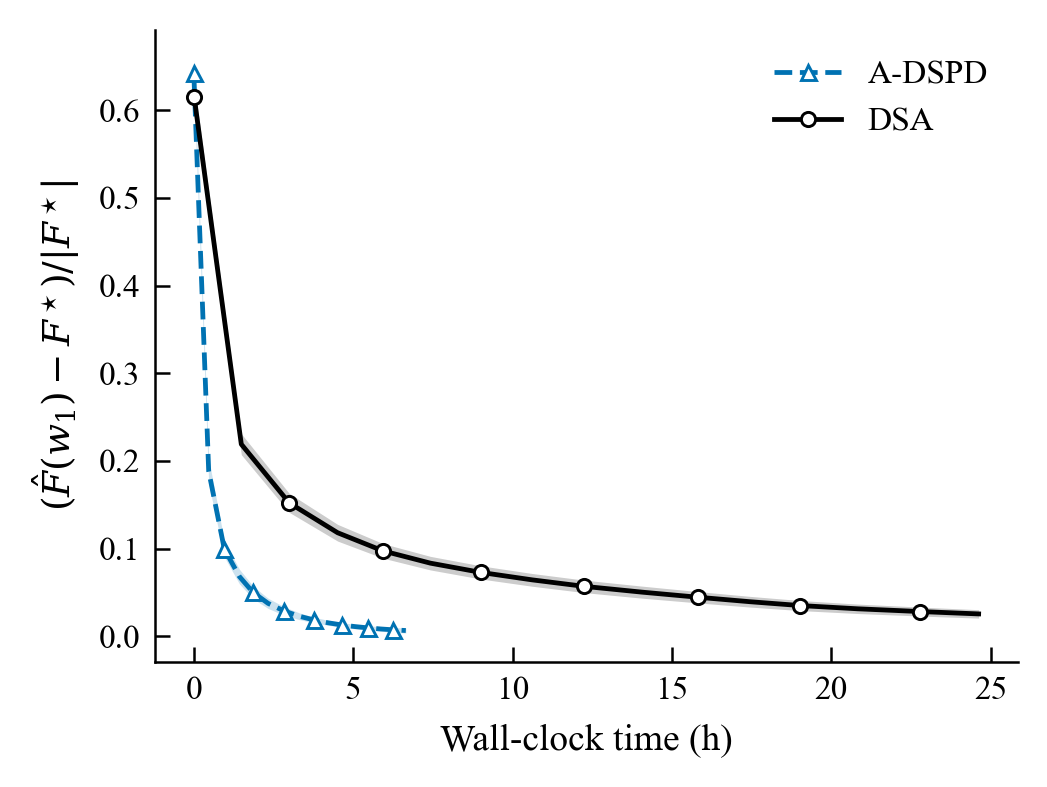}\hfill
  \includegraphics[width=0.32\textwidth]{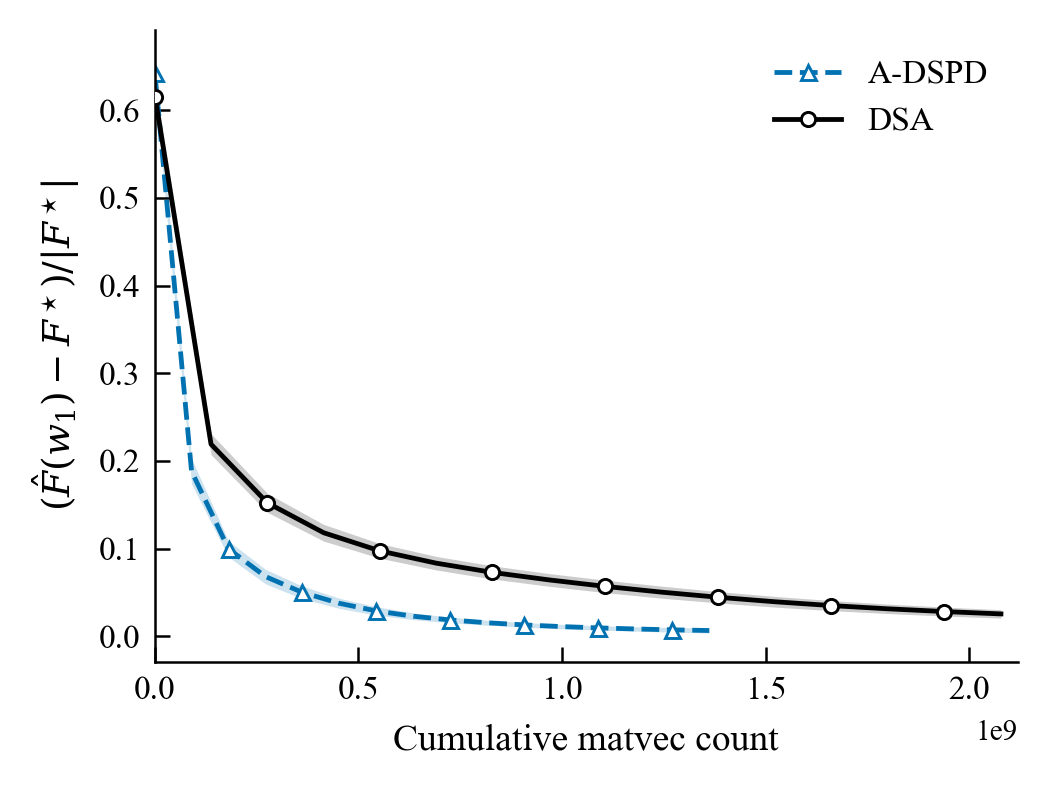}\\[-0.4em]
  {\footnotesize (c) Large instance ($n=200$ assets)}
  \caption{Fixed-$w_1$ relative optimality gap~\eqref{eq:fixed-w1-gap} versus outer iteration (left), wall-clock time (middle), and cumulative matvec count (right) across $n=20$ (top), $n=100$ (middle), and $n=200$ (bottom).}
  \label{fig:mv-relgap}
\end{figure}

\paragraph{Performance.}
Table~\ref{tab:mv-fixed-w1} reports the fixed-$w_1$ relative optimality gap~\eqref{eq:fixed-w1-gap} at the final checkpoint, together with covariance and linking matvec counts and wall-clock time from the dynamic runs.
Across the five seeds, \ADSPD{} attains a substantially smaller relative gap than \DSA{} on all three instances while using roughly one-third as many covariance matvecs and less wall-clock time.
Figure~\ref{fig:mv-relgap} plots the relative gap against the outer-iteration index, elapsed wall-clock time, and cumulative matvec count for $n=20$, $n=100$, and $n=200$.

\begin{table}[H]
  \centering
  \caption{Mean--variance portfolio model performance comparison.}
  \label{tab:mv-fixed-w1}
  \small
  \begin{tabular}{cccccc}
    \toprule
    $n$ & Method & Relative optimality gap (\%) & Covariance matvecs & Linking matvecs & Time (s) \\
    \midrule
    $20$ & \DSA{}   & $9.16\pm1.18$          & $2.300\times10^{7}$ & $1.816\times10^{7}$ & $299.2\pm2.4$ \\
         & \ADSPD{} & $\mathbf{3.83\pm0.28}$ & $\mathbf{6.040\times10^{6}}$ & $1.816\times10^{7}$ & $\mathbf{18.9\pm0.1}$ \\
    \addlinespace
    $100$ & \DSA{}   & $6.33\pm0.43$          & $2.908\times10^{8}$ & $2.906\times10^{8}$ & $8369\pm348$ \\
          & \ADSPD{} & $\mathbf{2.35\pm0.25}$ & $\mathbf{9.664\times10^{7}}$ & $2.906\times10^{8}$ & $\mathbf{1269\pm30}$ \\
    \addlinespace
    $200$ & \DSA{}   & $2.56\pm0.48$          & $1.056\times10^{9}$ & $1.022\times10^{9}$ & $88621\pm33223$ \\
          & \ADSPD{} & $\mathbf{0.67\pm0.29}$ & $\mathbf{3.398\times10^{8}}$ & $1.022\times10^{9}$ & $\mathbf{23951\pm14771}$ \\
    \bottomrule
  \end{tabular}
\end{table}

\paragraph{Inner proximal accuracy.}
To assess the accuracy of the inner proximal solves in \DSA{},
we monitor the NAG stopping diagnostics.
On the $n=20$ instance, the inner NAG subroutine terminates in an average of $3.81\pm0.04$ iterations per proximal solve with a mean relative stopping residual of $(4.52\pm0.02)\times10^{-6}$; on the $n=100$ instance, it terminates in an average of $3.01\pm0.00$ iterations with a mean stopping residual of $(2.37\pm0.20)\times10^{-6}$; and on the $n=200$ instance, the subroutine requires an average of $3.11\pm0.18$ iterations per solve.
Each NAG iteration requires a covariance matvec, whereas \ADSPD{} uses one such product per primal update.

\subsection{Mean--conditional value-at-risk (CVaR) portfolio model}
\label{subsec:exp-cvar}

The second experiment replaces the equity mean--variance stage cost with a mean--CVaR formulation on the Kenneth French \emph{100 Portfolios Formed on Size and Book-to-Market} (daily, value-weighted)~\footnote{Data retrieved from the Kenneth R. French Data Library (\url{https://mba.tuck.dartmouth.edu/pages/faculty/ken.french/data_library.html}).}.

At a node of stage~$t$, let $\{R_{t,\ell}\}_{\ell=1}^{S_t}\subset\mathbb{R}^n$ be the local return sample of size $S_t$, and let $R_t$ follow the uniform empirical distribution $\widehat{\mathbb P}_t=S_t^{-1}\sum_{\ell=1}^{S_t}\delta_{R_{t,\ell}}$, where $\delta_r$ denotes unit probability mass at $r$.
The scalar $\zeta$ is a loss threshold and the components of $u$ bound excess losses above that threshold. Then CVaR at level $\alpha\in(0,1)$ admits the Rockafellar--Uryasev (RU) representation~\citep{rockafellar2000optimization}
\begin{equation}
\label{eq:ru-cvar}
  \mathrm{CVaR}_\alpha\bigl(-R_t^\top w_t\bigr)
  =
  \min_{\zeta\in\bbr,\,u\in\bbr^{S_t}_+}
  \Bigl\{
  \zeta+\frac{1}{(1-\alpha)S_t}\sum_{\ell=1}^{S_t}u_\ell
  :\
  u_\ell\ge -R_{t,\ell}^\top w_t-\zeta
  \Bigr\}.
\end{equation}
From the local return sample, define
\begin{equation}
\begin{aligned}
\hat\mu_t
&=\frac{1}{S_t}\sum_{\ell=1}^{S_t}R_{t,\ell},\\
\widetilde\Sigma_t
&=\frac{1}{S_t-1}\sum_{\ell=1}^{S_t}
(R_{t,\ell}-\hat\mu_t)(R_{t,\ell}-\hat\mu_t)^\top,\\
\hat\Sigma_t
&=0.9\,\widetilde\Sigma_t+0.1\,\operatorname{Diag}(\widetilde\Sigma_t).
\end{aligned}
\end{equation}
Here $\hat\mu_t$ is the sample mean and $\hat\Sigma_t$ is the covariance estimate with $10\%$ diagonal shrinkage.

Let $\eta$ and $\lambda$ be the return and regularization weights. We take the stage cost $\mathrm{CVaR}_\alpha(-R_t^\top w_t)-\eta\hat\mu_t^\top w_t$ with~\eqref{eq:ru-cvar}, and add the quadratic regularizer $\tfrac{\lambda}{2}w_t^\top\hat\Sigma_t w_t$.
Adding a covariance or quadratic regularizer to an empirical CVaR objective is common in multi-stage risk-averse portfolio selection.
This yields
\begin{equation}
\label{eq:cvar-stage-cost}
\begin{aligned}
f_t(x_t,c_t)
&=
\frac{\lambda}{2}\,w_t^\top\hat\Sigma_t w_t
+\zeta_t
+\frac{1}{(1-\alpha)S_t}\sum_{\ell=1}^{S_t} u_{t,\ell}
-\eta\,\hat\mu_t^\top w_t,
\end{aligned}
\end{equation}
Here $c_t$ collects the node's return samples and cost coefficients, and the full stage decisions are $x_1=(w_1,\zeta_1,u_1)$ and $x_t=(w_t,s_t,\zeta_t,u_t)$ for $t=2,3$. The excess-loss vector $u_t$ is a local decision block, distinct from the fixed predecessor decision $u$ in~\eqref{consistent}.
The overall problem is therefore a nested expectation of~\eqref{eq:cvar-stage-cost}.
In the reported experiment, $\eta=1$ weights the expected-return reward and $\lambda=0.1$ weights the quadratic covariance regularizer.
Formulation~\eqref{eq:ru-cvar} leaves $\zeta_t$ free and $u_t$ bounded only from below, whereas the convergence guarantees of both methods require compact decision sets $X_t$.
We therefore add redundant box constraints derived from the return samples.
For a stage-$t$ node $\nu$ with local sample $\{R^{\nu}_{t,\ell}\}_{\ell=1}^{S_t}$, define the nodewise loss extrema
\begin{equation*}
  \underline{\zeta}^{\nu}_t:=\min_{\ell,i}\bigl(-R^{\nu}_{t,\ell,i}\bigr),
  \qquad
  \overline{\zeta}^{\nu}_t:=\max_{\ell,i}\bigl(-R^{\nu}_{t,\ell,i}\bigr),
\end{equation*}
and the stage-level constants
$\underline{\zeta}_t:=\min_{\nu}\underline{\zeta}^{\nu}_t$,
$\overline{\zeta}_t:=\max_{\nu}\overline{\zeta}^{\nu}_t$, and
$\overline{u}_t:=\max_{\nu}\bigl(\overline{\zeta}^{\nu}_t-\underline{\zeta}^{\nu}_t\bigr)$,
where $\nu$ ranges over all stage-$t$ nodes of the tree.
Since $w_t\in\Delta_n$, every nodewise loss $L^{\nu}_{t,\ell}(w_t):=-(R^{\nu}_{t,\ell})^\top w_t$ is a convex combination of the entries of $-R^{\nu}_{t,\ell}$ and hence lies in $[\underline{\zeta}^{\nu}_t,\overline{\zeta}^{\nu}_t]$.
For fixed $w_t$, the minimum in~\eqref{eq:ru-cvar} is attained at any $\alpha$-quantile $\zeta^\star$ of the empirical loss distribution, which belongs to $[\underline{\zeta}^{\nu}_t,\overline{\zeta}^{\nu}_t]$, together with $u^\star_\ell=[L^{\nu}_{t,\ell}(w_t)-\zeta^\star]_+\le\overline{\zeta}^{\nu}_t-\underline{\zeta}^{\nu}_t\le\overline{u}_t$.
Because $(\zeta_t,u_t)$ enter only the stage cost~\eqref{eq:cvar-stage-cost} and not the linking constraints, the constraints
\begin{equation}
\label{eq:ru-box}
  \zeta_t\in[\underline{\zeta}_t,\overline{\zeta}_t],
  \qquad
  u_t\in[0,\overline{u}_t]^{S_t}
\end{equation}
are redundant for the nested problem: they do not change its optimal value or its set of optimal portfolios, while making each $X_t$ compact.
Both methods enforce~\eqref{eq:ru-box} in their projection steps, so $X_1=\Delta_n\times[\underline{\zeta}_1,\overline{\zeta}_1]\times[0,\overline{u}_1]^{S_1}$ and, for $t=2,3$, $X_t=\Delta_n\times\{s_t\ge 0:\mathbf{1}^\top s_t\le\psi_t\}\times[\underline{\zeta}_t,\overline{\zeta}_t]\times[0,\overline{u}_t]^{S_t}$.
The corresponding Euclidean diameter bounds used in the step-size rules are
\begin{equation}
\label{eq:cvar-diameter}
  \Omega_t^2
  =\Omega_{\Delta}^2+(2\psi_t)^2+\bigl(\overline{\zeta}_t-\underline{\zeta}_t\bigr)^2+S_t\,\overline{u}_t^{\,2},
  \qquad \psi_1:=0,
\end{equation}
where $\Omega_\Delta=\sqrt{2}$ upper-bounds the diameter of the simplex constraint $w_t\in\Delta_n$ and $2\psi_t$ upper-bounds the diameter of the turnover block $\{s_t\ge 0:\mathbf{1}^\top s_t\le\psi_t\}$.

\paragraph{Scenario tree and problem instance.}
We construct a scenario tree from historical return blocks with decision date $t_0=\mathtt{2019\textrm{-}12\textrm{-}31}$.
Stage~$1$ uses a $252$-day lookback ending at $t_0$; stages~$2$ and~$3$ use historical blocks of length $H_2=H_3=126$ drawn strictly before $t_0$, with uniform branching probabilities and stage-$3$ children sampled conditionally after the parent stage-$2$ block.
The reported instance has $30$ stage-$2$ nodes and $30$ stage-$3$ children per stage-$2$ node, giving $900$ root-to-leaf paths and $931$ nodes in total.
We use the CVaR confidence level $\alpha=0.95$.
Since no within-node subsampling is applied, $(S_1,S_2,S_3)=(252,126,126)$, and the turnover limits are $\psi_2=\psi_3=0.25$.
Both methods use the budget $(N_1,N_2,N_3)=(1000,1000,150)$.

\paragraph{Metrics and baseline.}
The comparison uses the same fixed-$w_1$ protocol as in Section~\ref{subsec:exp-protocol}: COPT reoptimizes all recourse variables given $\widehat w_{1,s}$, yielding $Q_s(\widehat w_{1,s})$ and $\operatorname{Gap}_s$ relative to $F_s^\star$.
As a baseline, we use the equally weighted (EW) first-stage portfolio
$\hat{\mathbf{w}}_1=\mathbf{1}/n$ and evaluate it using the same
fixed-$w_1$ evaluation setting.
This policy is feasible with zero turnover ($s_2=s_3=0$), and feasible $(\zeta_t,u_t)$ are obtained from~\eqref{eq:ru-cvar}.
EW requires no optimization and therefore requires no matrix--vector products, and we also obtain its $Q_s$ value using COPT.

\begin{table}[H]
  \centering
  \caption{Mean--CVaR portfolio model performance comparison.}
  \label{tab:cvar-30x30}
  \small
  \begin{tabular}{@{}llccccc@{}}
    \toprule
    & Method & Mean--CVaR cost $Q(\widehat w_1)$ & Relative gap (\%) & Covariance matvecs & Linking matvecs & Time (s) \\
    \midrule
    & EW       & $0.05485\pm0.00052$ & $26.37\pm1.38$ & $0$ & $0$ & not applicable \\
    & \DSA{}   & $\mathbf{0.04471\pm0.00065}$ & $2.99\pm0.41$ & $4.575\times10^{8}$ & $4.540\times10^{8}$ & $26631\pm139$ \\
    & \ADSPD{} & $\mathbf{0.04471\pm0.00073}$ & $3.01\pm0.19$ & $\mathbf{1.510\times10^{8}}$ & $4.540\times10^{8}$ & $\mathbf{5147\pm12}$ \\
    \bottomrule
  \end{tabular}
\end{table}

\begin{figure}[H]
  \centering
  \includegraphics[width=0.32\textwidth]{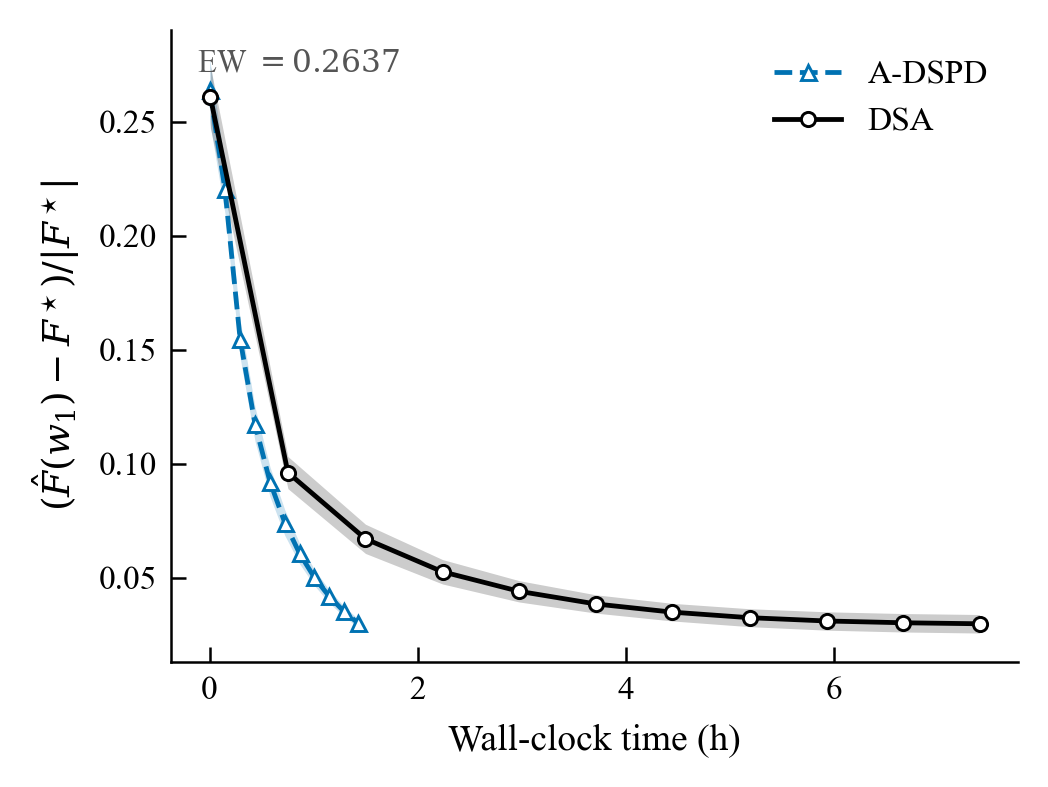}\hfill
  \includegraphics[width=0.32\textwidth]{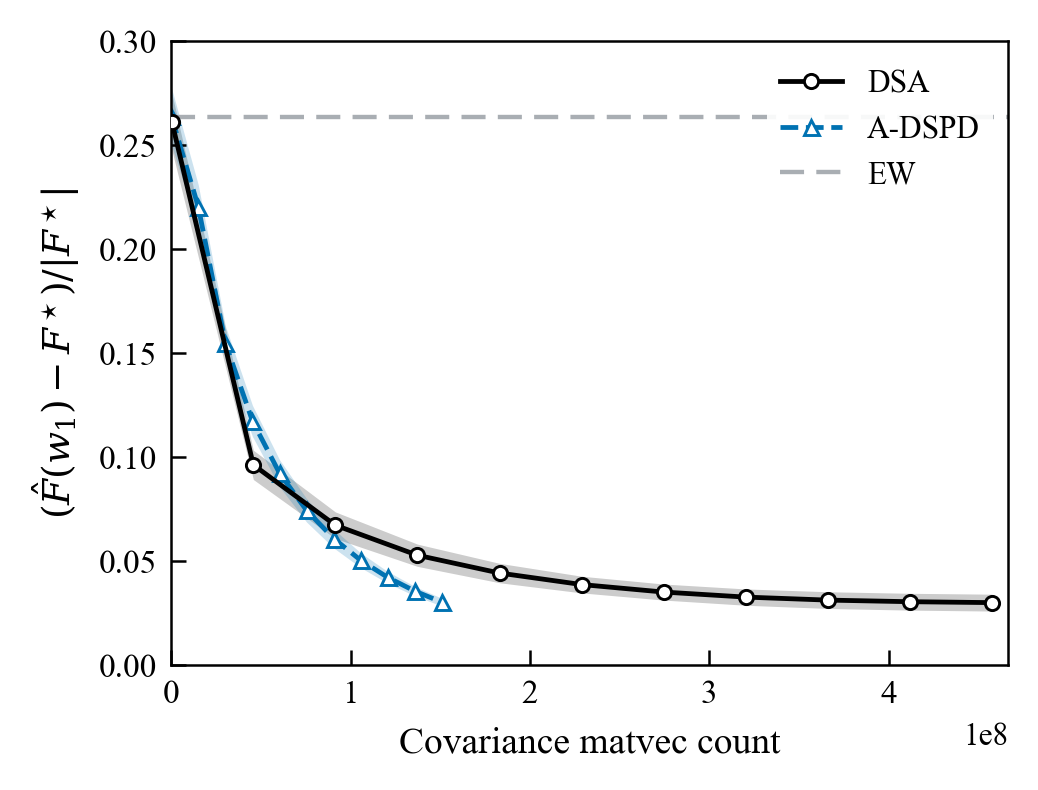}\hfill
  \includegraphics[width=0.32\textwidth]{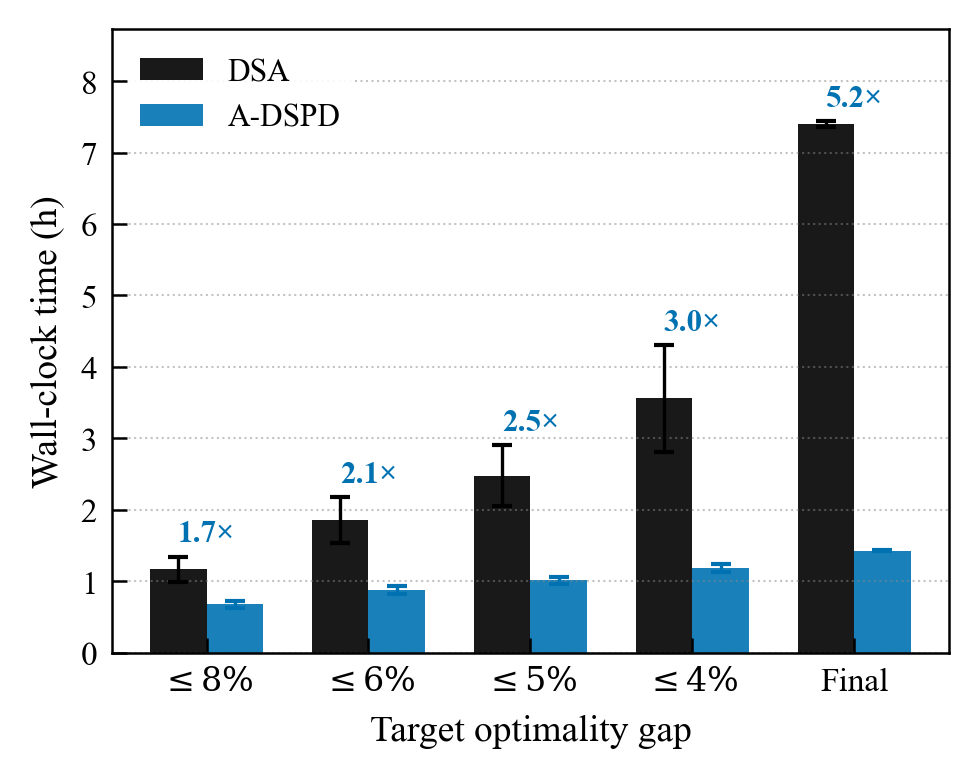}\\[-0.4em]
  \makebox[0.32\textwidth]{\footnotesize (a) Gap vs.\ wall-clock time}\hfill
  \makebox[0.32\textwidth]{\footnotesize (b) Gap vs.\ covariance matvecs}\hfill
  \makebox[0.32\textwidth]{\footnotesize (c) Time-to-accuracy milestones}
  \caption{
  Panels~(a) and~(b) plot the fixed-$w_1$ relative optimality gap~\eqref{eq:fixed-w1-gap} against elapsed wall-clock time and cumulative covariance matvec count.
  Panel~(c) compares the average wall-clock time required to achieve prescribed gap tolerances (error bars are 95\% confidence intervals based on Student's $t$ distribution), with speedup factors of \ADSPD{} over \DSA{} annotated above each pair.}
  \label{fig:cvar-relgap-b30}
  \vspace{-0.4em}
\end{figure}

\paragraph{Performance.}
Table~\ref{tab:cvar-30x30} compares solution quality and computational cost at the common iteration budget. The two methods obtain similar mean relative gaps: $2.99\%\pm0.41\%$ for \DSA{} and $3.01\%\pm0.19\%$ for \ADSPD{}, compared with $26.37\%\pm1.38\%$ for EW. The mean running time decreases from $26631$~s to $5147$~s, a factor of $5.2$. \ADSPD{} uses about one-third as many covariance matvecs.
Figure~\ref{fig:cvar-relgap-b30} compares the time and covariance matvecs needed to reduce the fixed-$w_1$ gap. For each of the five paired seeds, \ADSPD{} reaches all four reported thresholds ($8\%$, $6\%$, $5\%$, and $4\%$) sooner than \DSA{}. 

\section{Conclusions}\label{sec:conclusion}
We have introduced \ADSPD{}, a framework for multi-stage stochastic optimization problems with smooth convex stage costs and linear inequality constraints.
The method retains the recursive value-function structure of \DSA{} and replaces the proximal subproblems of the stage costs by gradient updates obtained from linearization.
The subproblem solver \iaspd{} attains an $O(1/N^2)$ rate for the smooth component and the $O(1/\sqrt{N})$ rate for the stochastic nonsmooth component.
Explicit bounds on the inexactness and on the second moments of the recursively generated stochastic $\epsilon$-subgradients, together with bounds on the dual iterates, give non-asymptotic guarantees on optimality, feasibility, and total complexity for three-stage problems and for general $T$-stage problems.
On two three-stage portfolio selection models, we compare the first-stage solution quality of \ADSPD{} and \DSA{} by fixing the computed first-stage portfolio and reoptimizing recourse. We report convergence against outer iteration, matrix--vector product count, and wall-clock time.
One direction for future work is to reduce the dependence of the total computational cost on the number of stages.
Another is to extend the framework to nonlinear convex constraints.

\bibliographystyle{abbrvnat}
\bibliography{reference}

\appendix

\appendix
\section{Deferred Proofs for Stage-Specific Corollaries}\label{app:proofs-corollaries}

\subsection{Proof of Corollary~\ref{cor:last-stage}}\label{app:proof-cor-last-stage}
\begin{proof}
		Direct substitution verifies that the parameter choices in \eqref{paraco1} satisfy the conditions of Theorem~\ref{thzb}. The required summation identities are
	\[
    \sum_{k=1}^{N}k = \frac{N(N+1)}{2}, \qquad
	\sum_{k=1}^{N}k^2 = \frac{N(N+1)(2N+1)}{6},
	\]
		Moreover,
	\begin{equation*}
	\begin{aligned}
	\sum_{k=1}^{N}\frac{2\gamma_k M^2}{\eta^{-1}_k-L_f\beta^{-1}_k - 2\norm{A}^2\tau_k}
		\leq 2M^2\sum_{k=1}^{N}k^2 = 2M^2\left(\frac{N(N+1)(2N+1)}{6}\right).
	\end{aligned}
	\end{equation*}
		Since $\alpha_k = 1$ and $\gamma_k = k$,
	\begin{equation*}
	\begin{aligned}
	\sum_{k=1}^{N}\frac{2\gamma^2_k M^2}{\alpha_k}
		= 2M^2\sum_{k=1}^{N}k^2 = 2M^2\left(\frac{N(N+1)(2N+1)}{6}\right).
	\end{aligned}
	\end{equation*}
		Substituting these estimates into Theorem~\ref{thzb} yields
	\begin{equation*}
	\begin{aligned}
	\bbe\left[\operatorname{gap}(\tilde{z}^{N+1})\right]
	&\leq \frac{2\Omega^2_X(L_f+N\norm{A})}{(N+1)N}+\frac{\Omega^2_X}{(N+1)N}
    +\frac{\Omega^2_X}{(N+1)N}\\
    &+\frac{\norm{A}\norm{y^*}^2}{N+1}
	+\frac{4M^2(2N+1)}{3}+\epsilon\\
	&\leq \frac{2\Omega^2_XL_f}{N^2}+\frac{2\Omega_X^2\norm{A}}{N}+\frac{2\Omega^2_X}{N^2}+\frac{\norm{A}\norm{y^*}^2}{N}
		+\frac{4M^2(2N+1)}{3}+\epsilon,
	\end{aligned}
	\end{equation*}
	and
	\begin{equation*}
	\begin{aligned}
	\bbe\left[\norm{\tilde{y}^{N+1}}^2\right]
	&\leq \frac{8\Omega_X^2(L_f+N\norm{A})}{N\norm{A}}+\frac{4\Omega_X^2}{\norm{A}N}
	+\frac{4\Omega_X^2}{\norm{A}N}+6\norm{y^*}^2\\
	&+\frac{16M^2(N+1)(2N+1)}{3\norm{A}}+\frac{4(N+1)\epsilon}{\norm{A}}\\
		&\leq 8\Omega_X^2+\frac{8\Omega_X^2L_f}{N\norm{A}}
	+\frac{8\Omega_X^2}{N\norm{A}}+6\norm{y^*}^2\\
	    &+\frac{16M^2(N+1)(2N+1)}{3\norm{A}}+\frac{4(N+1)\epsilon}{\norm{A}}.
	\end{aligned}
	\end{equation*}
		Applying Theorem~\ref{optgapth} with the same estimates yields the bounds in Corollary~\ref{cor:last-stage}(c).
\end{proof}

\subsection{Proof of Corollary~\ref{cor:middle-stage}}\label{app:proof-cor-middle-stage}
\begin{proof}
		Direct substitution verifies that the parameter choices in \eqref{paraco2} satisfy the conditions of Theorem~\ref{thzb}. Applying Theorem~\ref{thzb} together with the estimates
	\[
	\sum_{k=1}^{N} \sqrt{k} < \left(\frac{2N}{3} + \frac{1}{2}\right)\sqrt{N} \leq \frac{2\sqrt{N}(N+1)}{3}, \qquad
	\sum_{k=1}^{N} k^2 = \frac{N(N+1)(2N+1)}{6},
	\]
	as well as
	\begin{equation*}
	\begin{aligned}
	\sum_{k=1}^{N}\frac{2\gamma_k M^2}{\eta^{-1}_k-L_f\beta^{-1}_k-2\norm{A}^2\tau_{k}}\leq 2M\sum_{k=1}^{N}\sqrt{k}\leq\frac{4M(N+1)\sqrt{N}}{3},
	\end{aligned}
	\end{equation*}
		yields the following bounds:
	\begin{equation*}
	\begin{aligned}
	\bbe\left[\operatorname{gap}(\tilde{z}^{N+1})\right]
	&\leq \frac{2\Omega^2_X(L_f+N\norm{A})}{(N+1)N}+\frac{\Omega^2_XM\sqrt{N}}{N+1}+\frac{\Omega_X^2M\sqrt{N}}{N+1}\\
    &+\frac{\sqrt{N}\norm{A}\norm{y^*}^2}{N+1}
	   +\frac{8M}{3\sqrt{N}}+\frac{2M(2N+1)}{3N^{\frac{3}{2}}}+\epsilon\\
	&\leq \frac{2\Omega^2_XL_f}{N^2}+\frac{2\Omega_X^2\norm{A}}{N}+\frac{2M\Omega^2_X}{\sqrt{N}}+\frac{\norm{A}\norm{y^*}^2}{\sqrt{N}}\\
	&+\frac{8M}{3\sqrt{N}}+\frac{4M}{3\sqrt{N}}+\frac{2M}{3N^{\frac{3}{2}}}+\epsilon\\
		&= \frac{2\Omega^2_XL_f}{N^2}+\frac{2\Omega_X^2\norm{A}}{N}+\frac{2M\Omega^2_X}{\sqrt{N}}+\frac{\norm{A}\norm{y^*}^2}{\sqrt{N}}
		+\frac{12M}{3\sqrt{N}}+\frac{2M}{3N^{\frac{3}{2}}}+\epsilon,
	\end{aligned}
	\end{equation*}
	and	
	\begin{equation*}
	\begin{aligned}
	\bbe\left[\norm{\tilde{y}^{N+1}}^2\right]
	&\leq \frac{8\Omega_X^2(L_f+N\norm{A})}{N^{\frac{3}{2}}\norm{A}}+\frac{4M\Omega_X^2}{\norm{A}}
	+\frac{4M\Omega_X^2}{\norm{A}}+6\norm{y^*}^2\\
	&\quad + \frac{32M(N+1)}{3\norm{A}N}+\frac{8M(N+1)(2N+1)}{3\norm{A}N^2}+\frac{4(N+1)\epsilon}{\sqrt{N}\norm{A}}\\
	&\leq \frac{8L_f\Omega_X^2}{N^{\frac{3}{2}}\norm{A}}+\frac{8\Omega_X^2}{\sqrt{N}}+\frac{8M\Omega_X^2}{\norm{A}}+6\norm{y^*}^2
	+\frac{64M}{3\norm{A}}+\frac{16M}{\norm{A}}+\frac{8\sqrt{N}\epsilon}{\norm{A}},
	\end{aligned}
	\end{equation*}
	where the second inequality follows from the facts that $\frac{N+1}{N}\leq 2$, $\frac{(N+1)(2N+1)}{N^2}\leq 6$, and $\frac{N+1}{\sqrt{N}}\leq 2\sqrt{N}$.
		Applying Theorem~\ref{optgapth} with the same estimates yields the bounds in Corollary~\ref{cor:middle-stage}(c).
\end{proof}

\subsection{Proof of Corollary~\ref{cor:first-stage}}\label{app:proof-cor-first-stage}
\begin{proof}
		Direct substitution verifies that the parameter choices in \eqref{paraco3} satisfy the conditions of Theorem~\ref{thzb}. Applying Theorem~\ref{thzb} and using the estimates
		\begin{equation*}
		\sum_{k=1}^{N} \sqrt{k}<\left(\frac{2 N}{3}+\frac{1}{2}\right) \sqrt{N}\leq\left(\frac{2 (N+1)}{3}\right) \sqrt{N}, \quad  \quad
		\sum_{k=1}^{N}k^2=\frac{N(N+1)(2N+1)}{6},
		\end{equation*}
		together with the estimate
		\begin{equation*}
	\begin{aligned}
	\sum_{k=1}^{N}\frac{2\gamma_k M^2}{\eta^{-1}_k-L_f\beta^{-1}_k-2\norm{A}^2\tau_{k}}\leq 2M\sum_{k=1}^{N}\sqrt{k}\leq\frac{4M(N+1)\sqrt{N}}{3},
	\end{aligned}
	\end{equation*}
		yields
	\begin{equation*}
	\begin{aligned}
	\bbe\left[\operatorname{gap}(\tilde{z}^{N+1})\right]
	&\leq \frac{2\Omega^2_X(L_f+N\norm{A})}{(N+1)N}+\frac{\Omega^2_XM\sqrt{N}}{N+1}+\frac{\Omega_X^2M\sqrt{N}}{N+1}\\
    &+\frac{\norm{A}\norm{y^*}^2}{N+1}
	 +\frac{8M}{3\sqrt{N}}+\frac{2M(2N+1)}{3N^{\frac{3}{2}}}+\epsilon\\
	&\leq \frac{2\Omega^2_XL_f}{N^2}+\frac{2\Omega^2_X\norm{A}}{N}+\frac{2M\Omega^2_X}{\sqrt{N}}+\frac{\norm{A}\norm{y^*}^2}{N}\\
	&+\frac{8M}{3\sqrt{N}}+\frac{4M}{3\sqrt{N}}+\frac{2M}{3N^{\frac{3}{2}}}+\epsilon\\
		&= \frac{2\Omega^2_XL_f}{N^2}+\frac{2\Omega^2_X\norm{A}}{N}+\frac{2M\Omega^2_X}{\sqrt{N}}+\frac{\norm{A}\norm{y^*}^2}{N}
	 +\frac{12M}{3\sqrt{N}}+\frac{2M}{3N^{\frac{3}{2}}}+\epsilon,
	\end{aligned}
	\end{equation*}
	and	
	\begin{equation*}
	\begin{aligned}
	\bbe\left[\norm{\tilde{y}^{N+1}}^2\right]
	&\leq \frac{8\Omega^2_X(L_f+N\norm{A})}{N\norm{A}}+\frac{4\Omega_X^2M\sqrt{N}}{\norm{A}}
	+\frac{4\Omega_X^2M\sqrt{N}}{\norm{A}}+6\norm{y^*}^2\\
	&\quad + \frac{32M(N+1)}{3\sqrt{N}\norm{A}}
	+\frac{8M(N+1)(2N+1)}{3N^{\frac{3}{2}}\norm{A}}+\frac{4(N+1)\epsilon}{\norm{A}}\\
	&\leq 8\Omega^2_X+\frac{8\Omega^2_X L_f}{N\norm{A}}+\frac{8\Omega_X^2M\sqrt{N}}{\norm{A}}+6\norm{y^*}^2\\
	&\quad + \frac{32M(N+1)}{3\sqrt{N}\norm{A}}
	+\frac{8M(N+1)(2N+1)}{3N^{\frac{3}{2}}\norm{A}}+\frac{4(N+1)\epsilon}{\norm{A}}\\
	\end{aligned}
	\end{equation*}
		Applying Theorem~\ref{optgapth} with the preceding estimates yields the bounds in Corollary~\ref{cor:first-stage}(c).
\end{proof}

\end{document}